\documentclass[11pt]{article}

\usepackage[authoryear,round]{natbib}
\usepackage{amsmath,amssymb,amsthm,amsfonts}
\usepackage{array,siunitx,graphicx}
\usepackage{authblk}
\usepackage{enumerate}
\usepackage{placeins}
\usepackage{caption}
\usepackage{titlesec}
\usepackage{newtxtext}
\usepackage[subscriptcorrection]{newtxmath}
\usepackage[margin=1in]{geometry}
\usepackage[colorlinks=true,linkcolor=blue,citecolor=blue,urlcolor=blue]{hyperref}
\hypersetup{
  pdftitle={Signal-to-Noise Ratio Inference for Multivariate High-Dimensional Linear Models},
  pdfauthor={Xiaohan Hu, Zhentao Li, and Xiaodong Li}
}

\titleformat{\section}[hang]
  {\normalfont\bfseries\fontsize{11}{13pt plus .8pt minus .6pt}\selectfont}
  {\thesection.}{1em}{}
\titleformat{\subsection}[hang]
  {\normalfont\bfseries\fontsize{11}{13pt plus .8pt minus .6pt}\selectfont}
  {\thesubsection}{1em}{}

\allowdisplaybreaks[3]

\theoremstyle{plain}
\newtheorem{theorem}{Theorem}
\newtheorem{lemma}{Lemma}

\theoremstyle{definition}

\theoremstyle{remark}
\newtheorem*{remark}{Remark}

\newcommand{\vct}[1]{#1}
\newcommand{\mtx}[1]{#1}
\newcommand{\tr}{\operatorname{tr}}
\newcommand{\Var}{\operatorname{var}}
\newcommand{\E}{\operatorname{E}}
\newcommand{\Cov}{\operatorname{cov}}
\newcommand{\diag}{\operatorname{diag}}
\newcommand{\PP}{\operatorname{pr}}
\DeclareMathOperator{\vecop}{vec}

\newif\ifmainpaper
\newif\iftechnicalappendix
\mainpapertrue
\technicalappendixtrue

\title{Signal-to-Noise Ratio Inference for Multivariate\\
High-Dimensional Linear Models}
\author[1]{Xiaohan Hu}
\author[2]{Zhentao Li}
\author[2]{Xiaodong Li}
\affil[1]{Department of Mathematics, University of California, Davis}
\affil[2]{Department of Statistics, University of California, Davis}
\date{}

\begin{document}
\maketitle

\ifmainpaper

\setlength{\abovedisplayskip}{7.5pt plus 2pt minus 2pt}
\setlength{\belowdisplayskip}{7.5pt plus 2pt minus 2pt}
\setlength{\abovedisplayshortskip}{2pt plus 2pt}
\setlength{\belowdisplayshortskip}{5pt plus 2pt minus 2pt}

\begin{abstract}
We study inference for the fraction of total response variation explained by
multivariate high-dimensional linear models. Two moment equations yield
explicit estimates without fitting the regression coefficients. We establish
asymptotically valid confidence intervals for fixed- and random-effects
models, allowing the response dimension to grow with the sample size. The
theory describes how response dependence affects precision and provides a
correction for unequal noise levels in the random-design setting.
Simulations assess finite-sample performance, and an analysis of yeast growth
across environments illustrates joint inference for additive marker-based
signal.
\end{abstract}


\section{Introduction}

How much of the variation in a response can a linear model explain? In
genetics, the corresponding question concerns heritability. Estimates of the
noise level are also used to choose tuning parameters in regularized
regression. These questions remain meaningful when there are as many
predictors as observations, or more. In that case, estimating all the
regression coefficients may be difficult or impossible, but the total signal
and noise can still be estimated. We study this problem when each observation
has several responses that may be correlated.

Variance-component analysis provides a broader framework for separating
sources of variation in mixed models. \citet{jiang1996reml} established
conditions for consistency and asymptotic normality of restricted maximum
likelihood (REML) variance-component estimators, including for non-Gaussian
data. \citet{jiang2005partially} developed consistent covariance estimation
for inference in non-Gaussian mixed linear models. Moment methods also
provide scalable estimation and uncertainty assessment for variance
components in large, unbalanced crossed random-effects models
\citep{gao2017efficient,gao2020estimation}. In those models, observations
are linked by shared factor levels, such as customers and items.

For a single response, the literature provides several approaches to
estimating signal and noise in high dimensions. These include moment
estimators \citep{dicker2014variance}, confidence intervals based on the
singular values of the predictor matrix \citep{janson2017eigenprism}, and
fixed-effects inference under different assumptions on the design and the
coefficients
\citep{verzelen2018adaptive,cai2020semisupervised,song2024debiased}.
Related work studies likelihood-based variance estimation in high dimensions
\citep{dicker2016maximum}, REML under misspecified genetic mixed models
\citep{jiang2016unified}, and corresponding confidence intervals for error
variance and heritability \citep{dao2022variance}. Inference with dense
coefficients and unequal noise variances is studied by \citet{hu2022optimal}.

Several correlated responses occur in neuroimaging, molecular studies across
tissues, and breeding studies across environments
\citep{zhao2022bayesian,flutre2013statistical,burgueno2012genomic}.
Multivariate methods estimate genetic and residual covariance matrices
\citep{zhou2014efficient}, heritability for a linear combination of responses
chosen in advance \citep{li2022heritability}, or an overall heritability based
on covariance traces \citep{ge2016estimating}. High-dimensional MANOVA theory
studies the eigenvalues of covariance estimators \citep{fan2019eigenvalue}.
The noise structure also matters in applications: differences in residual
variances across environments and production groups have been documented in
genomic prediction \citep{ou2016genomic}.

This paper studies one overall measure of explained variation, together with
its standard error and confidence interval. The target is the total signal
variance across responses divided by the total signal and noise variance,
which we call the signal-to-noise ratio (SNR). It is a variance-weighted
average of the separate explained fractions and is unchanged by orthogonal
rotations of the responses. Under random-effects models, this target is a
function of signal and noise variance components, and its standard error
depends on their joint uncertainty. Two moment equations give an estimate
without fitting the full coefficient matrix. Our contribution concerns
inference for this target in three settings, each connected to existing work.

For fixed coefficients, the starting point is \citet{dicker2014variance}, who
established consistency and asymptotic normality of moment estimators for one
response without sparsity assumptions. We retain his Gaussian random-design
setting with known population predictor covariance; our estimator agrees
with his in the single-response case. We extend inference to several
correlated responses, allowing their number to increase with the sample
size, and prove consistency of the estimated standard error. The signal in
this setting is defined through prediction for new observations.

For random coefficients, the closest multivariate work is
\citet{ge2016estimating}. They developed trace-based moment estimation,
sampling-variance formulas, and Wald and permutation inference, and showed
how combining traits can improve precision. After matching normalizations,
their point estimator agrees with ours under common noise covariance.
Our contribution in Section~3 is a conditional asymptotic theory for this
estimator with the predictors held fixed and the number of responses $q$
allowed to grow with $n$. We prove a studentized normal limit, relative
consistency of the estimated variance, and asymptotic confidence-interval
coverage under explicit design and response-covariance conditions. These
conditions also allow some response-covariance eigenvalues to grow.

Our third setting allows random predictors and coefficients with noise
covariance that varies across observations. Related single-response inference
under unequal noise variances was studied by \citet{hu2022optimal}.
We quantify the additional uncertainty in the multivariate model and construct
adjusted confidence intervals when the noise covariances differ only by scale.

Across the three settings, the variance formulas show how the number of
responses and their correlations affect precision. We give conditions under
which the standard error decreases at the inverse-square-root rate as
responses are added, and show how stronger shared variation limits this
gain. The confidence intervals apply to explained fractions bounded away
from zero and one.

\subsection{Organization of the paper}

Section~2 studies fixed coefficients with random predictors. Section~3 studies
random coefficients with the predictor matrix held fixed. Section~4 allows
both random predictors and noise covariances that differ between observations.
Section~5 reports simulations for these three settings. Section~6 applies the
fixed-design method to yeast growth across environments, and Section~7
discusses the findings and limitations. The Supplementary Material contains
proofs, supporting calculations, and details of the yeast analysis.

\section{Fixed-effects inference}
\label{sec:fixed_mt}

\subsection{Model and estimation}
\label{sec:fixed_model}
\label{sec:fixed_estimator}

Write $\mtx Y=\mtx X\mtx B+\mtx E$, where $\mtx X$ records $p$
predictors for $n$ observations and $\mtx Y$ records $q$ responses for each
observation. Here the coefficient matrix $\mtx B$ is fixed. For a new
observation, the predictor vector follows
$\vct x\sim\mathcal N(\vct0,\mtx\Sigma_x)$, with known positive
definite covariance $\mtx\Sigma_x$. The response is
$\vct y=\mtx B^\top\vct x+\vct e$, where
$\vct e\sim\mathcal N(\vct0,\mtx\Sigma_e)$ is independent of
$\vct x$. The signal $\mtx B^\top\vct x$ has covariance
$\mtx\Sigma_b=\mtx B^\top\mtx\Sigma_x\mtx B$. We define
\[
\rho^2=\frac1q\tr(\mtx\Sigma_b),\qquad
\sigma^2=\frac1q\tr(\mtx\Sigma_e),\qquad
r^2=\frac{\rho^2}{\rho^2+\sigma^2},
\]
where $\rho^2$ and $\sigma^2$ are the average signal and noise variances
across the $q$ responses. The fraction $r^2$ measures how much of the total
variation is explained by the linear signal. It combines responses by adding
their variances, rather than by averaging their separate $R^2$ values.
Responses with greater total variance therefore receive greater weight. An
orthogonal rotation of the response coordinates leaves $r^2$ unchanged, but
changing their relative measurement scales need not do so.

For $q=1$, this is the fixed-effects target studied by
\citet{dicker2014variance}, \citet{janson2017eigenprism}, and
\citet{verzelen2018adaptive}. With several responses, it summarizes the
signal without choosing one response or one linear combination, unlike
inference for a combination chosen in advance \citep{li2022heritability}. The single-response target is also closely
related to the explained variance in \citet{cai2020semisupervised}.
Knowledge of the population covariance $\mtx\Sigma_x$ is required for both
our estimator and its inference theory. \citet{dicker2014variance} also
studied unknown covariance, using a sufficiently accurate covariance estimate
or additional conditions relating the covariance to the coefficients.
Those extensions are outside the present theorem; replacing
$\mtx\Sigma_x$ by an estimate would require a separate error analysis.

The estimator uses two summaries of the data. The first is the usual response
cross-product, which measures total variation. The second,
$\mtx Y^\top\mtx X\mtx\Sigma_x^{-1}\mtx X^\top\mtx Y$, weights
products of responses by the similarity of their predictors, accounting for
predictor covariance through $\mtx\Sigma_x^{-1}$. These summaries have
different expected contributions from signal and noise:
\[
\E\!\left(\frac{\mtx Y^\top\mtx Y}{n}\right)
=\mtx\Sigma_b+\mtx\Sigma_e,\qquad
\E\!\left(\frac{\mtx Y^\top\mtx X\mtx\Sigma_x^{-1}\mtx X^\top\mtx Y}{n^2}\right)
=\frac{p+n+1}{n}\mtx\Sigma_b+\frac pn\mtx\Sigma_e.
\]
Solving these equations gives
\begin{equation}
\label{eq:fixed_estimators_mt}
\widehat{\mtx\Sigma}_b
=\frac{\mtx Y^\top(\mtx X\mtx\Sigma_x^{-1}\mtx X^\top-p\mtx I_n)\mtx Y}{n(n+1)},\qquad
\widehat{\mtx\Sigma}_e
=\frac{(n+p+1)\mtx Y^\top\mtx Y-
\mtx Y^\top\mtx X\mtx\Sigma_x^{-1}\mtx X^\top\mtx Y}{n(n+1)}.
\end{equation}
To estimate the explained fraction, we need only the sums of their diagonal
entries, that is, their traces:
\begin{equation}
\label{eq:fixed_scalar_estimators_mt}
\hat\rho^2
=\frac{\tr\{\mtx Y^\top(\mtx X\mtx\Sigma_x^{-1}\mtx X^\top-p\mtx I_n)\mtx Y\}}{qn(n+1)},\qquad
\hat\sigma^2=\frac{\tr(\mtx Y^\top\mtx Y)}{nq}-\hat\rho^2,\qquad
\hat r^2=\frac{\hat\rho^2}{\hat\rho^2+\hat\sigma^2}.
\end{equation}
Thus the point estimate requires only two scalar summaries, not estimates of
the $pq$ regression coefficients. The full matrix estimates in
\eqref{eq:fixed_estimators_mt} are needed for the standard error, which
must account for correlations between responses. For one response, these are
exactly Dicker's two-moment estimators in the known-covariance case
\citep{dicker2014variance}. The transformation
$\mtx X\mtx\Sigma_x^{-1/2}$ gives identity predictor covariance, with
coefficient matrix $\mtx\Sigma_x^{1/2}\mtx B$ and unchanged signal
covariance $\mtx\Sigma_b$. This is the same whitening step used in his
single-response construction.

\subsection{Standard errors and confidence intervals}
\label{sec:fixed_asymptotics}

For a matrix $A$, $\|A\|$ and $\|A\|_F$ denote its operator and Frobenius
norms; $\|\cdot\|_2$ and $\|\cdot\|_{\psi_2}$ denote the Euclidean and
sub-Gaussian norms. Also, $A\preceq B$ means that $B-A$ is positive
semidefinite. For positive sequences, $a_n\asymp b_n$ means that their ratio
is bounded above and away from zero; $\to_P$ and $\Longrightarrow$ denote
convergence in probability and in distribution; and $z_u$ is the $u$th
standard-normal quantile. All limits are as $n\to\infty$. Constants are
uniform over $(n,p,q)$, and unspecified matrix-valued stochastic orders use
the operator norm.

The following function gives the variance needed for inference. For symmetric
$q\times q$ matrices $A,E$ and scalars $s,t$ with $s+t\ne0$, define
\begin{align}
\label{eq:fe_scalar_variance_main}
\mathcal V_{\mathrm{fe},n}(A,E;s,t)
&=\frac{2n^2}{q^2(n+1)^2(s+t)^4}\Big[
\|A\|_F^2\big\{(s+2t)^2+\tfrac pn(s+t)^2\big\}
+2\tr(EA)\big\{t^2+\tfrac pn(s+t)^2\big\}\notag\\
&\hspace{38mm}+\tr(E^2)\big\{s^2+\tfrac pn(s+t)^2\big\}\Big].
\end{align}
The true covariances and variances define the leading variance of
$\sqrt n(\hat r^2-r^2)$. Replacing them by their estimates gives its
data-based estimate:
\[
\sigma_{r,\mathrm{fe}}^2
=\mathcal V_{\mathrm{fe},n}(\mtx\Sigma_b,\mtx\Sigma_e;\rho^2,\sigma^2),
\qquad
\hat\sigma_{r,\mathrm{fe}}^2
=\mathcal V_{\mathrm{fe},n}(\widehat{\mtx\Sigma}_b,\widehat{\mtx\Sigma}_e;\hat\rho^2,\hat\sigma^2).
\]

\begin{theorem}[Fixed-effects inference]
\label{thm:asymptotics_fixed_mt}
Let $p/n=O(1)$ and $q=o(n)$. Suppose that the rows of
$\mtx X\in\mathbb R^{n\times p}$ are i.i.d.
$\mathcal N(\vct0,\mtx\Sigma_x)$, where $\mtx\Sigma_x$ is deterministic,
known and positive definite, that $\mtx B\in\mathbb R^{p\times q}$ is
deterministic, and that, independently of $\mtx X$,
$\vecop(\mtx E)\sim
\mathcal N(\vct0,\mtx\Sigma_e\otimes\mtx I_n)$.
Set $\mtx\Sigma_b=\mtx B^\top\mtx\Sigma_x\mtx B$ and suppose
$\rho^2,\sigma^2\in[c,C]$ and
\[
\frac qn\{1\vee\|\mtx\Sigma_b\|\vee\|\mtx\Sigma_e\|\}\to0.
\]
Then
\[
\frac{\sqrt n(\hat r^2-r^2)}{\hat\sigma_{r,\mathrm{fe}}}
\Longrightarrow\mathcal N(0,1).
\]
\end{theorem}

The estimated standard error of $\hat r^2$ is
$\hat\sigma_{r,\mathrm{fe}}/\sqrt n$. A normal-based, or Wald, interval
is centred at $\hat r^2$ and extends $z_{1-\alpha/2}$ times this standard
error on either side. For every fixed $\alpha\in(0,1)$, its probability of containing
$r^2$ approaches $1-\alpha$ as the sample size increases.

\begin{remark}
\citet{dicker2014variance} already established single-response asymptotic
normality. Theorem~\ref{thm:asymptotics_fixed_mt} extends inference to
several correlated responses and proves consistency of the estimated standard
error while their number may increase with the sample size. The covariance terms in
\eqref{eq:fe_scalar_variance_main} show how much precision is gained.
When all three terms grow in proportion to $q$, the standard error has order
$(nq)^{-1/2}$. Strong shared variation between responses can make these
terms larger, reducing the gain. The result concerns the overall explained
fraction, not a linear combination selected in advance
\citep{li2022heritability}.
\end{remark}

The Gaussian assumption justifies the moment equations and normal
approximation used here. The estimator can be computed for other designs,
but the theorem does not cover them. The interval is intended for $r^2$
bounded away from zero and one. Supplementary
Section~\ref{sec:supp_fixed_theorem} gives the moment calculations, normal
approximation, and proof that the estimated standard error is consistent.

\section{Fixed-design random-effects inference}
\label{sec:random_fixed_homo_mt}

\subsection{Model and estimation}
\label{sec:random_fixed_model}
\label{sec:random_fixed_estimator}

We now hold the predictor matrix $\mtx X$ fixed and treat the coefficients
as random. The rows of $\mtx B$ are independent
$\mathcal N(\vct0,\mtx\Sigma_\beta/p)$ vectors, and the noise rows are
independent $\mathcal N(\vct0,\mtx\Sigma_e)$ vectors. Repeated samples
redraw the coefficients and noise, but use the same predictors. Averaging the
signal covariance across the observed rows gives
\[
\mtx\Sigma_b=\frac{\|\mtx X\|_F^2}{np}\mtx\Sigma_\beta,
\qquad
\rho^2=\frac1q\tr(\mtx\Sigma_b),\qquad
\sigma^2=\frac1q\tr(\mtx\Sigma_e),\qquad
r^2=\frac{\rho^2}{\rho^2+\sigma^2}.
\]
The signal fraction now refers to the supplied predictor matrix. Its overall
scale determines how much variation the random coefficients generate. For
example, multiplying all predictors by the same constant, while keeping the
coefficient distribution fixed, changes the variance of $\mtx X\mtx B$
and therefore changes the signal fraction. The definition above accounts for
this change; the predictors need not be normalized to a prescribed overall
scale.

This is the fixed-design version of variance-component models based on
predictor similarity, such as GCTA \citep{yang2011gcta}. It also underlies
the moment estimation and uncertainty assessment for multivariate heritability
in \citet{ge2016estimating}. This section establishes conditional asymptotic
guarantees for their moment estimator, after matching normalizations,
including when the number of responses increases with the sample size.

As before, we compare the overall response covariance with a response
cross-product weighted by predictor similarity. Equating these two summaries
to their expectations with $\mtx X$ held fixed gives
\begin{equation}
\label{eq:homo_estimators_mt}
\widehat{\mtx\Sigma}_b
=\frac{
\mtx Y^\top\mtx X\mtx X^\top\mtx Y/\|\mtx X\|_F^2
-\mtx Y^\top\mtx Y/n}
{n\tr\{(\mtx X\mtx X^\top)^2\}/\|\mtx X\|_F^4-1},
\qquad
\widehat{\mtx\Sigma}_e=\frac{\mtx Y^\top\mtx Y}{n}-\widehat{\mtx\Sigma}_b.
\end{equation}
Consequently,
\begin{equation}
\label{eq:fixed_design_scalar_estimators_mt}
\hat\rho^2
=\frac{q^{-1}\left\{
\tr(\mtx Y^\top\mtx X\mtx X^\top\mtx Y)/\|\mtx X\|_F^2
-\tr(\mtx Y^\top\mtx Y)/n\right\}}
{n\tr\{(\mtx X\mtx X^\top)^2\}/\|\mtx X\|_F^4-1},
\qquad
\hat\sigma^2=\frac{\tr(\mtx Y^\top\mtx Y)}{nq}-\hat\rho^2,
\qquad
\hat r^2=\frac{\hat\rho^2}{\hat\rho^2+\hat\sigma^2}.
\end{equation}
The denominator equals
\[
\frac1n\left\|
\frac{n\mtx X\mtx X^\top}{\|\mtx X\|_F^2}-\mtx I_n
\right\|_F^2.
\]
This quantity measures how far the normalized predictor-similarity matrix is
from the identity. If it is small, the two summaries are nearly the same and
cannot separate signal from noise accurately. The theorem below requires it
to stay bounded away from zero. Computing the estimator requires neither a
likelihood fit nor the inverse of a $p\times p$ matrix, even when $p$ is
as large as or larger than $n$.

\subsection{Standard errors and confidence intervals}
\label{sec:random_fixed_asymptotics}

The standard error uses three summaries of the predictor-similarity matrix:
\begin{equation}
\label{eq:fixed_design_gk_main}
g_k=\frac1n\tr\!\left\{\left(
\frac{n\mtx X\mtx X^\top}{\|\mtx X\|_F^2}
\right)^k\right\},\qquad k=2,3,4.
\end{equation}
These numbers are unchanged if all predictors are multiplied by the same
constant. They describe how the variation in the predictor matrix is
distributed across its directions. In particular,
$g_2-1=n^{-1}\|n\mtx X\mtx X^\top/\|\mtx X\|_F^2-\mtx I_n\|_F^2$
is the denominator that separates the two moments. The higher moments $g_3$
and $g_4$ also enter the standard error.

For symmetric $q\times q$ matrices $A,E$, positive scalars $s,t$, and
$u_2>1,u_3,u_4\ge0$, define
\begingroup\small
\begin{equation}
\label{eq:fixed_homo_variance_function}
\begin{aligned}
\mathcal V_{\mathrm{re,h}}(A,E;s,t;u_2,u_3,u_4)
&=\frac{2}{q^2(s+t)^4}\Bigg[\|A\|_F^2\Bigg\{u_2t^2+2\frac{u_3-u_2^2}{u_2-1}t(s+t)+\frac{u_4-2u_2u_3+u_2^3}{(u_2-1)^2}(s+t)^2\Bigg\}\\
&\qquad+2\tr(EA)\Bigg\{t^2+\frac{u_3-u_2^2}{(u_2-1)^2}(s+t)^2\Bigg\}+\|E\|_F^2\Bigg\{s^2+\frac{(s+t)^2}{u_2-1}\Bigg\}\Bigg].
\end{aligned}
\end{equation}
\endgroup
With the predictor matrix held fixed, the leading variance of
$\sqrt n(\hat r^2-r^2)$ and its estimate are
\begin{equation}
\label{eq:fixed_homo_scalar_variance_main}
\sigma_{r,\mathrm{re,h}}^2
=\mathcal V_{\mathrm{re,h}}(\mtx\Sigma_b,\mtx\Sigma_e;\rho^2,\sigma^2;g_2,g_3,g_4),
\qquad
\hat\sigma_{r,\mathrm{re,h}}^2
=\mathcal V_{\mathrm{re,h}}(\widehat{\mtx\Sigma}_b,\widehat{\mtx\Sigma}_e;\hat\rho^2,\hat\sigma^2;g_2,g_3,g_4).
\end{equation}

\begin{theorem}[Conditional random-effects inference]
\label{thm:asymptotics_random_fixed_homo_mt}
Let $\mtx X\in\mathbb R^{n\times p}$ be deterministic and nonzero, with
$p\asymp n$. Suppose, for constants $0<c<C<\infty$, that
\[
\left\|\frac{n\mtx X\mtx X^\top}{\|\mtx X\|_F^2}\right\|\le C,
\qquad
\frac1n\left\|
\frac{n\mtx X\mtx X^\top}{\|\mtx X\|_F^2}-\mtx I_n
\right\|_F^2\ge c.
\]
Let $\mtx B$ and $\mtx E$ be independent with
$\vecop(\mtx B)\sim
\mathcal N(\vct0,p^{-1}\mtx\Sigma_\beta\otimes\mtx I_p)$ and
$\vecop(\mtx E)\sim
\mathcal N(\vct0,\mtx\Sigma_e\otimes\mtx I_n)$, and set
$\mtx\Sigma_b=\|\mtx X\|_F^2\mtx\Sigma_\beta/(np)$. Suppose
$\rho^2,\sigma^2\in[c,C]$ and
\[
\frac qn\{1\vee\|\mtx\Sigma_b\|\vee\|\mtx\Sigma_e\|\}\to0.
\]
Then, conditionally on $\mtx X$,
\[
\frac{\sqrt n(\hat r^2-r^2)}{\hat\sigma_{r,\mathrm{re,h}}}
\Longrightarrow\mathcal N(0,1).
\]
\end{theorem}

The estimated standard error is $\hat\sigma_{r,\mathrm{re,h}}/\sqrt n$.
The corresponding normal-based interval has coverage approaching $1-\alpha$
when the predictor matrix is held fixed and the coefficients and noise are
redrawn, for every fixed $\alpha\in(0,1)$.

\begin{remark}
Equations~(16)--(19) of \citet{ge2016estimating} already give the general
quadratic-form covariance calculation and delta-method variance; their
simplification for unrelated subjects follows in equations~(20)--(21).
Our variance function is an expanded form of their unsimplified calculation
after matching normalizations. The theorem establishes when estimating this
variance yields valid normal inference, with fixed or diverging $q$.
In particular, bounded response-covariance eigenvalues permit $q=o(n)$;
growing eigenvalues are allowed under the stated joint growth condition.
The design conditions bound the largest normalized eigenvalue and keep the
two moment equations distinct, without requiring a Gaussian predictor
matrix.
\end{remark}

The interval is intended for $r^2$ away from zero and one.
Proof details are given in Supplementary
Section~\ref{sec:supp_fixed_design_theorem}.

\section{Random-design random-effects inference}
\label{sec:random_mt}

\subsection{Model, target, and point estimation}
\label{sec:random_model}
\label{sec:random_estimator}

In the third setting, repeated samples also redraw the predictors. Their
population covariance is $\mtx\Sigma$, and each row of the random
coefficient matrix has covariance $\mtx\Sigma_\beta/p$. We also allow the
noise covariance to differ between observations: the $i$th noise vector has
covariance $\mtx\Sigma_i$. The signal covariance, averaged over the
predictor population, and the average noise covariance are
\[
\mtx\Sigma_b=\frac{\tr(\mtx\Sigma)}p\mtx\Sigma_\beta,
\qquad
\bar{\mtx\Sigma}_e=\frac1n\sum_{i=1}^n\mtx\Sigma_i,
\qquad
\rho^2=\frac1q\tr(\mtx\Sigma_b),\quad
\sigma^2=\frac1q\tr(\bar{\mtx\Sigma}_e),\quad
r^2=\frac{\rho^2}{\rho^2+\sigma^2}.
\]
The difference from Section~\ref{sec:random_fixed_homo_mt} is in the signal
definition. There we used the average squared size of the observed predictors.
Here we use its population counterpart, $p^{-1}\tr(\mtx\Sigma)$. The
noise contribution still depends only on the average covariance. Thus two
sets of observation-specific noise covariances with the same average have the
same target, although they may lead to different standard errors.

When all noise covariances are equal, this is the overall trace-based
heritability parameter of \citet{ge2016estimating}. Allowing them to differ
extends the single-response treatment of unequal noise variances in
\citet{hu2022optimal} to several correlated responses.

We use the same two data summaries as in the fixed-design model and the same
formulas to estimate signal and noise:
\begin{equation}
\label{eq:random_estimators_mt}
\widehat{\mtx\Sigma}_b
=\frac{
\mtx Y^\top\mtx X\mtx X^\top\mtx Y/\|\mtx X\|_F^2
-\mtx Y^\top\mtx Y/n}
{n\tr\{(\mtx X\mtx X^\top)^2\}/\|\mtx X\|_F^4-1},
\qquad
\widehat{\bar{\mtx\Sigma}}_e
=\frac{\mtx Y^\top\mtx Y}{n}-\widehat{\mtx\Sigma}_b.
\end{equation}
Consequently,
\begin{equation}
\label{eq:random_design_scalar_estimators_mt}
\hat\rho^2
=\frac{q^{-1}\left\{
\tr(\mtx Y^\top\mtx X\mtx X^\top\mtx Y)/\|\mtx X\|_F^2
-\tr(\mtx Y^\top\mtx Y)/n\right\}}
{n\tr\{(\mtx X\mtx X^\top)^2\}/\|\mtx X\|_F^4-1},
\qquad
\hat\sigma^2=\frac{\tr(\mtx Y^\top\mtx Y)}{nq}-\hat\rho^2,
\qquad
\hat r^2=\frac{\hat\rho^2}{\hat\rho^2+\hat\sigma^2}.
\end{equation}
The denominator is again the squared distance of the normalized
predictor-similarity matrix from the identity, divided by $n$:
\[
\frac1n\left\|
\frac{n\mtx X\mtx X^\top}{\|\mtx X\|_F^2}-\mtx I_n
\right\|_F^2.
\]
The same data therefore give the same point estimate in Sections~3 and~4,
but the two analyses answer different sampling questions. Section~3 holds
predictors fixed; this section redraws them and defines the signal through
their population distribution. When all noise covariances are equal, the
leading variance formula below reduces to the fixed-design formula at the
same covariance and design inputs. Under our assumptions, the difference
between the observed and population predictor scales is negligible for this
normal approximation. Noise covariances that differ between observations add
a separate term to the variance.

\subsection{General heterogeneity: normal approximation}
\label{sec:random_asymptotics}

As in Section~3, the variance calculation uses
\begin{equation}
\label{eq:random_design_gk_main}
g_k=\frac1n\tr\!\left\{\left(
\frac{n\mtx X\mtx X^\top}{\|\mtx X\|_F^2}
\right)^k\right\},\qquad k=2,3,4,
\end{equation}
together with the average squared difference between each noise covariance
and their common average:
\[
\kappa_{\mathrm{tot}}
=\frac1n\sum_{i=1}^n
\|\mtx\Sigma_i-\bar{\mtx\Sigma}_e\|_F^2.
\]
The numbers $g_2,g_3,g_4$ describe the observed predictor matrix.
The quantity $\kappa_{\mathrm{tot}}$ measures how much the noise
covariances differ between observations; it is zero when they are all equal.
These differences are what we mean by noise heterogeneity.

The point estimator needs neither the individual $\mtx\Sigma_i$ nor
$\kappa_{\mathrm{tot}}$. Its standard error, however, needs the latter.
The two second-moment equations do not, in general, determine this additional
quantity. A fourth moment of the responses also involves other features of
the individual covariance matrices, so it does not solve the problem without
further assumptions. We therefore first prove a normal approximation using
the true, unknown variance. This is an oracle result: it describes the
estimator's uncertainty, but does not yet give a confidence interval that can
be computed entirely from the data.

For symmetric $q\times q$ matrices $A,E$, scalars $s,t$ with $s+t\ne0$,
$\kappa\ge0$, and $u_2>1$, define
\begingroup\small
\begin{equation}
\label{eq:random_variance_function}
\begin{aligned}
\mathcal V_{\mathrm{re}}(A,E;s,t;\kappa,u_2,u_3,u_4)
&=\frac{2}{q^2(s+t)^4}\Bigg[\|A\|_F^2\Bigg\{u_2t^2+2\frac{u_3-u_2^2}{u_2-1}t(s+t)+\frac{u_4-2u_2u_3+u_2^3}{(u_2-1)^2}(s+t)^2\Bigg\}\\
&\qquad+2\tr(EA)\Bigg\{t^2+\frac{u_3-u_2^2}{(u_2-1)^2}(s+t)^2\Bigg\}+\|E\|_F^2\Bigg\{s^2+\frac{(s+t)^2}{u_2-1}\Bigg\}+\kappa s^2\Bigg].
\end{aligned}
\end{equation}
\endgroup
Using the observed $g_2,g_3,g_4$ and the true covariance quantities, define
the leading variance of $\sqrt n(\hat r^2-r^2)$ by
\begin{equation}
\label{eq:random_scalar_variance_main}
\sigma_{r,\mathrm{re}}^2
:=\mathcal V_{\mathrm{re}}(\mtx\Sigma_b,\bar{\mtx\Sigma}_e;
\rho^2,\sigma^2;\kappa_{\mathrm{tot}},g_2,g_3,g_4).
\end{equation}

\begin{theorem}[Oracle normality under general heterogeneity]
\label{thm:asymptotics_mt}
Let $p\asymp n$ and
$\mtx X=\mtx Z\mtx\Sigma^{1/2}$, where the entries of $\mtx Z$ are i.i.d.,
centred, standardized and uniformly sub-Gaussian. Suppose
$p^{-1}\tr(\mtx\Sigma)>0$ and
$\|\mtx\Sigma\|/\{p^{-1}\tr(\mtx\Sigma)\}\le C$. Independently of $\mtx X$, let
$\vecop(\mtx B)\sim
\mathcal N(\vct0,p^{-1}\mtx\Sigma_\beta\otimes\mtx I_p)$ and let
$\vct e_i\sim\mathcal N(\vct0,\mtx\Sigma_i)$ be mutually independent and
independent of $\mtx B$. Set
$\mtx\Sigma_b=\{p^{-1}\tr(\mtx\Sigma)\}\mtx\Sigma_\beta$ and suppose
$\rho^2,\sigma^2\in[c,C]$. If
\[
\{1\vee\|\mtx\Sigma_b\|\vee\max_i\|\mtx\Sigma_i\|\}\sqrt{q/n}\to0,
\qquad
\{1\vee\|\mtx\Sigma_b\|\vee\max_i\|\mtx\Sigma_i\|\}^2
\sqrt{\log p/n}\to0,
\]
then
\[
\frac{\sqrt n(\hat r^2-r^2)}{\sigma_{r,\mathrm{re}}}
\Longrightarrow\mathcal N(0,1).
\]
\end{theorem}

\begin{remark}
Theorem~\ref{thm:asymptotics_mt} does not require
$\kappa_{\mathrm{tot}}=o(q)$: its oracle variance includes the heterogeneity
contribution. The covariance assumptions control its size through
$\kappa_{\mathrm{tot}}\le\{\max_i\|\mtx\Sigma_i\|\}
\tr(\bar{\mtx\Sigma}_e)$, allowing order $q$ even when the noise eigenvalues
are bounded. For a fixed average covariance, the target of the moment
estimator \citep{ge2016estimating} remains unchanged. The theorem alone
does not provide an estimator of this contribution for a data-based
confidence interval.

The heterogeneity correction is negligible relative to the remaining variance
if $\kappa_{\mathrm{tot}}=o(\|\mtx\Sigma_b\|_F^2+
\|\bar{\mtx\Sigma}_e\|_F^2)$. With bounded eigenvalues of
$\mtx\Sigma_b$ and $\bar{\mtx\Sigma}_e$, the remaining leading variance of
$\sqrt n(\hat r^2-r^2)$ is of order $q^{-1}$, while the correction is of
order $\kappa_{\mathrm{tot}}/q^2$; thus $\kappa_{\mathrm{tot}}=o(q)$
makes it negligible. For example, $\mtx\Sigma_b=\mtx I_q$ and noise
covariances $(1+a)\mtx I_q$ and $(1-a)\mtx I_q$ in equal proportions,
with fixed $0<a<1$, give $\kappa_{\mathrm{tot}}=a^2q$. This example
satisfies the response growth conditions when $q=o(n)$, yet the correction
and the remaining variance are both of order $q^{-1}$. Hence diverging $q$
alone does not justify omitting the correction.
\end{remark}

\subsection{Scalar heterogeneity: confidence intervals}
\label{sec:random_scalar_inference}

To estimate the extra variance, we now assume that each noise covariance is a
nonnegative multiple of the same matrix. This allows the overall noise level
to differ between observations, while keeping the relative variances and
correlations across responses the same wherever the noise level is nonzero.
We call this scalar heterogeneity. The added assumption is needed for the
standard error, not for the definition or point estimation of $r^2$.
Specifically, suppose
\begin{equation}
\label{eq:scalar_heterogeneity_mt}
\mtx\Sigma_i=\nu_i\mtx\Sigma_e,\qquad
\nu_i\ge0,\qquad \frac1n\sum_i\nu_i=1,
\qquad
\eta=\frac1n\sum_i(\nu_i-1)^2.
\end{equation}
The normalization gives $\bar{\mtx\Sigma}_e=\mtx\Sigma_e$, and
$\kappa_{\mathrm{tot}}=\eta\|\mtx\Sigma_e\|_F^2$. Thus we only need
to estimate one additional number, $\eta$, which measures variation in the
noise multipliers. For each observation, we sum the squared responses and
then square that sum. Averaging across observations gives the fourth-moment
statistic $n^{-1}\sum_i(\vct y_i^\top\vct y_i)^2$. Its expectation separates into
signal, signal--noise, and noise contributions. Subtracting the estimated
signal and signal--noise terms and solving for $\eta$ gives

\begingroup\small
\begin{equation}
\label{eq:eta_hat_scalar_main}
\begin{aligned}
\hat\eta
&=\{2\|\widehat{\bar{\mtx\Sigma}}_e\|_F^2+q^2\hat\sigma^4\}^{-1}\Bigg[\frac1n\sum_{i=1}^n(\vct y_i^\top\vct y_i)^2-\left\{\frac1n\sum_{i=1}^n\left(\frac{n\|\vct x_i\|_2^2}{\|\mtx X\|_F^2}\right)^2\right\}\{2\|\widehat{\mtx\Sigma}_b\|_F^2+q^2\hat\rho^4\}\\
&\qquad-2\{2\tr(\widehat{\mtx\Sigma}_b\widehat{\bar{\mtx\Sigma}}_e)+q^2\hat\sigma^2\hat\rho^2\}\Bigg]-1.
\end{aligned}
\end{equation}
\endgroup
Because $\eta$ is nonnegative, we replace a negative estimate by zero:
$\hat\eta_+=\max(\hat\eta,0)$. We then set
$\hat\kappa_{\mathrm{tot}}=\hat\eta_+\|\widehat{\bar{\mtx\Sigma}}_e\|_F^2$
and replace the unknown quantities in the variance formula by their estimates:
\begin{equation}
\label{eq:random_plugin_variance_main}
\hat\sigma_{r,\mathrm{re}}^2
:=\mathcal V_{\mathrm{re}}(\widehat{\mtx\Sigma}_b,
\widehat{\bar{\mtx\Sigma}}_e;\hat\rho^2,\hat\sigma^2;
\hat\kappa_{\mathrm{tot}},g_2,g_3,g_4).
\end{equation}

\begin{theorem}[Studentized inference under scalar heterogeneity]
\label{thm:studentized_random_mt}
Under the conditions of Theorem~\ref{thm:asymptotics_mt}, suppose in addition
that \eqref{eq:scalar_heterogeneity_mt} holds, $\max_i\nu_i\le C$, and
\[
\{1\vee\|\mtx\Sigma_b\|\vee\max_i\|\mtx\Sigma_i\|\}^3
\sqrt{(q+4\log p)/n}\to0.
\]
Then
\[
\frac{\sqrt n(\hat r^2-r^2)}{\hat\sigma_{r,\mathrm{re}}}
\Longrightarrow\mathcal N(0,1).
\]
\end{theorem}

The standard error $\hat\sigma_{r,\mathrm{re}}/\sqrt n$ can now be
computed from the data. The corresponding normal-based interval has coverage
approaching $1-\alpha$ for every fixed $\alpha\in(0,1)$.

\begin{remark}
The two theorems distinguish what is needed to estimate the signal fraction
from what is needed to estimate its standard error. General noise covariances
are allowed in the first result. Requiring a common covariance shape lets us
estimate the extra variance and obtain a confidence interval from the data.
This extends the single-response setting of \citet{hu2022optimal} to
correlated responses. Scalar heterogeneity does not mean independent
responses: their noise may remain strongly correlated. Compared with the
common-noise model of \citet{ge2016estimating}, the procedure accounts for
the additional uncertainty caused by different noise levels between
observations.
\end{remark}

The target uses $p^{-1}\tr(\mtx\Sigma)$, but the procedure needs only
the observed average squared predictor value,
$\|\mtx X\|_F^2/(np)$. There is no need to rescale all predictors to a
prescribed overall level. The growth conditions allow increasing response
dependence, and the variance formula determines the resulting precision.
As before, the interval is intended for $r^2$ bounded away from zero and one.
Supplementary Section~\ref{sec:supp_random_design_theorem} gives the
covariance calculation, normal approximation, and variance-estimation proof.

\section{Numerical studies}
\label{sec:simulation}

\subsection{Simulation conventions}
\label{sec:simulation_conventions}

We generate $\mtx Y=\mtx X\mtx B+\mtx E$ and assess point estimation,
standard errors, and 95\% interval coverage. Unless a setting specifies
otherwise, $q=20$, $\rho^2=1$, and $\sigma^2=0.5$, giving
$r^2=\rho^2/(\rho^2+\sigma^2)=2/3$. When varying $r^2$, we keep
$\sigma^2=0.5$ and set $\rho^2=\sigma^2r^2/(1-r^2)$. The tables and
figure captions give the dimension and sensitivity grids. We use 500 runs
per fixed-coefficient setting and 2000 per random-coefficient setting in the
main text.

The common response-covariance construction is
\[
\begin{gathered}
\mtx R_b=(0.8^{|j-k|})_{j,k\le q},\qquad
\mtx R_e=(0.5^{|j-k|})_{j,k\le q},\qquad v_j=j^{-1/2},\\
\mtx D_\pi=\diag(v_{\pi(1)}^{1/2},\ldots,v_{\pi(q)}^{1/2}),\qquad
\mtx\Sigma_e=c_e\mtx D_\pi\mtx R_e\mtx D_\pi,\qquad
c_e=\frac{q\sigma^2}{\sum_{j=1}^q v_j}.
\end{gathered}
\]
Here $\pi$ is a uniformly random permutation of $\{1,\ldots,q\}$,
drawn once per setting, so $q^{-1}\tr(\mtx\Sigma_e)=\sigma^2$.
For random coefficients, set $\mtx\Sigma_b=c_b\mtx R_b$ with
$c_b=\rho^2$, giving $q^{-1}\tr(\mtx\Sigma_b)=\rho^2$.
All population covariance matrices and permutations remain fixed across runs
within a setting. Noise rows are independent Gaussian vectors with the
covariances specified below. The predictors, random coefficients when
present, and noise are mutually independent. The noise-correlation and
response-dimension comparisons specify their alternative $\mtx R_e$ in
Sections~\ref{sec:simulation_noise_correlation}
and~\ref{sec:simulation_q_efficiency}.

For each fitted method, write $\hat V=\hat\sigma_r^2/n$ for its estimated
variance of $\hat r^2$. If $\hat V$ is finite and strictly positive, we
form the interval
\[
\bigl[\hat r^2-z_{0.975}\sqrt{\hat V},\;
      \hat r^2+z_{0.975}\sqrt{\hat V}\bigr].
\]
Both the point estimate and the interval are untruncated: neither is
restricted to $[0,1]$. We also leave negative variance-component estimates
and indefinite covariance estimates unchanged. If
$\hat V\le0$ or is nonfinite, the standard error and interval are recorded
as unavailable. Such runs are omitted from that method's average standard
error, average interval length, and coverage, but retained in the Monte Carlo
mean and empirical standard error of $\hat r^2$; no replacement runs are
generated. The simulation output records the proportion of valid standard
errors. All methods in the reported settings had valid standard errors in
every run. The nonnegative correction to $\hat\eta$ described below is
separate from this interval convention.

The Monte Carlo mean is the average estimate across runs, and the empirical
standard error is their sample standard deviation. The average estimated
standard error is the mean of $\sqrt{\hat V}$, and coverage is the
percentage of available intervals containing the true $r^2$. At 95\%
coverage, the simulation standard error is about 0.97 percentage points
with 500 runs and 0.49 with 2000 runs, so small coverage differences may
reflect simulation error alone.

\subsection{Fixed coefficients}
\label{sec:simulation_FE}

We first keep the regression coefficients fixed and generate new predictors
and noise in each run. Draw $\mtx B_0$ once per setting with independent
$\mathcal N(\vct0,\mtx R_b)$ rows and set
\[
\mtx B=\frac{\sqrt{q\rho^2}}{\|\mtx B_0\|_F}\mtx B_0,
\qquad q^{-1}\tr(\mtx B^\top\mtx B)=\rho^2.
\]
The resulting $\mtx B$ stays fixed. Each run draws a new $\mtx X$ with
independent $\mathcal N(0,1)$ entries and independent
$\mathcal N(\vct0,\mtx\Sigma_e)$ noise rows, using the covariance
convention above.

The mean estimates in Table~\ref{tab:dense_fixed} differ from the true
$r^2$ by at most 0.009. Coverage ranges from 94.0\% to 97.0\%, with no
clear worsening as $n$, $p$, $q$, or $r^2$ changes. These experiments do
not cover $r^2=0$ or $r^2=1$.

\begin{table}[!htbp]
\captionsetup{justification=justified,singlelinecheck=false}
\caption{Fixed-coefficient experiment with Gaussian predictors, using 500 runs per
setting and the scaling in Section~\ref{sec:simulation_conventions}.
The panel headings and columns give the dimensions and true signal fractions.}
\label{tab:dense_fixed}
\centering
\normalsize
\renewcommand{\arraystretch}{0.95}
\setlength{\tabcolsep}{8pt}
\begin{tabular}{
  c
  S[table-format=3.0]
  S[table-format=1.3]
  S[table-format=1.3]
  S[table-format=2.1]
}
\multicolumn{1}{c}{$(n,p)$}
  & \multicolumn{1}{c}{$q$}
  & \multicolumn{1}{c}{$r^2$}
  & \multicolumn{1}{c}{\shortstack{Monte Carlo\\mean}}
  & \multicolumn{1}{c}{\shortstack{Coverage\\(\%)}} \\
\noalign{\vskip2pt}
\multicolumn{5}{l}{\itshape Panel A: varying dimensions, $q=20$ and $r^2=0.667$} \\
$(100,100)$   & 20 & 0.667 & 0.658 & 94.2 \\
$(200,100)$   & 20 & 0.667 & 0.665 & 95.2 \\
$(400,100)$   & 20 & 0.667 & 0.664 & 94.4 \\
$(100,1000)$  & 20 & 0.667 & 0.661 & 96.6 \\
$(250,1000)$  & 20 & 0.667 & 0.666 & 96.6 \\
$(500,1000)$  & 20 & 0.667 & 0.666 & 94.8 \\
$(1000,1000)$ & 20 & 0.667 & 0.665 & 95.2 \\
$(2000,1000)$ & 20 & 0.667 & 0.666 & 95.4 \\
$(4000,1000)$ & 20 & 0.667 & 0.666 & 95.6 \\
$(5000,1000)$ & 20 & 0.667 & 0.667 & 96.6 \\
\noalign{\vskip2pt}
\multicolumn{5}{l}{\itshape Panel B: varying response dimension at $(n,p)=(400,100)$ and $r^2=0.667$} \\
$(400,100)$ &   1 & 0.667 & 0.662 & 95.8 \\
$(400,100)$ &   2 & 0.667 & 0.664 & 95.2 \\
$(400,100)$ &   5 & 0.667 & 0.667 & 95.0 \\
$(400,100)$ &  10 & 0.667 & 0.664 & 94.6 \\
$(400,100)$ &  20 & 0.667 & 0.664 & 94.4 \\
$(400,100)$ &  50 & 0.667 & 0.667 & 95.2 \\
$(400,100)$ & 100 & 0.667 & 0.666 & 95.6 \\
\noalign{\vskip2pt}
\multicolumn{5}{l}{\itshape Panel C: varying signal fraction at $(n,p,q)=(400,100,20)$} \\
$(400,100)$ & 20 & 0.10 & 0.100 & 97.0 \\
$(400,100)$ & 20 & 0.20 & 0.201 & 96.0 \\
$(400,100)$ & 20 & 0.35 & 0.350 & 95.4 \\
$(400,100)$ & 20 & 0.50 & 0.499 & 94.6 \\
$(400,100)$ & 20 & 0.65 & 0.648 & 94.6 \\
$(400,100)$ & 20 & 0.80 & 0.801 & 95.8 \\
$(400,100)$ & 20 & 0.90 & 0.901 & 94.0 \\
\end{tabular}
\end{table}
\subsection{Fixed predictors and random coefficients}
\label{sec:simulation_RE_fixed_homo}

We next keep the predictors fixed and redraw the coefficients and noise.
To choose a non-Gaussian predictor matrix for each $(n,p)$ pair, generate
\[
\mtx X_0=\mtx Z\mtx L_x^\top,\qquad
\mtx L_x\mtx L_x^\top=\mtx\Sigma_x
=(0.5^{|j-k|})_{j,k\le p},\qquad
Z_{ij}\overset{\mathrm{iid}}{\sim}\sqrt{5/7}\,t_7,
\]
where $\mtx L_x$ is the lower-triangular Cholesky factor and the innovations
have variance one. We normalize the realized matrix as
$\mtx X=\sqrt{np}\mtx X_0/\|\mtx X_0\|_F$, giving
$\|\mtx X\|_F^2/(np)=1$, and use this matrix in every run and in
the sensitivity settings with the same $(n,p)$. The theorem does not require
this scaling. Under it, $\mtx\Sigma_\beta=\mtx\Sigma_b$.
Each run draws independent coefficient rows from
$\mathcal N(\vct0,\mtx\Sigma_b/p)$ and independent noise rows from
$\mathcal N(\vct0,\mtx\Sigma_e)$, with both covariances fixed as in
Section~\ref{sec:simulation_conventions}.

The mean estimates in Table~\ref{tab:homo} differ from the true $r^2$
by at most 0.008, and coverage ranges from 94.4\% to 97.7\%.
For the smallest designs, the intervals contain the truth more often than
95\%; coverage is closer to 95\% for larger designs. The point estimator
matches \citet{ge2016estimating} after accounting for the different
normalizations. This experiment checks the finite-sample accuracy of the
standard error and confidence interval justified by
Theorem~\ref{thm:asymptotics_random_fixed_homo_mt}.

\begin{table}[!htbp]
\captionsetup{justification=justified,singlelinecheck=false}
\caption{Fixed-predictor experiment with random coefficients and the same noise
covariance for every observation, using 2000 runs per setting and the scaling
in Section~\ref{sec:simulation_conventions}. Panel A uses nine matrices generated from
$t_7$ variables. For these matrices the estimator's denominator stays away
from zero, and $\| X\|/\sqrt n$ ranges from 1.97 to 3.33. Panels B--C
use the same matrix at $(n,p)=(1000,1000)$.}
\label{tab:homo}
\centering
\normalsize
\renewcommand{\arraystretch}{0.95}
\setlength{\tabcolsep}{8pt}
\begin{tabular}{
  c
  S[table-format=3.0]
  S[table-format=1.3]
  S[table-format=1.3]
  S[table-format=2.1]
}
$(n,p)$ & {$q$} & {$r^2$}
  & \multicolumn{1}{c}{\shortstack{Monte Carlo\\mean}}
  & \multicolumn{1}{c}{\shortstack{Coverage\\(\%)}} \\
\noalign{\vskip2pt}
\multicolumn{5}{l}{\itshape Panel A: varying dimensions, $q=20$ and $r^2=0.667$} \\
$(100,100)$   & 20 & 0.667 & 0.659 & 97.7 \\
$(200,100)$   & 20 & 0.667 & 0.665 & 97.1 \\
$(400,100)$   & 20 & 0.667 & 0.663 & 96.9 \\
$(250,1000)$  & 20 & 0.667 & 0.666 & 95.5 \\
$(500,1000)$  & 20 & 0.667 & 0.667 & 95.2 \\
$(1000,1000)$ & 20 & 0.667 & 0.666 & 95.7 \\
$(2000,1000)$ & 20 & 0.667 & 0.666 & 94.6 \\
$(4000,1000)$ & 20 & 0.667 & 0.666 & 95.5 \\
$(5000,1000)$ & 20 & 0.667 & 0.666 & 95.3 \\
\noalign{\vskip2pt}
\multicolumn{5}{l}{\itshape Panel B: varying response dimension at $(n,p)=(1000,1000)$ and $r^2=0.667$} \\
$(1000,1000)$ &   1 & 0.667 & 0.663 & 95.8 \\
$(1000,1000)$ &   2 & 0.667 & 0.665 & 94.4 \\
$(1000,1000)$ &   5 & 0.667 & 0.666 & 95.0 \\
$(1000,1000)$ &  10 & 0.667 & 0.665 & 94.8 \\
$(1000,1000)$ &  20 & 0.667 & 0.666 & 95.7 \\
$(1000,1000)$ &  50 & 0.667 & 0.666 & 95.5 \\
$(1000,1000)$ & 100 & 0.667 & 0.667 & 96.5 \\
\noalign{\vskip2pt}
\multicolumn{5}{l}{\itshape Panel C: varying signal fraction at $(n,p,q)=(1000,1000,20)$} \\
$(1000,1000)$ & 20 & 0.10 & 0.100 & 95.7 \\
$(1000,1000)$ & 20 & 0.20 & 0.200 & 95.4 \\
$(1000,1000)$ & 20 & 0.35 & 0.350 & 95.9 \\
$(1000,1000)$ & 20 & 0.50 & 0.500 & 95.4 \\
$(1000,1000)$ & 20 & 0.65 & 0.650 & 95.2 \\
$(1000,1000)$ & 20 & 0.80 & 0.799 & 95.3 \\
$(1000,1000)$ & 20 & 0.90 & 0.898 & 95.9 \\
\end{tabular}
\end{table}
\subsection{Random predictors and unequal noise levels}
\label{sec:simulation_RE_random_scalar}

Here we ask whether the intervals work when some observations have much more
noise than others. The baseline is $(n,p,q)=(1000,1000,20)$,
$r^2=0.5$, and $\eta=30$, with $\sigma^2=0.5$ and
$\rho^2=\sigma^2r^2/(1-r^2)$. We generate new predictors, coefficients,
and noise in every run. Specifically, generate correlated, genotype-like
predictors from population-standardized binomial counts:
\[
\begin{gathered}
\mtx X=\mtx Z\mtx L_x^\top,\qquad
\mtx L_x\mtx L_x^\top=\mtx\Sigma_x
=(0.5^{|j-k|})_{j,k\le p},\\
Z_{ij}=\frac{G_{ij}-0.6}{\sqrt{0.42}},\qquad
G_{ij}\overset{\mathrm{iid}}{\sim}\operatorname{Binomial}(2,0.3).
\end{gathered}
\]
Here $\mtx L_x$ is the lower-triangular Cholesky factor. The innovations
are bounded, centred, and have variance one, as required by the design
assumptions in Section~4. We apply no sample column standardization or
overall rescaling to $\mtx X$. Set
$\mtx\Sigma_\beta=\mtx\Sigma_b=\rho^2\mtx R_b$ and
$\mtx\Sigma_e=c_e\mtx D_\pi\mtx R_e\mtx D_\pi$, with the response
shapes and trace normalization in Section~\ref{sec:simulation_conventions}.
The equality $\mtx\Sigma_\beta=\mtx\Sigma_b$ follows from
$p^{-1}\tr(\mtx\Sigma_x)=1$. Independently of $\mtx X$, draw
independent coefficient rows from $\mathcal N(\vct0,\mtx\Sigma_b/p)$
and independent noise rows from
$\mathcal N(\vct0,\nu_i\mtx\Sigma_e)$, independently of the
coefficients. Thus $\nu_i$ changes observation $i$'s noise level while
preserving the pattern of correlation across responses.

To construct the multipliers, set
\[
w_i=\frac{n a_n^{i-1}}{\sum_{k=0}^{n-1}a_n^k},\qquad
i=1,\ldots,n,
\]
choose $a_n\in(0,1]$ so that $n^{-1}\sum_i(w_i-1)^2=\eta$, and
randomly permute $(w_i)$ once to obtain $(\nu_i)$. When $\eta=0$,
$a_n=1$ and every multiplier equals one. Consequently,
$n^{-1}\sum_i\nu_i=1$, $\bar{\mtx\Sigma}_e=\mtx\Sigma_e$, and
$\kappa_{\mathrm{tot}}=\eta\|\mtx\Sigma_e\|_F^2$. The covariance
matrices and multipliers stay fixed within each setting, and the estimator
does not observe the individual multipliers. Larger $\eta$ means greater
differences in noise levels. In the notation of \citet{hu2022optimal},
$\kappa_0=1+\eta$; our $\eta=30$ case is comparable to their setting
with large noise differences.

For each simulated data set, compute the moment estimate $\hat r^2$ and
the observed design moments $g_2,g_3,g_4$ in
\eqref{eq:random_design_gk_main}. We compare three versions of the
variance formula \eqref{eq:random_plugin_variance_main}, setting
\[
\hat\kappa_{\mathrm{tot}}
=\eta^\star\|\widehat{\bar{\mtx\Sigma}}_e\|_F^2,
\qquad \eta^\star\in\{0,\hat\eta_+,\eta\},\qquad
\hat\eta_+=\max(\hat\eta,0),
\]
where $\hat\eta$ is the fourth-moment estimate in
\eqref{eq:eta_hat_scalar_main}. These choices respectively ignore
heterogeneity, estimate it, or use its true value as a benchmark. All three
use the same point estimate, estimated response covariances, and observed
design moments; the benchmark supplies only the true $\eta$, not the full
true variance. The intervals and treatment of invalid variances follow
Section~\ref{sec:simulation_conventions}. Table~\ref{tab:scalar_main}
varies $\eta$, $n$, $q$, and $r^2$ one at a time around the baseline.

Ignoring unequal noise levels makes the intervals miss the truth too often
(Table~\ref{tab:scalar_main}). At the baseline $\eta=30$, coverage is
84.8\% without the correction, 95.0\% with estimated $\eta$, and 95.7\%
with true $\eta$. The loss of coverage is largest at intermediate signal
fractions. The estimate of $\eta$ tends to be too low when noise levels
differ greatly, but improves as $n$ increases. Using estimated rather than
true $\eta$ changes coverage by at most 0.9 percentage points in panels
B--C. Neither version is always close to 95\%: both fall below 95\% at
$q=1,2$ and exceed 98\% at $q=100$. Estimation of $\eta$ is also less
accurate at $r^2=0.9$. Thus the correction helps, but does not remove all
error in these settings.

\begin{table}[!htbp]
\captionsetup{justification=justified,singlelinecheck=false}
\caption{Interval coverage when noise levels differ between observations, using
correlated genotype-like predictors and 2000 runs per setting. The baseline
is $(n,p,q)=(1000,1000,20)$, $r^2=0.5$, and $\eta=30$. Panels A--D
vary $\eta$, $n$, $q$, and $r^2$, respectively; other parameters stay at
baseline, with the scaling in Section~\ref{sec:simulation_conventions}.
The column $\overline{\hat\eta_+}$ is the
average estimate of the noise-level differences measured by $\eta$.}
\label{tab:scalar_main}
\centering
\normalsize
\renewcommand{\arraystretch}{0.95}
\setlength{\tabcolsep}{4.5pt}
\begin{tabular}{
  c
  S[table-format=3.0]
  S[table-format=1.2]
  S[table-format=2.1]
  S[table-format=2.2]
  S[table-format=2.1]
  S[table-format=2.1]
  S[table-format=2.1]
}
$(n,p)$ & {$q$} & {$r^2$} & {$\eta$} & {$\overline{\hat\eta_+}$}
  & \multicolumn{1}{c}{\shortstack{$\eta=0$\\coverage (\%)}}
  & \multicolumn{1}{c}{\shortstack{Estimated $\eta$\\coverage (\%)}}
  & \multicolumn{1}{c}{\shortstack{True $\eta$\\coverage (\%)}} \\
\noalign{\vskip2pt}
\multicolumn{8}{l}{\itshape Panel A: varying heterogeneity} \\
$(1000,1000)$ & 20 & 0.5 &  0.0 &  0.02 & 95.1 & 95.1 & 95.1 \\
$(1000,1000)$ & 20 & 0.5 &  1.5 &  1.46 & 94.8 & 95.3 & 95.4 \\
$(1000,1000)$ & 20 & 0.5 &  5.0 &  4.87 & 93.4 & 95.1 & 95.4 \\
$(1000,1000)$ & 20 & 0.5 & 15.0 & 14.24 & 90.3 & 94.8 & 95.1 \\
$(1000,1000)$ & 20 & 0.5 & 30.0 & 27.38 & 84.8 & 95.0 & 95.7 \\
\noalign{\vskip2pt}
\multicolumn{8}{l}{\itshape Panel B: varying sample size} \\
$(250,1000)$  & 20 & 0.5 & 30.0 & 22.57 & 90.8 & 96.6 & 97.1 \\
$(500,1000)$  & 20 & 0.5 & 30.0 & 25.39 & 86.4 & 94.5 & 95.4 \\
$(1000,1000)$ & 20 & 0.5 & 30.0 & 27.38 & 84.8 & 95.0 & 95.7 \\
$(2000,1000)$ & 20 & 0.5 & 30.0 & 28.61 & 86.6 & 95.9 & 96.0 \\
$(4000,1000)$ & 20 & 0.5 & 30.0 & 29.19 & 87.2 & 94.9 & 95.0 \\
\noalign{\vskip2pt}
\multicolumn{8}{l}{\itshape Panel C: varying response dimension} \\
$(1000,1000)$ &   1 & 0.5 & 30.0 & 27.68 & 84.1 & 93.7 & 93.5 \\
$(1000,1000)$ &   2 & 0.5 & 30.0 & 26.98 & 85.2 & 93.6 & 94.5 \\
$(1000,1000)$ &   5 & 0.5 & 30.0 & 26.88 & 85.0 & 94.9 & 95.5 \\
$(1000,1000)$ &  10 & 0.5 & 30.0 & 27.46 & 84.2 & 94.3 & 94.8 \\
$(1000,1000)$ &  20 & 0.5 & 30.0 & 27.38 & 84.8 & 95.0 & 95.7 \\
$(1000,1000)$ &  50 & 0.5 & 30.0 & 27.74 & 83.9 & 96.7 & 96.9 \\
$(1000,1000)$ & 100 & 0.5 & 30.0 & 27.81 & 87.6 & 98.5 & 98.7 \\
\noalign{\vskip2pt}
\multicolumn{8}{l}{\itshape Panel D: varying signal fraction} \\
$(1000,1000)$ & 20 & 0.1 & 30.0 & 27.78 & 92.7 & 95.7 & 96.0 \\
$(1000,1000)$ & 20 & 0.2 & 30.0 & 27.45 & 88.5 & 94.9 & 95.2 \\
$(1000,1000)$ & 20 & 0.35 & 30.0 & 27.42 & 86.4 & 95.4 & 96.0 \\
$(1000,1000)$ & 20 & 0.5 & 30.0 & 27.38 & 84.8 & 95.0 & 95.7 \\
$(1000,1000)$ & 20 & 0.65 & 30.0 & 27.70 & 88.3 & 95.4 & 96.1 \\
$(1000,1000)$ & 20 & 0.8 & 30.0 & 27.18 & 91.3 & 96.0 & 96.5 \\
$(1000,1000)$ & 20 & 0.9 & 30.0 & 21.96 & 94.6 & 96.1 & 97.3 \\
\end{tabular}
\end{table}
\subsection{Accounting for noise correlation}
\label{sec:simulation_noise_correlation}

Does the standard error need to account for correlations between response
errors? We compare the standard error in
\eqref{eq:fixed_homo_scalar_variance_main} with one that sets the
off-diagonal entries of $\widehat{\mtx\Sigma}_e$ to zero. Both use the
same estimate of $r^2$ and the same estimated signal covariance; only the
treatment of noise correlation differs.

As $\rho_e$ increases from 0 to 0.9, the observed standard deviation of
$\hat r^2$ rises from 0.0272 to 0.0316
(Figure~\ref{fig:noise_correlation}). The standard error that includes
noise correlations follows this rise, with coverage between 95.7\% and
96.0\%. The version that ignores them stays near 0.0289, and its coverage
falls from 96.0\% to 92.4\%. Noise correlation leaves the signal fraction
unchanged here, but ignoring it understates uncertainty.

\begin{figure}[!htbp]
\centering
\includegraphics[width=0.86\textwidth]{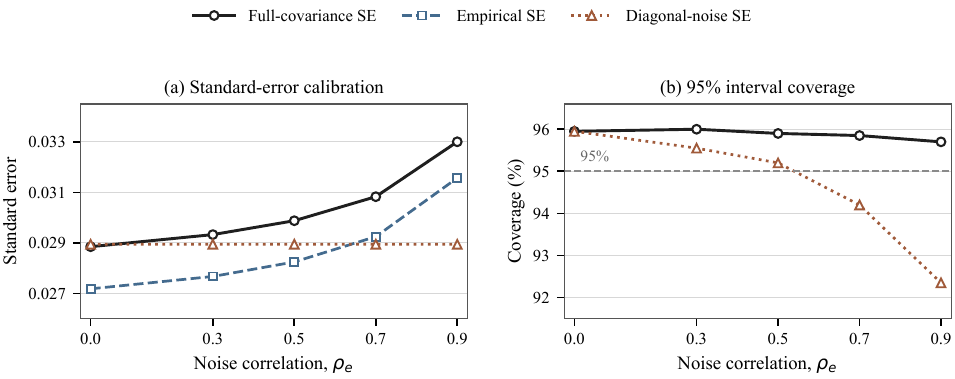}
\captionsetup{justification=justified,singlelinecheck=false}
\caption{Changing noise correlation while keeping each response's noise variance
fixed. We use the fixed $t_7$ design at $(n,p,q)=(1000,1000,20)$,
$\rho^2=1$, and $\sigma^2=0.5$, with 2000 runs per setting. We set
$\Sigma_e=c_e D_\pi(\rho_e^{|j-k|})_{j,k\le q} D_\pi$
for $\rho_e=0,0.3,0.5,0.7,0.9$ and keep $\Sigma_b=c_b R_b$ fixed.
Panel (a) compares the observed standard deviation with the average estimated
standard errors. Panel (b) shows coverage; the horizontal line marks 95\%.}

\label{fig:noise_correlation}
\end{figure}
\subsection{Adding responses}
\label{sec:simulation_q_efficiency}

We compare interval lengths as the number of responses increases. The three
noise models have independent responses, correlations that decrease with
the distance between response indices (AR(1), $\rho_e=0.5$), or a noise
component shared by all responses. The third uses
\[
\mtx R_e^{\mathrm{sp}}
=\frac{q}{q+19}\left(\mtx I_q+19\vct u_q\vct u_q^\top\right),
\qquad \vct u_q=q^{-1/2}\vct 1_q.
\]
Its diagonal entries are one, and its largest eigenvalue is 20 times the
others. Increasing $q$ from 1 to 100 reduces the average 95\% interval
length by 79.1\%, 78.4\%, and 78.3\% in the three settings, respectively
(Figure~\ref{fig:q_efficiency}). The shared component raises the average
standard error by about 10\% at $q=5,20$ relative to independent noise.
Even so, adding responses gives substantially narrower intervals in all
three settings.

\begin{figure}[!htbp]
\centering
\includegraphics[width=0.98\textwidth]{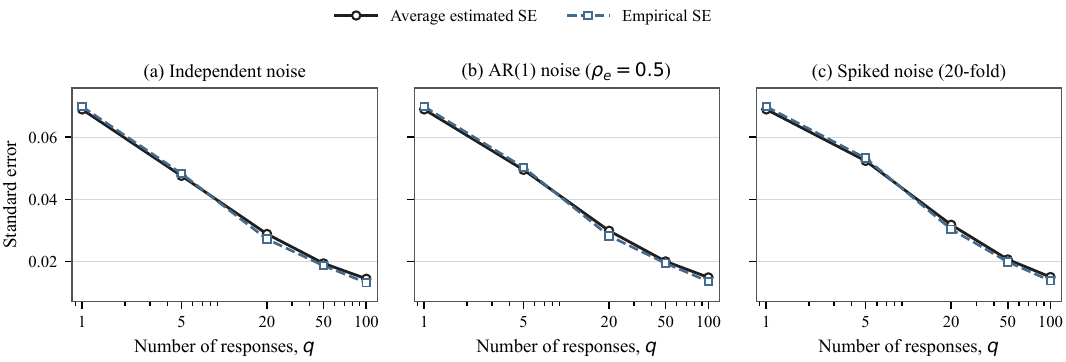}
\captionsetup{justification=justified,singlelinecheck=false}
\caption{Adding responses, with $q=1,5,20,50,100$. Panels (a)--(c) use independent
noise, AR(1) noise with $\rho_e=0.5$, and noise with a shared component.
All use the same fixed $t_7$ design at $(n,p)=(1000,1000)$,
$\rho^2=1$, and $\sigma^2=0.5$, with 2000 runs per setting. We keep
$\Sigma_b=c_b R_b$ and set $\Sigma_e=c_e D_\pi R_e D_\pi$, matching
each response's noise variance across panels. Solid curves show the average
estimated standard error; dashed curves show the standard deviation of the
estimates across runs. The interval length, without truncation, is
$2z_{0.975}\hat\sigma_r/\sqrt n$, so average length follows the solid
curve.}

\label{fig:q_efficiency}
\end{figure}
\section{Yeast growth across environments}
\label{sec:yeast_application}
\vspace{0.25\baselineskip}
\newcommand{\YeastEtaRaw}{0.108}
\newcommand{\YeastHeteroSEBase}{0.045636}
\newcommand{\YeastHeteroSE}{0.045643}
\newcommand{\YeastHeteroSEChange}{0.015}
\newcommand{\YeastHeteroPoint}{0.5428824}
\newcommand{\YeastHeteroCILower}{0.453}
\newcommand{\YeastHeteroCIUpper}{0.632}

We illustrate the overall additive marker-based signal fraction using yeast
growth across 20 environments. The BY$\times$RM study of
\citet{bloom2015interactions} provides 4390 strains, 28,220 markers, and
repeated growth measurements \citep{bloom2015data}. We average the available
measurements within each strain and environment and retain the 3474 strains
with at least one valid measurement in every environment. The released
measurements follow the source study's filtering and normalization protocol,
including spatial plate adjustment \citep{bloom2013finding,bloom2015interactions}.
The target concerns these strain means \citep{kruijer2015marker}: it is
the sum of the 20 signal variances divided by the sum of their total
variances, or equivalently a total-variance-weighted average of the
individual signal fractions. We retain the released response scales,
choosing a summary in which environments with more variation contribute more.

With the observed marker matrix held fixed, the Section~3 random-effects
analysis gives
\[
\hat r^2=0.543\quad\text{(s.e. }0.046\text{)},
\qquad 95\%\ \text{interval }[0.453,0.632].
\]
The point estimate combines the individual environments' signal and total
variance estimates. Its standard error also accounts for their dependence,
using Theorem~\ref{thm:asymptotics_random_fixed_homo_mt}.
We also illustrate Section~4's scalar heterogeneity correction, adopting
its random-design model and a common noise covariance shape with
strain-specific noise levels. This gives $\hat\eta_+=\hat\eta=\YeastEtaRaw$.
At the same inputs, the standard error changes from \YeastHeteroSEBase\
with $\kappa=0$ to \YeastHeteroSE\ with the correction, an increase of
\YeastHeteroSEChange\%; both analyses give the displayed point estimate
and interval after rounding. Thus the correction has little numerical
effect here, which does not establish equal noise covariances. Supplementary
Section~\ref{sec:yeast_supp} gives the data construction, centering, and
normalization used in these calculations.

Table~\ref{tab:yeast_environment_contributions} shows how each environment
contributes to the original-scale summary. YPD accounts for 31.4\% of the
estimated total variance and has signal fraction 0.584, contributing 0.184
to the joint ratio. The five largest weights, for YPD, YNB, Zeocin,
E6-Berbamine, and ManganeseSulfate, sum to 54.9\%. The weighted signal
fractions give 0.543; their arithmetic mean is 0.486. These are descriptive
point estimates under the additive marker model. The genetic interactions
identified by \citet{bloom2015interactions} may affect this variance
decomposition and induce residual dependence between strains. Both working
models assume independent residuals across strains; unequal variances and
correlations across environments are allowed within their common noise
covariance shape. Available replicate counts vary, so averaging alone does
not ensure equal noise covariance across strains. The scalar correction
addresses different noise levels, but not residual dependence between strains.

\begin{table}[!htbp]
\caption{Environment contributions on the released response scales for the
same 3474 strains and 28,220-marker panel. Weight is an environment's
percentage of the summed total variance; signal fraction is its fitted
signal variance divided by its total variance. Entries are point estimates,
rounded only for display.}
\label{tab:yeast_environment_contributions}
\centering\small
\begin{tabular}{lrrr}
\hline
Environment & Total variance & Weight (\%) & Signal fraction \\
\hline
CobaltChloride & 0.8546 & 5.23 & 0.395 \\
CopperSulfate & 0.4975 & 3.05 & 0.521 \\
Diamide & 0.6929 & 4.24 & 0.758 \\
E6-Berbamine & 0.9241 & 5.66 & 0.504 \\
Ethanol & 0.2077 & 1.27 & 0.388 \\
Formamide & 0.5615 & 3.44 & 0.370 \\
Hydroxyurea & 0.1149 & 0.70 & 0.092 \\
IndolaceticAcid & 0.3455 & 2.12 & 0.273 \\
Lactate & 0.5928 & 3.63 & 0.639 \\
Lactose & 0.5771 & 3.53 & 0.389 \\
MagnesiumChloride & 0.1749 & 1.07 & 0.386 \\
ManganeseSulfate & 0.8955 & 5.48 & 0.785 \\
Menadione & 0.6306 & 3.86 & 0.565 \\
Neomycin & 0.6910 & 4.23 & 0.563 \\
Raffinose & 0.4197 & 2.57 & 0.544 \\
Trehalose & 0.4954 & 3.03 & 0.429 \\
Xylose & 0.5007 & 3.07 & 0.414 \\
YNB & 1.0679 & 6.54 & 0.440 \\
YPD & 5.1327 & 31.44 & 0.584 \\
Zeocin & 0.9491 & 5.81 & 0.680 \\
\hline
\end{tabular}
\end{table}

\section{Discussion}
\label{sec:discussion}
The three analyses all estimate an overall explained fraction, but they differ
in what is held fixed in repeated samples. With fixed coefficients, the
signal is defined through new predictor values. With random coefficients and
a fixed predictor matrix, repeated samples redraw only the coefficients and
noise. The third setting also redraws the predictors and allows noise
covariance to differ between observations. The choice of model should match
the sampling question in the application. Across all three settings,
correlation affects precision: more responses provide the usual square-root
gain only when their variation is not too concentrated in shared directions.

Estimating the explained fraction requires less information than estimating
its standard error. Two moment equations separate average signal from average
noise without fitting the full coefficient matrix. The standard error must
also account for response correlations and the predictor matrix. When noise
covariances differ between observations, it needs an additional measure of
those differences. The two moment equations do not generally determine that
quantity. Assuming a common covariance shape allows us to estimate it from a
fourth moment. The measurement scales of the responses also matter for the
target and should be specified before analysis, rather than chosen from the
observed outcomes.

The assumptions limit where the confidence intervals apply. The fixed-effects
theorem requires Gaussian predictors with known covariance. The random-effects
theorems require Gaussian coefficients and noise; the random-design theorem
allows sub-Gaussian predictors, but not heavier-tailed ones. The normal
approximations apply to signal fractions bounded away from zero and one.
The argument for removing nuisance covariates uses a fixed projection and
common noise covariance, and does not cover a projection chosen from the
responses or general heterogeneous noise. Further work could address testing
at the endpoints, broader distributions, and more general differences in
noise covariance, possibly using repeated measurements.

\section*{Declaration of the use of generative AI and AI-assisted technologies}
\vspace{0.25\baselineskip}

During the preparation of this work, the authors used OpenAI Codex to assist
with improving mathematical proofs, implementing numerical experiments, and improving the exposition and organization of the
manuscript. The authors reviewed, edited and verified all AI-assisted material
and take full responsibility for the content of the publication.

\enlargethispage{\baselineskip}
\section*{Acknowledgement}
\vspace{0.25\baselineskip}

The authors thank Ethan Anderes, Jiming Jiang and Miles Lopes for their helpful
comments and suggestions. This work was supported by the U.S. National Science
Foundation under CAREER award DMS-1848575.

\ifdefined\mainsubmissiononly
\section*{Supplementary material}

The supplementary material contains all proofs, supporting moment
calculations, normal approximations and matrix concentration bounds.
\fi

\fi

\bibhang=1.7pc
\bibsep=0pt

\iftechnicalappendix

\ifmainpaper
\normalsize
\clearpage
\begin{center}
  {\Large\bfseries Supplementary Material}\par
  \vspace{0.5em}
  {\large Signal-to-Noise Ratio Inference for Multivariate
  High-Dimensional Linear Models}
\end{center}
\vspace{1.5em}
\fi
\setcounter{section}{0}
\setcounter{equation}{0}
\setcounter{table}{0}
\setcounter{figure}{0}
\renewcommand{\thesection}{S\arabic{section}}
\renewcommand{\thesubsection}{\thesection.\arabic{subsection}}
\renewcommand{\theequation}{S\arabic{equation}}
\renewcommand{\thetable}{S\arabic{table}}
\renewcommand{\thefigure}{S\arabic{figure}}
\renewcommand{\theHsection}{supp.\arabic{section}}
\renewcommand{\theHsubsection}{supp.\arabic{section}.\arabic{subsection}}
\renewcommand{\theHequation}{supp.\arabic{equation}}
\renewcommand{\theHtable}{supp.\arabic{table}}
\renewcommand{\theHfigure}{supp.\arabic{figure}}
\newcommand{\proofblock}[1]{\par\medskip\noindent\textit{#1}\par\nobreak\smallskip\nobreak}
\section*{Organization of the supplementary material}
Section~S1 gives the ratio identities and probability bounds used in the
proofs. Sections~S2--S4 establish inference for fixed effects, fixed-design
homogeneous random effects, and random-design heterogeneous random effects,
respectively. Section~S3 also justifies the fixed projection used to remove
an intercept. Section~S5 describes data processing and implementation details
for the yeast analysis.

\section{Preliminaries}
\label{sec:supp_preliminaries}
The proofs share the same final step: a normal approximation for one linear
combination of the variance-component errors is transferred to the SNR
estimator, and the oracle variance is then replaced by its estimate. We
first isolate this step using exact ratio identities. The remaining tools
are used to obtain the model-specific covariance and concentration bounds.

\subsection{Exact ratio identity and studentization}
\label{sec:common_ratio_reduction}
In all three regimes, order the variance components as noise followed by
signal, and write
\[
\vct\theta=\begin{pmatrix}\sigma^2\\ \rho^2\end{pmatrix},
\qquad
\widehat{\vct\theta}=\begin{pmatrix}\hat\sigma^2\\ \hat\rho^2\end{pmatrix},
\qquad
\vct\Delta=\widehat{\vct\theta}-\vct\theta,
\qquad
\tau=\rho^2+\sigma^2,\quad
\hat\tau=\hat\rho^2+\hat\sigma^2.
\]
For $h(s,t)=t/(s+t)$, the gradients at $\vct\theta$ and
$\widehat{\vct\theta}$ are
\begin{equation}
\label{eq:snr_gradient}
\vct d=\frac1{\tau^2}
\begin{pmatrix}-\rho^2\\ \sigma^2\end{pmatrix},
\qquad
\hat{\vct d}=\frac1{\hat\tau^2}
\begin{pmatrix}-\hat\rho^2\\ \hat\sigma^2\end{pmatrix}.
\end{equation}
Direct subtraction gives the exact identities, whenever $\hat\tau\ne0$,
\begin{align}
\hat r^2-r^2
&=\frac{\sigma^2(\hat\rho^2-\rho^2)
-\rho^2(\hat\sigma^2-\sigma^2)}{\tau\hat\tau}
=\frac{\tau}{\hat\tau}\,\vct d^\top\vct\Delta,
\label{eq:exact_ratio_identity}\\
\hat{\vct d}
&=\frac{\tau}{\hat\tau}\vct d
-\frac{\hat r^2-r^2}{\hat\tau}\vct 1_2,
\qquad \vct 1_2=(1,1)^\top.
\label{eq:exact_gradient_identity}
\end{align}
The second identity follows by writing
$\vct d=\tau^{-1}(-r^2,1-r^2)^\top$ and similarly for
$\hat{\vct d}$. The first identity transfers a normal approximation for
$\vct d^\top\vct\Delta$ to the ratio as soon as the total variance is
consistently estimated. The second controls the effect of estimating the
gradient in the standard error. Thus neither step requires a Taylor
remainder or a joint central limit theorem.

The next lemma allows a design-dependent covariance matrix, as needed in
S4. Its plug-in error is measured relative to the lower variance scale
$q^{-1}$, rather than merely required to vanish.

\begin{lemma}[Exact ratio reduction and studentization]
\label{lem:common_studentization}
Let $q=q_n\ge1$, let $\ell_n\ge1$, and let $\mtx V_n$ be a symmetric,
possibly random, $2\times2$ matrix. Suppose $\vct\theta$ remains in a fixed
compact subset of $(0,\infty)^2$ and, with probability tending to one,
\begin{equation}
\label{eq:common_variance_scale}
\frac cq\mtx I_2\preceq\mtx V_n\preceq\frac{C\ell_n}{q}\mtx I_2.
\end{equation}
Assume
\begin{equation}
\label{eq:common_total_consistency}
\frac{\hat\tau}{\tau}\to_P1,
\qquad \frac{\ell_n}{nq}\to0,
\end{equation}
and, for $\sigma_{r,n}^2=\vct d^\top\mtx V_n\vct d$,
\begin{equation}
\label{eq:common_linear_clt}
\frac{\sqrt n\,\vct d^\top\vct\Delta}{\sigma_{r,n}}
\Longrightarrow\mathcal N(0,1).
\end{equation}
Then, with probability tending to one,
$c/q\le\sigma_{r,n}^2\le C\ell_n/q$, and
\begin{equation}
\label{eq:common_oracle_clt}
\frac{\sqrt n(\hat r^2-r^2)}{\sigma_{r,n}}
\Longrightarrow\mathcal N(0,1).
\end{equation}
If, in addition, a symmetric matrix $\widehat{\mtx V}_n$ satisfies
\begin{equation}
\label{eq:common_plugin_matrix}
\|\widehat{\mtx V}_n-\mtx V_n\|=o_P(q^{-1}),
\end{equation}
then $\hat\sigma_{r,n}^2=\hat{\vct d}^{\top}
\widehat{\mtx V}_n\hat{\vct d}$ satisfies
\begin{equation}
\label{eq:common_variance_consistency}
\frac{\hat\sigma_{r,n}^2}{\sigma_{r,n}^2}\to_P1,
\qquad \PP(\hat\sigma_{r,n}^2>0)\to1,
\end{equation}
and
\begin{equation}
\label{eq:common_studentized_clt}
\frac{\sqrt n(\hat r^2-r^2)}{\hat\sigma_{r,n}}
\Longrightarrow\mathcal N(0,1).
\end{equation}
Consequently, the corresponding two-sided Wald interval has asymptotic
coverage $1-\alpha$ for every fixed $\alpha\in(0,1)$.
\end{lemma}

\begin{proof}
We first establish the oracle limit. Compactness of $\vct\theta$ gives
$\tau\asymp1$ and $c\le\|\vct d\|_2\le C$, so
\eqref{eq:common_variance_scale} gives the asserted variance bounds. The
exact identity \eqref{eq:exact_ratio_identity}, together with
\eqref{eq:common_total_consistency} and \eqref{eq:common_linear_clt},
then proves \eqref{eq:common_oracle_clt} by Slutsky's theorem. This also gives
the error rate needed to estimate the gradient:
\[
\frac{\hat r^2-r^2}{\sigma_{r,n}}=O_P(n^{-1/2}),
\qquad
\hat r^2-r^2=O_P\!\left(\sqrt{\frac{\ell_n}{nq}}\right)=o_P(1),
\]
and \eqref{eq:exact_gradient_identity} implies
$\|\hat{\vct d}\|_2=O_P(1)$. These two consistency statements also give
consistency of the individual components, without a joint limit theorem:
\begin{align*}
\hat\rho^2-\rho^2
&=\hat\tau(\hat r^2-r^2)+r^2(\hat\tau-\tau)=o_P(1),\\
\hat\sigma^2-\sigma^2
&=(1-r^2)(\hat\tau-\tau)-\hat\tau(\hat r^2-r^2)=o_P(1).
\end{align*}
In particular, both estimated components are positive with probability
tending to one. This justifies evaluating variance functions whose scalar
arguments are required to be positive; positivity of the estimated
covariance matrices themselves is not needed.

For studentization, first hold $\mtx V_n$ fixed and replace only the
gradient. Applying the reverse triangle inequality to
\eqref{eq:exact_gradient_identity} in the norm induced by $\mtx V_n$ gives,
on the event $\hat\tau>0$ whose probability tends to one,
\begin{align*}
\left|
\frac{(\hat{\vct d}^{\top}\mtx V_n\hat{\vct d})^{1/2}}{\sigma_{r,n}}
-\frac{\tau}{\hat\tau}\right|
&\le
\frac{|\hat r^2-r^2|}{\hat\tau\sigma_{r,n}}
(\vct 1_2^\top\mtx V_n\vct 1_2)^{1/2}\\
&=O_P\!\left(\sqrt{\frac{\ell_n}{nq}}\right)=o_P(1).
\end{align*}
Thus $\hat{\vct d}^{\top}\mtx V_n\hat{\vct d}/\sigma_{r,n}^2\to_P1$.
Replacing the covariance matrix itself contributes only
\[
\frac{|\hat{\vct d}^{\top}(\widehat{\mtx V}_n-\mtx V_n)
\hat{\vct d}|}{\sigma_{r,n}^2}
\le Cq\|\widehat{\mtx V}_n-\mtx V_n\|\,\|\hat{\vct d}\|_2^2
=o_P(1).
\]
This proves ratio consistency and positivity of the estimated variance.
Slutsky's theorem proves \eqref{eq:common_studentized_clt}, and coverage
follows by continuity of the standard normal law.
\end{proof}

\subsection{Gaussian and concentration tools}
\label{sec:prelim_gaussian_tools}

\proofblock{Gaussian quadratic forms and Wishart moments.}
For the random-effects proofs, conditioning on the design reduces each
linear combination of the two moments to a centered Gaussian quadratic
form. The following distribution-function bound will therefore provide the
normal approximation in S3 and S4. Write $\Phi$ for the standard-normal
distribution function.
Let $\vct z\sim\mathcal N(\vct 0,\mtx I_m)$ and let $\mtx A$ be a
deterministic symmetric matrix with $\|\mtx A\|_F>0$. Since
$\Var(\vct z^\top\mtx A\vct z)=2\|\mtx A\|_F^2$, the classical
Berry--Esseen inequality gives
\begin{equation}
\label{eq:gaussian_quadratic_be}
\sup_{t\in\mathbb R}
\left|
\PP\left(
\frac{\vct z^\top\mtx A\vct z-\tr(\mtx A)}
{\sqrt2\|\mtx A\|_F}\le t\right)-\Phi(t)
\right|
\le C\frac{\|\mtx A\|}{\|\mtx A\|_F}.
\end{equation}
Indeed, if $\lambda_1,\ldots,\lambda_m$ are the eigenvalues of
$\mtx A$, orthogonal invariance gives
\[
\vct z^\top\mtx A\vct z-\tr(\mtx A)\overset d=
\sum_{r=1}^m\lambda_r(Z_r^2-1),
\]
where the $Z_r$ are independent standard normal variables.  Therefore
Berry--Esseen and $\E|Z_1^2-1|^3<\infty$ imply
\[
\sup_t\left|
\PP\left(
\frac{\vct z^\top\mtx A\vct z-\tr(\mtx A)}
{\sqrt2\|\mtx A\|_F}\le t\right)-\Phi(t)
\right|
\le
C\frac{\sum_r|\lambda_r|^3}
        {(\sum_r\lambda_r^2)^{3/2}}
\le C\frac{\max_r|\lambda_r|}
          {(\sum_r\lambda_r^2)^{1/2}},
\]
which is \eqref{eq:gaussian_quadratic_be}. The same argument applies
conditionally when $\mtx A$ is measurable with respect to the conditioning
sigma-field and $\vct z$ is conditionally standard Gaussian. The bound is
uniform over directions whenever the displayed norm ratio is uniformly
controlled; in S4 this will allow a direction chosen from the observed
design.

The fixed-effects covariance calculation instead requires moments of
Wishart traces. We use the following formula to collect these moments
without expanding individual entries.

\begin{lemma}[A standard trace formula]\label{lem:exp_tr_WH}
\leavevmode\par\noindent
For $d\times d$ symmetric matrices $H_1,\ldots,H_k$ and a
Wishart$(n,\mtx{\Sigma})$ random matrix
$\mtx{W}=\mtx{X}^\top\mtx{X}$,
the moment formula of \citet{letac2004wishart} and
\citet{graczyk2005wishart} gives

\begin{align}\label{eq:exp_tr_WH}
 \E \left\{ \operatorname{tr}(W H_1) \cdots \operatorname{tr}(W H_k) \right\} =
&\sum_{\pi\,\text{a permutation of }\{1,\ldots,k\}}
\notag\\[-2pt]
&\quad
 2^{k-\#\{\text{cycles of }\pi\}}n^{\#\{\text{cycles of }\pi\}}
 \prod_{(i_1\cdots i_\ell)\,\text{a cycle of }\pi}
 \operatorname{tr}\left(\prod_{r=1}^{\ell}\Sigma H_{i_r}\right).
\end{align}
\end{lemma}

\begin{proof}
For completeness, the coefficients can be read directly from the
Gaussian moment-generating function. Put
$\mtx L(\vct t)=\sum_{j=1}^k t_j\mtx H_j$. For $\vct t$ near zero,
independence of the $n$ Gaussian rows gives
\begin{align*}
\E\exp\{\tr(\mtx W\mtx L(\vct t))\}
&=\det\{\mtx I_d-2\mtx\Sigma^{1/2}\mtx L(\vct t)
                          \mtx\Sigma^{1/2}\}^{-n/2},\\
\log\E\exp\{\tr(\mtx W\mtx L(\vct t))\}
&=n\sum_{r\ge1}\frac{2^{r-1}}r
       \tr\{(\mtx\Sigma\mtx L(\vct t))^r\}.
\end{align*}
For a set of $r$ distinct labels, the corresponding mixed derivative of
the second line is the sum over cycles on those labels, each with
coefficient $2^{r-1}n$. The factor $r^{-1}$ cancels the $r$ cyclic
rotations of each product inside the trace. Differentiating the exponential
partitions $\{1,\ldots,k\}$ into such cycles. A permutation with $c$
cycles therefore has coefficient $2^{k-c}n^c$, proving
\eqref{eq:exp_tr_WH}. This derivation also allows singular
$\mtx\Sigma$, since no covariance inverse is used.
\end{proof}

We also use the first-order Gaussian Poincar\'e inequality
\citep{chatterjee2009fluctuations}. If
$\vct\xi=\vct\mu+\mtx\Omega^{1/2}\vct z$, with
$\vct z\sim\mathcal N(\vct0,\mtx I_m)$ and $\mtx\Omega\succeq0$,
then the chain rule gives, for the smooth polynomials used below,
\begin{equation}
\label{eq:gaussian_poincare_covariance}
\Var\{f(\vct\xi)\}
\le \E\{\nabla f(\vct\xi)^\top\mtx\Omega\nabla f(\vct\xi)\}
\le\|\mtx\Omega\|\,\E\|\nabla f(\vct\xi)\|_2^2.
\end{equation}
The same inequality applies conditionally when the mean and covariance
are fixed by the conditioning variable. It requires neither independent
coordinates of $\vct\xi$ nor invertibility of $\mtx\Omega$.

\proofblock{Sub-Gaussian concentration and matrix bounds.}

We need moment bounds both for individual design rows and for the operator
norm of the full design matrix. The row bound below is a standard
consequence of sub-Gaussian concentration; see
\citet{vershynin2012nonasymptotic}.
Let $\vct z\in\mathbb R^p$ have i.i.d.\ sub-Gaussian entries satisfying
$\E z_i=0$, $\Var(z_i)=1$, and $\|z_i\|_{\psi_2}=O(1)$. For every
deterministic positive semidefinite $p\times p$ matrix $\mtx\Sigma$,
\[
\E\|\mtx\Sigma^{1/2}\vct z\|^k
\le \|\mtx\Sigma\|^{k/2}C^k k^{k/2}p^{k/2},
\qquad k=1,2,\ldots.
\]

For the full matrix, the following tail bound also gives the fixed-order
moments used in S2 and S4. Let $\mtx Z\in\mathbb R^{n\times m}$
have i.i.d.\ centred, standardized, uniformly sub-Gaussian entries, let
$\mtx\Sigma\succeq0$ be deterministic and $m\times m$, and set
$\mtx X=\mtx Z\mtx\Sigma^{1/2}$. The rows of $\mtx Z$ are independent
isotropic sub-Gaussian vectors with uniformly bounded sub-Gaussian norm.
Theorem~5.39 of \citet{vershynin2012nonasymptotic} therefore gives
\[
\PP\{\|\mtx Z\|>\sqrt n+C\sqrt m+t\}\le2e^{-ct^2},
\qquad t\ge0.
\]
Using $\|\mtx X\|\le\|\mtx\Sigma\|^{1/2}\|\mtx Z\|$, we obtain
\begin{equation}
\label{eq:subg_matrix_norm_tail}
\PP\!\left\{\|\mtx X\|>
\|\mtx\Sigma\|^{1/2}(\sqrt n+C\sqrt m+t)\right\}
\le2e^{-ct^2},\qquad t\ge0.
\end{equation}
Integrating the preceding tail bound yields, for every fixed integer
$r\ge1$,
\begin{equation}
\label{eq:subg_matrix_norm_moments}
\{\E\|\mtx X\|^r\}^{1/r}\le C_r\sqrt n,
\qquad
\E\left\|\frac1n\mtx X^\top\mtx X\right\|^r
=\frac1{n^r}\E\|\mtx X\|^{2r}\le C_r,
\end{equation}
whenever $m=O(n)$ and $\|\mtx\Sigma\|=O(1)$. The constants may depend on
the fixed moment order and the uniform bounds, but not on the dimensions.
No inverse of $\mtx\Sigma$ is used, so the conclusion also holds for
singular $\mtx\Sigma$.

\begin{lemma}[Hanson--Wright inequality]\label{lem:Hanson-Wright}
The following sub-Gaussian form is due to \citet{rudelson2013hanson}.
Let $\vct z=(z_1,\ldots,z_m)^\top$ have independent coordinates satisfying
$\E z_i=0$ and $\max_i\|z_i\|_{\psi_2}=O(1)$, and let $\mtx A$ be an
$m\times m$ deterministic matrix. Then, for any $t>0$,
\[
\PP\{|\vct z^\top\mtx A\vct z-\E(\vct z^\top\mtx A\vct z)|>t\}
\le 2\exp\left\{-c\min\left(
\frac{t^2}{\|\mtx A\|_F^2},\frac{t}{\|\mtx A\|}
\right)\right\},
\]
where $c>0$ depends only on the uniform sub-Gaussian bound.
\end{lemma}

The next net bounds turn fixed-direction concentration into an
operator-norm bound. We use them with the preceding Hanson--Wright
inequality when estimating response covariance matrices.

We will also use the corresponding quadratic-form variance bound.
For i.i.d.\ centered, standardized entries and a symmetric matrix
$\mtx A$, the diagonal and off-diagonal terms in
\[
\vct z^\top\mtx A\vct z-\tr(\mtx A)
=\sum_i A_{ii}(z_i^2-1)+2\sum_{i<j}A_{ij}z_iz_j
\]
are pairwise uncorrelated. Hence
\begin{equation}
\label{eq:subg_quadratic_variance}
\Var(\vct z^\top\mtx A\vct z)
=2\|\mtx A\|_F^2+(\E z_1^4-3)\sum_i A_{ii}^2
\le C\|\mtx A\|_F^2.
\end{equation}
The constant is uniform under the sub-Gaussian assumption. For
$\vct x=\mtx\Sigma^{1/2}\vct z$, apply this bound with
$\mtx A=\mtx\Sigma^{1/2}\mtx M\mtx\Sigma^{1/2}$ to control
$\Var(\vct x^\top\mtx M\vct x)$.

\begin{lemma}[Net bounds for operator norms]
\label{lem:operator norm bound by net cover}
The following covering-number and net bounds are standard; see
\citet{vershynin2012nonasymptotic}.
For every $\varepsilon>0$, the unit sphere in $\mathbb R^d$ has an
$\varepsilon$-net $\mathcal N_d$ with
$|\mathcal N_d|\le(1+2/\varepsilon)^d$.  If $\mtx A$ is symmetric and
$\varepsilon<1/2$, then
\[
\|\mtx A\|
\le\frac1{1-2\varepsilon}
\max_{\vct u\in\mathcal N_d}|\vct u^\top\mtx A\vct u|.
\]
For a general $m\times d$ matrix $\mtx A$ and nets
$\mathcal N_m,\mathcal N_d$ with $\varepsilon<1$, one has
\[
\|\mtx A\|
\le\frac1{(1-\varepsilon)^2}
\max_{\substack{\vct u\in\mathcal N_m\\
                 \vct v\in\mathcal N_d}}
|\vct u^\top\mtx A\vct v|.
\]
\end{lemma}

\section{Proof of Theorem~\ref{thm:asymptotics_fixed_mt}}
\label{sec:supp_fixed_theorem}
After whitening the design, the two response moments are polynomials in
independent standard Gaussian variables. We compute their covariance using
Wishart identities, establish a normal approximation using derivative
bounds, and reuse the first derivatives to control the estimated
covariance. The ratio lemma in S1 then completes the proof.

\subsection{Moment calculations and covariance structure}
\label{sec:fixed_supporting}

Throughout S2, put
\[
\ell_n=1\vee\|\mtx\Sigma_b\|\vee\|\mtx\Sigma_e\|.
\]
The trace bounds imply $\ell_n\le Cq$, so the theorem's condition
$\ell_n q/n\to0$ gives
\begin{equation}
\label{eq:fe_spectral_growth_implications}
\frac{\ell_n^2}{n}\to0.
\end{equation}
This implication is used in the normal approximation. The original
condition $\ell_nq/n\to0$ controls the plug-in error, and also implies
$\ell_n/(nq)\to0$ for the ratio reduction.

To put the design in standard Gaussian form, define
\[
\widetilde{\mtx X}=\mtx X\mtx\Sigma_x^{-1/2},
\qquad
\widetilde{\mtx B}=\mtx\Sigma_x^{1/2}\mtx B,
\qquad
\mtx W=\widetilde{\mtx X}^{\top}\widetilde{\mtx X}.
\]
Then $\mtx W\sim\operatorname{Wishart}(n,\mtx I_p)$,
$\widetilde{\mtx B}^{\top}\widetilde{\mtx B}=\mtx\Sigma_b$, and
\[
\mtx Y=\widetilde{\mtx X}\widetilde{\mtx B}+\mtx E,
\qquad
\mtx X\mtx\Sigma_x^{-1}\mtx X^\top
=\widetilde{\mtx X}\widetilde{\mtx X}^{\top}.
\]
Whitening leaves the responses, the estimator and its target unchanged.
For the rest of S2, write $\mtx X$ and $\mtx B$ for
$\widetilde{\mtx X}$ and $\widetilde{\mtx B}$, respectively. Thus
$\mtx W=\mtx X^\top\mtx X$ is Wishart with identity covariance and
$\mtx B^\top\mtx B=\mtx\Sigma_b$.

The covariance calculation below involves products of quadratic forms in
$\mtx W$ and $\mtx W^2$, together with the matrix moment $\E(\mtx W^3)$.
We collect these identities first.

\begin{lemma}[Wishart identities through the fourth trace moment]
\label{lem:fixed_wishart_identities}
Let $\mtx W\sim\operatorname{Wishart}(n,\mtx I_p)$.  For any
$\vct\alpha,\vct\beta\in\mathbb R^p$,
\begin{align*}
\E\!\left[(\vct\alpha^\top\mtx W\vct\alpha)
                (\vct\beta^\top\mtx W\vct\beta)\right]
&=2n(\vct\alpha^\top\vct\beta)^2
+n^2\|\vct\alpha\|_2^2\|\vct\beta\|_2^2,\\
\E\!\left[(\vct\alpha^\top\mtx W\vct\alpha)
                (\vct\beta^\top\mtx W^2\vct\beta)\right]
&=2n(2n+p+2)(\vct\alpha^\top\vct\beta)^2 +n(n^2+np+n+2)
\|\vct\alpha\|_2^2\|\vct\beta\|_2^2,\\
\E\!\left[(\vct\alpha^\top\mtx W^2\vct\alpha)
                (\vct\beta^\top\mtx W^2\vct\beta)\right]
&=2n(4n^2+5np+10n+p^2+5p+10)
   (\vct\alpha^\top\vct\beta)^2\\
&\quad+\bigl(n^4+2n^3p+2n^3+n^2p^2+2n^2p
 +11n^2+6np+10n\bigr) \times\|\vct\alpha\|_2^2\|\vct\beta\|_2^2.
\end{align*}
In addition,
\begin{equation}
\label{eq:wishart_third_matrix}
\E(\mtx W^3)
=n\{n^2+3np+3n+p^2+3p+4\}\mtx I_p.
\end{equation}
\end{lemma}

\begin{proof}
The first identity is the $k=2$ specialization of
Lemma~\ref{lem:exp_tr_WH} with
$H_1=\vct\alpha\vct\alpha^\top$ and
$H_2=\vct\beta\vct\beta^\top$.
For the two identities involving $\mtx W^2$, choose an orthonormal basis
$\vct u_1,\ldots,\vct u_p$ and put
\[
\mtx H_i(\vct v)=\frac12(\vct v\vct u_i^\top+\vct u_i\vct v^\top).
\]
This expresses a quadratic form in $\mtx W^2$ as a sum of products of
linear Wishart traces:
\[
\vct v^\top\mtx W^2\vct v
=\sum_i\tr^2\{\mtx W\mtx H_i(\vct v)\}.
\]
Applying Lemma~\ref{lem:exp_tr_WH} with $k=3$ and $k=4$, respectively,
and summing over the basis indices reduces the calculation to
\begin{align*}
\sum_i\mtx H_i(\vct v)^2
&=\frac14\{\|\vct v\|_2^2\mtx I_p+(p+2)\vct v\vct v^\top\},\\
\sum_i\mtx H_i(\vct v)\mtx M\mtx H_i(\vct v)
&=\frac14\{\vct v\vct v^\top\mtx M+\mtx M\vct v\vct v^\top
+\tr(\mtx M)\vct v\vct v^\top
+(\vct v^\top\mtx M\vct v)\mtx I_p\},
\end{align*}
for symmetric $\mtx M$. In the fourth trace moment, label the two repeated
factors involving $\vct\alpha$ as one pair and those involving
$\vct\beta$ as the other. The 24 permutations then form the eight classes
below. For each permutation in a class, the last two columns give the
coefficients of $(\vct\alpha^\top\vct\beta)^2$ and
$\|\vct\alpha\|_2^2\|\vct\beta\|_2^2$ after summing over the two basis
indices; these coefficients are multiplied by the class size and Wishart
weight.
\begin{center}
\small
\begin{tabular}{@{}lcccc@{}}
cycle class & number & Wishart weight
& $(\vct\alpha^\top\vct\beta)^2$
& $\|\vct\alpha\|_2^2\|\vct\beta\|_2^2$ \\[2pt]
identity & $1$ & $n^4$ & $0$ & $1$ \\
within-pair transposition & $2$ & $2n^3$ & $0$ & $(p+1)/2$ \\
across-pair transposition & $4$ & $2n^3$ & $1$ & $0$ \\
three-cycle & $8$ & $4n^2$ & $(p+2)/4$ & $1/4$ \\
within-pair double transposition & $1$ & $4n^2$ & $0$ & $(p+1)^2/4$ \\
cross-pair double transposition & $2$ & $4n^2$ & $(p+2)/4$ & $1/4$ \\
adjacent four-cycle & $4$ & $8n$ & $(p+2)^2/16$ & $(3p+4)/16$ \\
alternating four-cycle & $2$ & $8n$ & $(p+6)/8$ & $1/8$ \\
\end{tabular}
\end{center}
Multiplying the class sizes and weights and collecting the two trace
functionals gives, respectively,
\begin{align*}
&8n^3+10n^2(p+2)+2n(p+2)^2+2n(p+6)\\
&\qquad=2n(4n^2+5np+10n+p^2+5p+10),\\
&n^4+2n^3(p+1)+n^2(p+1)^2+10n^2+2n(3p+4)+2n\\
&\qquad=n^4+2n^3p+2n^3+n^2p^2+2n^2p+11n^2+6np+10n.
\end{align*}
For the mixed $\mtx W,\mtx W^2$ moment, abbreviate
$a=\|\vct\alpha\|_2^2\|\vct\beta\|_2^2$ and
$b=(\vct\alpha^\top\vct\beta)^2$. Apply
\eqref{eq:exp_tr_WH} to
$\vct\alpha\vct\alpha^\top,\mtx H_i(\vct\beta),
\mtx H_i(\vct\beta)$ and sum over $i$. The identity, three
transpositions and two three-cycles give, respectively,
\begin{align*}
\E\{(\vct\alpha^\top\mtx W\vct\alpha)
      (\vct\beta^\top\mtx W^2\vct\beta)\}
&=n^3a+2n^2\left\{\frac{p+1}{2}a+2b\right\}
  +8n\left\{\frac{a+(p+2)b}{4}\right\}\\
&=n(n^2+np+n+2)a+2n(2n+p+2)b.
\end{align*}
Here the three sums used are
$\sum_i\tr(\mtx H_i(\vct\beta))^2=\|\vct\beta\|_2^2$,
$\sum_i\tr\{\mtx H_i(\vct\beta)^2\}=(p+1)\|\vct\beta\|_2^2/2$,
and the displayed formula for $\sum_i\mtx H_i(\vct\beta)^2$.
The cross-transposition sum is
$\sum_i\tr\{\vct\alpha\vct\alpha^\top\mtx H_i(\vct\beta)\}
\tr\{\mtx H_i(\vct\beta)\}=b$, by orthonormality of the basis.

For the third matrix moment, write
$\mtx W=\sum_{r=1}^n\vct z_r\vct z_r^\top$, with independent
$\vct z_r\sim\mathcal N(\vct0,\mtx I_p)$. In the expansion of
$\mtx W^3$, three distinct row indices contribute $\mtx I_p$ per
ordered triple. Each of the three configurations with exactly two equal
indices contributes $(p+2)\mtx I_p$, using
$\E(\|\vct z_r\|_2^2\vct z_r\vct z_r^\top)=(p+2)\mtx I_p$.
When all three indices coincide, the contribution is
$(p+2)(p+4)\mtx I_p$. Thus
\begin{align*}
\E(\mtx W^3)
&=\{n(n-1)(n-2)+3n(n-1)(p+2)+n(p+2)(p+4)\}\mtx I_p\\
&=n\{n^2+3np+3n+p^2+3p+4\}\mtx I_p.
\end{align*}
The two single-row identities follow by rotational invariance and the
fourth and sixth moments of $\|\vct z_r\|_2$. This also gives
$n(n+2)(n+4)$ when $p=1$.
\end{proof}

Allowing $\vct\alpha\ne\vct\beta$ retains the cross-response terms needed
for multivariate inference and extends the corresponding calculation in
Proposition~S1 of \citet{dicker2014variance}. We now use these identities
to obtain the covariance of the variance-component estimators. The
$O(n^{-1})$ remainder below will be negligible relative to the lower
variance scale $q^{-1}$ because $q=o(n)$.

\begin{lemma}[Fixed-effects covariance and scale]
\label{lem:fixed_variance}
Under the conditions of Theorem~\ref{thm:asymptotics_fixed_mt},
$\E(\hat\sigma^2)=\sigma^2$ and $\E(\hat\rho^2)=\rho^2$.
Moreover,
\begin{equation}
\label{eq:fe_uniform_cov_expansion}
\left\|n\,\Cov\!\left(
\begin{bmatrix}\hat\sigma^2\\ \hat\rho^2\end{bmatrix}
\right)-\mtx V_{\mathrm{fe}}\right\|
\le \frac{C}{n},
\end{equation}
where $\mtx V_{\mathrm{fe}}$ is the leading covariance matrix
\begin{equation}
\label{eq:fe_V_matrix_decomp}
\begin{aligned}
\mtx V_{\mathrm{fe}}
&=\frac{2n^2}{q^2(n+1)^2}
\Bigg\{
\|\mtx\Sigma_b\|_F^2
\begin{pmatrix}1+p/n&-(2+p/n)\\-(2+p/n)&4+p/n\end{pmatrix}\\
&\qquad+\tr(\mtx\Sigma_e\mtx\Sigma_b)
\begin{pmatrix}2p/n&-2p/n\\-2p/n&2+2p/n\end{pmatrix} +\tr(\mtx\Sigma_e^2)
\begin{pmatrix}1+p/n&-p/n\\-p/n&p/n\end{pmatrix}
\Bigg\}.
\end{aligned}
\end{equation}
The remainder bound is uniform over all parameter sequences satisfying the
conditions of Theorem~\ref{thm:asymptotics_fixed_mt}. Consequently, for all
sufficiently large $n$,
\begin{equation}
\label{eq:fe_exact_and_leading_scale}
\frac{c}{q}\mtx I_2
\preceq \mtx V_{\mathrm{fe}}
\preceq \frac{C\ell_n}{q}\mtx I_2,
\qquad
\frac{c}{q}\mtx I_2
\preceq n\,\Cov\!\left(
\begin{bmatrix}\hat\sigma^2\\ \hat\rho^2\end{bmatrix}\right)
\preceq \frac{C\ell_n}{q}\mtx I_2.
\end{equation}
\end{lemma}

\begin{proof}
\proofblock{Conditional moments and averaging over the design.}
We first compute the covariance of the two scalar response moments, then
apply the linear transformation that defines the estimators. Recall
\[
\mtx W=\mtx X^\top\mtx X.
\]
Let
\[
T_1=\frac1n\tr(\mtx Y^\top\mtx Y),
\qquad
T_2=\frac1{n^2}\tr(\mtx Y^\top\mtx X\mtx X^\top\mtx Y).
\]
Conditional on $\mtx X$, $\vecop(\mtx Y)$ is Gaussian with mean
$\vecop(\mtx X\mtx B)$ and covariance $\mtx\Sigma_e\otimes\mtx I_n$.
For a Gaussian vector $\vct z$ with mean $\vct\mu$ and symmetric matrices
$\mtx A,\mtx C$,
\begin{align}
\Cov(\vct z^\top\mtx A\vct z,\vct z^\top\mtx C\vct z)
=2\tr\{\mtx A\Cov(\vct z)\mtx C\Cov(\vct z)\}
 +4\vct\mu^\top\mtx A\Cov(\vct z)\mtx C\vct\mu.
\label{eq:gaussian_quad_cov_fe}
\end{align}
Applying \eqref{eq:gaussian_quad_cov_fe} with the row-space matrices
$\mtx I_n$ and $\mtx X\mtx X^\top$ gives
\begin{align}
\E(T_1\mid\mtx X)&=\frac1n\tr(\mtx W\mtx B\mtx B^\top)
+\tr(\mtx\Sigma_e),\notag\\
\E(T_2\mid\mtx X)&=\frac1{n^2}\{\tr(\mtx W^2\mtx B\mtx B^\top)
+\tr(\mtx\Sigma_e)\tr(\mtx W)\},
\label{eq:fe_conditional_means}\\
\Var(T_1\mid\mtx X)&=\frac2n\tr(\mtx\Sigma_e^2)
+\frac4{n^2}\tr(\mtx W\mtx B\mtx\Sigma_e\mtx B^\top),\notag\\
\Cov(T_1,T_2\mid\mtx X)&=\frac2{n^3}\tr(\mtx\Sigma_e^2)\tr(\mtx W)
+\frac4{n^3}\tr(\mtx W^2\mtx B\mtx\Sigma_e\mtx B^\top),\notag\\
\Var(T_2\mid\mtx X)&=\frac2{n^4}\tr(\mtx\Sigma_e^2)\tr(\mtx W^2)
+\frac4{n^4}\tr(\mtx W^3\mtx B\mtx\Sigma_e\mtx B^\top).
\label{eq:fe_conditional_covariances}
\end{align}
After averaging over the design, the mixed terms involve
$\tr(\mtx B\mtx\Sigma_e\mtx B^\top)
=\tr(\mtx\Sigma_e\mtx\Sigma_b)$. The following Wishart identities
supply both the average conditional covariance and the covariance of the
conditional means:
\begin{align}
\E\tr(\mtx W\mtx B\mtx B^\top)&=n\tr(\mtx\Sigma_b),\notag\\
\E\tr(\mtx W^2\mtx B\mtx B^\top)&=n(n+p+1)\tr(\mtx\Sigma_b),\notag\\
\Var\{\tr(\mtx W\mtx B\mtx B^\top)\}&=2n\|\mtx\Sigma_b\|_F^2,\notag\\
\Cov\{\tr(\mtx W\mtx B\mtx B^\top),\tr(\mtx W)\}
&=2n\tr(\mtx\Sigma_b),\notag\\
\Cov\{\tr(\mtx W^2\mtx B\mtx B^\top),\tr(\mtx W)\}
&=4n(n+p+1)\tr(\mtx\Sigma_b),\notag\\
\Var\{\tr(\mtx W)\}&=2np,
\label{eq:fe_wishart_table_1}
\end{align}
\begin{align}
\Cov\{\tr(\mtx W\mtx B\mtx B^\top),
\tr(\mtx W^2\mtx B\mtx B^\top)\}
&=2n(2n+p+2)\|\mtx\Sigma_b\|_F^2
+2n\tr^2(\mtx\Sigma_b),\notag\\
\Var\{\tr(\mtx W^2\mtx B\mtx B^\top)\}
&=2n(4n^2+5np+10n+p^2+5p+10)\|\mtx\Sigma_b\|_F^2\notag\\
&\quad +(10n^2+6np+10n)\tr^2(\mtx\Sigma_b),
\label{eq:fe_wishart_table_2}
\end{align}
and
\begin{align}
\E\tr(\mtx W)&=np,\notag\\
\E\tr(\mtx W^2)&=np(n+p+1),\notag\\
\E\tr(\mtx W\mtx B\mtx\Sigma_e\mtx B^\top)
&=n\tr(\mtx\Sigma_e\mtx\Sigma_b),\notag\\
\E\tr(\mtx W^2\mtx B\mtx\Sigma_e\mtx B^\top)
&=n(n+p+1)\tr(\mtx\Sigma_e\mtx\Sigma_b),\notag\\
\E\tr(\mtx W^3\mtx B\mtx\Sigma_e\mtx B^\top)
&=n\{n^2+3np+3n+p^2+3p+4\}
\tr(\mtx\Sigma_e\mtx\Sigma_b).
\label{eq:fe_wishart_table_3}
\end{align}
The trace moments involving $\mtx B\mtx B^\top$ follow by diagonalizing
that matrix and applying Lemma~\ref{lem:fixed_wishart_identities}; the
identities involving $\tr(\mtx W)$ also follow from
Lemma~\ref{lem:exp_tr_WH}. The last line of
\eqref{eq:fe_wishart_table_3} uses \eqref{eq:wishart_third_matrix} and
linearity of the trace.

The law of total covariance now combines two contributions: the average
conditional covariance in \eqref{eq:fe_conditional_covariances} and the
covariance of the conditional means in \eqref{eq:fe_conditional_means}.
Substituting the preceding identities gives
\begin{align}
\Var(T_1)
&=\frac{2}{n}\{\|\mtx\Sigma_b\|_F^2
+2\tr(\mtx\Sigma_e\mtx\Sigma_b)+\tr(\mtx\Sigma_e^2)\},
\label{eq:fe_T_var1}\\
\Cov(T_1,T_2)
&=\frac{2}{n^2}\Bigl\{(2n+p+2)\|\mtx\Sigma_b\|_F^2
+\tr^2(\mtx\Sigma_b)\notag\\
&\qquad+2(n+p+1)\tr(\mtx\Sigma_e\mtx\Sigma_b)
+\tr(\mtx\Sigma_b)\tr(\mtx\Sigma_e)
+p\tr(\mtx\Sigma_e^2)\Bigr\},
\label{eq:fe_T_cov12}\\
\Var(T_2)
&=\frac{2}{n^3}\Bigl\{
 (4n^2+5np+10n+p^2+5p+10)\|\mtx\Sigma_b\|_F^2
 +(5n+3p+5)\tr^2(\mtx\Sigma_b)\notag\\
&\qquad +2(p^2+3np+3p+n^2+3n+4)
\tr(\mtx\Sigma_e\mtx\Sigma_b)\notag\\
&\qquad +4(n+p+1)\tr(\mtx\Sigma_b)\tr(\mtx\Sigma_e)
+p(n+p+1)\tr(\mtx\Sigma_e^2)+p\tr^2(\mtx\Sigma_e)
\Bigr\}.
\label{eq:fe_T_var2}
\end{align}
For example, the pure-noise contribution to \eqref{eq:fe_T_cov12} is
$2\tr(\mtx\Sigma_e^2)\E\tr(\mtx W)/n^3
=2p\tr(\mtx\Sigma_e^2)/n^2$. It involves the full noise covariance,
including its off-diagonal entries.

\proofblock{Leading covariance and remainder.}
The two moment estimators satisfy
\[
\begin{bmatrix}\hat\sigma^2\\ \hat\rho^2\end{bmatrix}
=\frac1{q(n+1)}
\begin{pmatrix}p+n+1&-n\\-p&n\end{pmatrix}
\begin{bmatrix}T_1\\T_2\end{bmatrix}.
\]
The means in \eqref{eq:fe_conditional_means} and the first moments in
\eqref{eq:fe_wishart_table_1} give unbiasedness.  Substituting
\eqref{eq:fe_T_var1}--\eqref{eq:fe_T_var2} into the covariance of the
preceding linear transformation and collecting terms gives
\[
n\,\Cov(\widehat{\vct\theta})=\mtx V_{\mathrm{fe}}+\mtx R_n,
\qquad
\widehat{\vct\theta}=(\hat\sigma^2,\hat\rho^2)^\top,
\]
where $\mtx V_{\mathrm{fe}}$ is exactly
\eqref{eq:fe_V_matrix_decomp}. The remainder $\mtx R_n$ consists of the
lower-order coefficients multiplying the same six trace functionals
appearing in \eqref{eq:fe_T_var1}--\eqref{eq:fe_T_var2}.

The remainder can be checked entrywise without expanding the leading
matrix again. Extract the common factor $2\{q^2(n+1)^2\}^{-1}$.
For each trace functional $F$ in the following table, let
$\mtx R_F$ be the symmetric matrix whose three entries are listed. Then
\[
\mtx R_n=\frac{2}{q^2(n+1)^2}\sum_F F\mtx R_F.
\]
These are the exact omitted coefficients obtained from the preceding
linear transformation:
\begin{center}
\small
\begin{tabular}{@{}lccc@{}}
trace functional $F$ & $(\mtx R_F)_{11}$ & $(\mtx R_F)_{12}$
 & $(\mtx R_F)_{22}$ \\[2pt]
$\|\mtx\Sigma_b\|_F^2$ & $4n+p+7$ & $-6n-p-8$ & $10n+p+10$ \\
$\tr(\mtx\Sigma_e\mtx\Sigma_b)$ & $2n+2p+6$ & $-2n-2p-6$ & $6n+2p+8$ \\
$\tr(\mtx\Sigma_e^2)$ & $2n+p+1$ & $-p$ & $p$ \\
$\tr^2(\mtx\Sigma_b)$ & $3n+p+3$ & $-4n-p-4$ & $5n+p+5$ \\
$\tr(\mtx\Sigma_b)\tr(\mtx\Sigma_e)$ & $2n+2p+2$ & $-3n-2p-3$ & $4n+2p+4$ \\
$\tr^2(\mtx\Sigma_e)$ & $p$ & $-p$ & $p$ \\
\end{tabular}
\end{center}
Every entry is $O(n+p)=O(n)$ because $p/n=O(1)$.
Under the conditions of
Theorem~\ref{thm:asymptotics_fixed_mt},
\[
\|\mtx\Sigma_b\|_F^2+\tr(\mtx\Sigma_e\mtx\Sigma_b)
+\tr(\mtx\Sigma_e^2)=O(\ell_n q),
\qquad
\tr^2(\mtx\Sigma_b)+\tr(\mtx\Sigma_b)\tr(\mtx\Sigma_e)
+\tr^2(\mtx\Sigma_e)=O(q^2).
\]
After multiplication by the common factor $q^{-2}(n+1)^{-2}$ in
$n\Cov(\hat\sigma^2,\hat\rho^2)$, the two classes contribute respectively
$O\{\ell_n/(nq)\}$ and $O(n^{-1})$. Since
$\ell_n\le Cq$ by the trace upper bounds, the total is still $O(n^{-1})$.  This proves the asserted $C/n$ bound.

\proofblock{Covariance scale.}
We establish the lower bound before using the remainder estimate. The trace
assumptions provide signal and noise contributions of order at least $q$,
and their coefficient matrices together are uniformly positive definite,
even when $p/n$ tends to zero. For every positive semidefinite $q\times q$
matrix $\mtx A$,
\[
\frac{\tr^2(\mtx A)}{q}\le\|\mtx A\|_F^2
\le\|\mtx A\|\tr(\mtx A).
\]
The trace bounds in Theorem~\ref{thm:asymptotics_fixed_mt} give
\[
\|\mtx\Sigma_b\|_F^2\ge cq,
\qquad
\tr(\mtx\Sigma_e^2)\ge cq,
\]
while positive semidefiniteness and the definition of $\ell_n$ give
\[
\|\mtx\Sigma_b\|_F^2+\tr(\mtx\Sigma_e^2)
+\tr(\mtx\Sigma_e\mtx\Sigma_b)\le C\ell_n q.
\]
For $t=p/n$, denote the first and third coefficient matrices in
\eqref{eq:fe_V_matrix_decomp} by $\mtx P_b(t)$ and $\mtx P_e(t)$.
Both are positive semidefinite, as is the middle coefficient matrix.
Set $b=q^{-1}\|\mtx\Sigma_b\|_F^2$ and
$e=q^{-1}\tr(\mtx\Sigma_e^2)$. Since $b,e\ge c_0>0$,
\[
b\mtx P_b(t)+e\mtx P_e(t)
\succeq c_0\{\mtx P_b(t)+\mtx P_e(t)\}
\succeq c_0
\begin{pmatrix}2&-2\\-2&4\end{pmatrix}
\succeq c\mtx I_2.
\]
The second inequality follows because the difference is
$2t\begin{pmatrix}1&-1\\-1&1\end{pmatrix}\succeq0$.
This proves the lower bound even when $b$ and $e$ diverge. For the upper
bound, all coefficient matrices have bounded norm because $p/n=O(1)$,
and all three trace functionals are at most $C\ell_n q$.
Substitution in \eqref{eq:fe_V_matrix_decomp} gives the first pair of
bounds in \eqref{eq:fe_exact_and_leading_scale}.  The uniform remainder
bound gives
\[
\left\|n\Cov\!\left(
\begin{bmatrix}\hat\sigma^2\\ \hat\rho^2\end{bmatrix}\right)
-\mtx V_{\mathrm{fe}}\right\|\le C/n.
\]
Since $q=o(n)$, this perturbation is $o(q^{-1})$, and Weyl's inequality
proves the second pair of bounds.
\end{proof}

To connect this covariance calculation to the standard error in the main
paper, let $\widehat{\mtx V}_{\mathrm{fe}}$ be the plug-in version of
\eqref{eq:fe_V_matrix_decomp}, obtained by replacing
$(\mtx\Sigma_b,\mtx\Sigma_e)$ with
$(\widehat{\mtx\Sigma}_b,\widehat{\mtx\Sigma}_e)$. With the gradients
$\vct d$ and $\hat{\vct d}$ defined in \eqref{eq:snr_gradient}, set
\begin{equation}
\label{eq:sigma_r_fixed_mt}
\sigma_{r,\mathrm{fe}}^2=\vct d^\top\mtx V_{\mathrm{fe}}\vct d,
\qquad
\hat\sigma_{r,\mathrm{fe}}^2
=\hat{\vct d}^\top\widehat{\mtx V}_{\mathrm{fe}}\hat{\vct d}.
\end{equation}
Expanding the first quadratic form using
\eqref{eq:fe_V_matrix_decomp} gives
$\mathcal V_{\mathrm{fe},n}(\mtx\Sigma_b,\mtx\Sigma_e;\rho^2,\sigma^2)$
in \eqref{eq:fe_scalar_variance_main}; the plug-in version is obtained by
evaluating the same function at the estimated covariance matrices and scalar
components. Since
$\rho^2,\sigma^2\in[c,C]$, $p/n=O(1)$ and
$\tau=\rho^2+\sigma^2\asymp1$, all three scalar coefficients multiplying
$\|\mtx\Sigma_b\|_F^2$,
$\tr(\mtx\Sigma_e\mtx\Sigma_b)$ and
$\tr(\mtx\Sigma_e^2)$ in
\eqref{eq:fe_scalar_variance_main} are bounded above and away from zero.
Hence
\begin{equation}
\label{eq:fe_effective_variance_scale}
\sigma_{r,\mathrm{fe}}^2
\asymp
\frac{\|\mtx\Sigma_b\|_F^2
+\tr(\mtx\Sigma_e\mtx\Sigma_b)
+\tr(\mtx\Sigma_e^2)}{q^2}.
\end{equation}

\subsection{Normal approximation, plug-in variance, and completion}
\label{sec:proof_fe_chaos_clt}
\label{sec:proof_plugin_trace_fe}
The unconditional moments are degree-two and degree-four polynomials in
Gaussian variables, rather than a single quadratic form. We therefore use
the second-order Gaussian Poincar\'e inequality for the normal
approximation. Its first-derivative calculation will also give the
Frobenius error bound needed for the plug-in variance. Write
$d_{\mathrm{TV}}$ for total variation distance.

\begin{lemma}[Normal approximation for fixed effects]
\label{lem:fe_chaos_clt}
Under the conditions of Theorem~\ref{thm:asymptotics_fixed_mt}, uniformly
over deterministic
$\vct a\in\mathbb R^2$ with $\|\vct a\|_2=1$,
\begin{equation}
\label{eq:fe_chaos_bound}
d_{\mathrm{TV}}\!\left(
\frac{\sqrt n\,\vct a^\top
\begin{pmatrix}\hat\sigma^2-\sigma^2\\
\hat\rho^2-\rho^2\end{pmatrix}}
{\left[n\Var\!\left\{\vct a^\top
\begin{pmatrix}\hat\sigma^2-\sigma^2\\
\hat\rho^2-\rho^2\end{pmatrix}\right\}\right]^{1/2}},
\mathcal N(0,1)
\right)
\le \frac{C\ell_n}{\sqrt n}=o(1).
\end{equation}
\end{lemma}

\begin{proof}
\proofblock{Gaussian polynomial representation.}
Write $\mtx E=\mtx G_e\mtx\Sigma_e^{1/2}$, where $\mtx G_e$ has i.i.d.\ standard Gaussian entries and is independent of $\mtx X$, and set
\[
\mtx Z=[\mtx X,\mtx G_e]\in\mathbb R^{n\times(p+q)}.
\]
Let
\[
\mtx P=\begin{pmatrix}\mtx I_p&0\\0&0\end{pmatrix},
\qquad
\mtx H=
\begin{bmatrix}\mtx B\\ \mtx\Sigma_e^{1/2}\end{bmatrix}
\begin{bmatrix}\mtx B^\top&\mtx\Sigma_e^{1/2}\end{bmatrix}.
\]
Here $\mtx P$ selects the predictor columns, and $\mtx H$ contains the
fixed signal and noise loadings. The first $p$ columns of $\mtx Z$ are
$\mtx X$, while
$\mtx Y=\mtx Z\begin{bmatrix}\mtx B\\
\mtx\Sigma_e^{1/2}\end{bmatrix}$. Hence
\begin{equation}
\label{eq:fe_T1_T2_matrix_form}
T_1=\frac1n\tr(\mtx Z^\top\mtx Z\mtx H),
\qquad
T_2=\frac1{n^2}\tr(\mtx Z^\top\mtx Z\mtx P\mtx Z^\top\mtx Z\mtx H).
\end{equation}
For a unit vector $\vct a$, define
\begin{equation}
\label{eq:fe_F_polynomial}
F_{n,\vct a}:=\sqrt n\,\vct a^\top
\begin{pmatrix}\hat\sigma^2-\sigma^2\\
\hat\rho^2-\rho^2\end{pmatrix}.
\end{equation}
The estimator transformation in Lemma~\ref{lem:fixed_variance} expresses
$F_{n,\vct a}$ as a centered polynomial in the independent Gaussian entries
of $\mtx Z$. For $\vct a=(a_1,a_2)^\top$, in the ordering
$(\hat\sigma^2,\hat\rho^2)$, it is
\begin{equation}
\label{eq:fe_F_coefficients}
F_{n,\vct a}
=\frac{\sqrt n\{(p+n+1)a_1-pa_2\}}{q(n+1)}(T_1-\E T_1)
+\frac{n^{3/2}(-a_1+a_2)}{q(n+1)}(T_2-\E T_2).
\end{equation}
Both coefficients are bounded in absolute value by $C\sqrt n/q$,
uniformly over $\|\vct a\|_2=1$, because $p/n=O(1)$.

\proofblock{Derivative bounds.}
We bound the gradient and Hessian of each moment before applying this
linear combination. Matrices are viewed as Euclidean vectors with the
Frobenius inner product. For a perturbation $\mtx U$,
\[
D T_1[\mtx U]=\frac2n\langle\mtx Z\mtx H,\mtx U\rangle_F,
\qquad
\nabla T_1=\frac2n\mtx Z\mtx H,
\qquad
D^2T_1[\mtx U]=\frac2n\mtx U\mtx H.
\]
For $T_2$, because
$D(\mtx Z^\top\mtx Z)[\mtx U]=\mtx Z^\top\mtx U+\mtx U^\top\mtx Z$,
\begin{align}
\nabla T_2
&=\frac2{n^2}\mtx Z(\mtx P\mtx Z^\top\mtx Z\mtx H
+\mtx H\mtx Z^\top\mtx Z\mtx P),
\label{eq:fe_grad_T2}\\
D^2T_2[\mtx U]
&=\frac2{n^2}\mtx U(\mtx P\mtx Z^\top\mtx Z\mtx H
+\mtx H\mtx Z^\top\mtx Z\mtx P)
+\frac2{n^2}\mtx Z\Bigl\{
\mtx P(\mtx Z^\top\mtx U+\mtx U^\top\mtx Z)\mtx H
+\mtx H(\mtx Z^\top\mtx U+\mtx U^\top\mtx Z)\mtx P
\Bigr\}.
\label{eq:fe_hess_T2}
\end{align}
Since $\|\mtx P\|=1$, these formulas imply
\begin{align}
\|\nabla T_1\|_F&\le \frac2n\|\mtx Z\|\,\|\mtx H\|_F,
&\|\nabla^2T_1\|_{\mathrm{op}}&\le \frac2n\|\mtx H\|,
\label{eq:fe_derivative_bounds_T1}\\
\|\nabla T_2\|_F&\le \frac4{n^2}\|\mtx Z\|^3\|\mtx H\|_F,
&\|\nabla^2T_2\|_{\mathrm{op}}&\le \frac C{n^2}\|\mtx Z\|^2\|\mtx H\|_F.
\label{eq:fe_derivative_bounds_T2}
\end{align}
Here the Hessian operator norm is taken with respect to the Frobenius
norm of the perturbation. Set
$\mtx R=\mtx Z^\top\mtx U+\mtx U^\top\mtx Z$; then
$\|\mtx R\|\le2\|\mtx Z\|\|\mtx U\|_F$. For example, two of the
terms in \eqref{eq:fe_hess_T2} satisfy
\begin{align*}
\|\mtx U\mtx P\mtx Z^\top\mtx Z\mtx H\|_F
&\le\|\mtx U\|_F\|\mtx Z\|^2\|\mtx H\|_F,\\
\|\mtx Z\mtx P\mtx R\mtx H\|_F
&\le\|\mtx Z\|\|\mtx R\|\|\mtx H\|_F
\le2\|\mtx Z\|^2\|\mtx H\|_F\|\mtx U\|_F.
\end{align*}
The terms with $\mtx P$ and $\mtx H$ interchanged obey the same bounds.
Taking the supremum over $\|\mtx U\|_F=1$ proves the Hessian bound without
an extra factor depending on the number of matrix entries.

The response dependence enters these bounds through the deterministic
matrix $\mtx H$, for which
\begin{align}
\|\mtx H\|&\le C\ell_n,\notag\\
\|\mtx H\|_F^2
&=\tr\{(\mtx\Sigma_b+\mtx\Sigma_e)^2\}=O(\ell_n q).
\label{eq:fe_H_bounds}
\end{align}
Moreover, because $p/n$ is bounded and $q=o(n)$, the Gaussian matrix
$\mtx Z$ has $p+q=O(n)$ columns. Applying
\eqref{eq:subg_matrix_norm_moments} with $m=p+q$ and the identity
population covariance gives, for every fixed integer $k\ge1$,
\begin{equation}
\label{eq:fe_Z_norm_moments}
\{\E\|\mtx Z\|^k\}^{1/k}=O(\sqrt n).
\end{equation}
Combining \eqref{eq:fe_F_polynomial}--\eqref{eq:fe_Z_norm_moments} yields, uniformly over unit vectors,
\begin{align}
\{\E\|\nabla F_{n,\vct a}\|_F^4\}^{1/4}
&\le C\sqrt{\frac{\ell_n}{q}},
\label{eq:fe_kappa1}\\
\{\E\|\nabla^2F_{n,\vct a}\|_{\mathrm{op}}^4\}^{1/4}
&\le C\sqrt{\frac{\ell_n}{nq}}.
\label{eq:fe_kappa2}
\end{align}
Specifically, multiplying the two moment-derivative bounds by their
coefficients in \eqref{eq:fe_F_coefficients} gives
\begin{align*}
\{\E\|\nabla F_{n,\vct a}\|_F^4\}^{1/4}
&\le \frac{C\sqrt n}{q}\frac{\sqrt{n\ell_n q}}n
     +\frac{C\sqrt n}{q}\frac{n^{3/2}\sqrt{\ell_n q}}{n^2}
\le C\sqrt{\frac{\ell_n}{q}},\\
\{\E\|\nabla^2F_{n,\vct a}\|_{\mathrm{op}}^4\}^{1/4}
&\le \frac{C\sqrt n}{q}\frac{\ell_n}n
     +\frac{C\sqrt n}{q}\frac{n\sqrt{\ell_n q}}{n^2}
\le C\sqrt{\frac{\ell_n}{nq}}.
\end{align*}
The last inequality uses $\ell_n\le Cq$.

\proofblock{Normal approximation.}
The derivative bounds are compared with the exact variance, not only with
an upper variance bound. Lemma~\ref{lem:fixed_variance} gives
\begin{equation}
\label{eq:fe_exact_variance_lower}
\frac{c}{q}\le
n\Var\!\left\{\vct a^\top
\begin{pmatrix}\hat\sigma^2-\sigma^2\\
\hat\rho^2-\rho^2\end{pmatrix}\right\}
\le\frac{C\ell_n}{q}
\end{equation}
for all sufficiently large $n$, uniformly over $\|\vct a\|_2=1$.

The second-order Gaussian Poincar\'e inequality of
\citet[Theorem~2.2]{chatterjee2009fluctuations}, specialized to independent
standard Gaussian inputs, gives, for a centred twice continuously
differentiable function $F$ with finite fourth moment and variance $s^2>0$,
\[
d_{\mathrm{TV}}(F/s,\mathcal N(0,1))
\le \frac{2\sqrt5}{s^2}
\{\E\|\nabla F\|_2^4\}^{1/4}
\{\E\|\nabla^2F\|_{\mathrm{op}}^4\}^{1/4}.
\]
Here $F_{n,\vct a}$ is a polynomial in Gaussian variables, so its fourth
moment is finite. Substituting
\eqref{eq:fe_kappa1}--\eqref{eq:fe_exact_variance_lower} proves
\[
d_{\mathrm{TV}}\!\left(
\frac{F_{n,\vct a}}
{\left[n\Var\!\left\{\vct a^\top
\begin{pmatrix}\hat\sigma^2-\sigma^2\\
\hat\rho^2-\rho^2\end{pmatrix}\right\}\right]^{1/2}},
\mathcal N(0,1)\right)
\le C\ell_n n^{-1/2}=o(1),
\]
which is \eqref{eq:fe_chaos_bound}.
\end{proof}

\proofblock{Plug-in variance estimate.}
The leading covariance depends on the response matrices only through their
two squared Frobenius norms and their mixed trace. We first bound the
Frobenius errors of the covariance estimators, then use one product
inequality to control all three functionals. No separate operator-norm
consistency argument is needed here.

\begin{lemma}[Accuracy of the fixed-effects variance estimate]
\label{lem:plugin_trace_fe}
Under the conditions of Theorem~\ref{thm:asymptotics_fixed_mt},
\begin{equation}
\label{eq:fe_matrix_mse}
q^{-2}\|\widehat{\mtx\Sigma}_b-\mtx\Sigma_b\|_F^2
+q^{-2}\|\widehat{\mtx\Sigma}_e-\mtx\Sigma_e\|_F^2
=O_P(n^{-1}).
\end{equation}
Moreover,
\begin{align}
\label{eq:fe_trace_relative_rates}
q^{-2}\bigl|\|\widehat{\mtx\Sigma}_b\|_F^2-\|\mtx\Sigma_b\|_F^2\bigr|
&=O_P\!\left(\sqrt{\frac{\ell_n}{nq}}+n^{-1}\right),\\
q^{-2}\bigl|\tr(\widehat{\mtx\Sigma}_e\widehat{\mtx\Sigma}_b)
-\tr(\mtx\Sigma_e\mtx\Sigma_b)\bigr|
&=O_P\!\left(\sqrt{\frac{\ell_n}{nq}}+n^{-1}\right),\notag\\
q^{-2}\bigl|\tr(\widehat{\mtx\Sigma}_e^2)-\tr(\mtx\Sigma_e^2)\bigr|
&=O_P\!\left(\sqrt{\frac{\ell_n}{nq}}+n^{-1}\right).\notag
\end{align}
Consequently,
\begin{equation}
\label{eq:fe_V_relative_rate}
\|\widehat{\mtx V}_{\mathrm{fe}}-\mtx V_{\mathrm{fe}}\|
=O_P\!\left(\sqrt{\frac{\ell_n}{nq}}+n^{-1}\right)
=o_P(q^{-1}),
\end{equation}
\end{lemma}

\begin{proof}
We reuse the first-derivative bounds from the normal-approximation proof,
now applied to individual entries of the two matrix moments. Use the Gaussian matrix $\mtx Z=[\mtx X,\mtx G_e]$ and projector
$\mtx P$ from the proof of Lemma~\ref{lem:fe_chaos_clt}, and let $\vct v_k$
be the $k$th column of
$[\mtx B^\top,\mtx\Sigma_e^{1/2}]^\top$, so that
$\mtx Y_{\cdot k}=\mtx Z\vct v_k$. Put
\[
\mtx H_{kl}=\frac12(\vct v_k\vct v_l^\top+\vct v_l\vct v_k^\top).
\]
Then
\[
\|\mtx H_{kl}\|_F^2
=\frac12\{\|\vct v_k\|_2^2\|\vct v_l\|_2^2
 +(\vct v_k^\top\vct v_l)^2\}
\le\|\vct v_k\|_2^2\|\vct v_l\|_2^2,
\qquad
\|\vct v_k\|_2^2=(\mtx\Sigma_b+\mtx\Sigma_e)_{kk}.
\]
The two sample-moment entries have the same polynomial forms as
\eqref{eq:fe_T1_T2_matrix_form}, with $\mtx H$ replaced by $\mtx H_{kl}$:
\begin{align*}
\frac1n\mtx Y_{\cdot k}^\top\mtx Y_{\cdot l}
&=\frac1n\tr(\mtx Z^\top\mtx Z\mtx H_{kl}),\\
\frac1{n^2}\mtx Y_{\cdot k}^\top\mtx X\mtx X^\top\mtx Y_{\cdot l}
&=\frac1{n^2}\tr(\mtx Z^\top\mtx Z\mtx P
                       \mtx Z^\top\mtx Z\mtx H_{kl}).
\end{align*}
The derivative calculations
\eqref{eq:fe_derivative_bounds_T1}--\eqref{eq:fe_derivative_bounds_T2}
hold for every symmetric $\mtx H_{kl}$, without a positivity requirement.
Thus Gaussian Poincar\'e and \eqref{eq:fe_Z_norm_moments} give
\begin{align*}
\Var\!\left(\frac1n\mtx Y_{\cdot k}^\top\mtx Y_{\cdot l}\right)
&\le\frac4{n^2}\E\|\mtx Z\|^2\,\|\mtx H_{kl}\|_F^2
\le\frac Cn(\mtx\Sigma_b+\mtx\Sigma_e)_{kk}
              (\mtx\Sigma_b+\mtx\Sigma_e)_{ll},\\
\Var\!\left(\frac1{n^2}\mtx Y_{\cdot k}^\top\mtx X\mtx X^\top
                  \mtx Y_{\cdot l}\right)
&\le\frac{16}{n^4}\E\|\mtx Z\|^6\,\|\mtx H_{kl}\|_F^2
\le\frac Cn(\mtx\Sigma_b+\mtx\Sigma_e)_{kk}
              (\mtx\Sigma_b+\mtx\Sigma_e)_{ll}.
\end{align*}
The matrix expectations used for centering follow in every response
direction: for $\vct u\in\mathbb R^q$, apply the scalar moment calculation
to $\mtx Y\vct u=\mtx X(\mtx B\vct u)+\mtx E\vct u$, whose signal
and noise variances are $\vct u^\top\mtx\Sigma_b\vct u$ and
$\vct u^\top\mtx\Sigma_e\vct u$. Equality of these quadratic forms
for every $\vct u$ gives the two matrix expectations below. Summing the
entrywise variances then gives the expected squared Frobenius error;
no covariance between different entries is needed:

\begin{align}
\label{eq:fe_sample_moment_mse}
&\E\left\|\frac1n\mtx Y^\top\mtx Y
 -(\mtx\Sigma_b+\mtx\Sigma_e)\right\|_F^2
+\E\left\|\frac1{n^2}\mtx Y^\top\mtx X\mtx X^\top\mtx Y
 -\left(\frac{p+n+1}{n}\mtx\Sigma_b+\frac pn\mtx\Sigma_e\right)
 \right\|_F^2\notag\\
&\qquad\le\frac Cn\sum_{k,l}
 (\mtx\Sigma_b+\mtx\Sigma_e)_{kk}
 (\mtx\Sigma_b+\mtx\Sigma_e)_{ll}
=\frac Cn\{\tr(\mtx\Sigma_b)+\tr(\mtx\Sigma_e)\}^2
=O(q^2/n).
\end{align}
The estimators satisfy the exact centered representations
\begin{align*}
\widehat{\mtx\Sigma}_b-\mtx\Sigma_b
&=-\frac p{n+1}
\left\{\frac1n\mtx Y^\top\mtx Y-(\mtx\Sigma_b+\mtx\Sigma_e)\right\}\\
&\quad+\frac n{n+1}
\left\{\frac1{n^2}\mtx Y^\top\mtx X\mtx X^\top\mtx Y
 -\left(\frac{p+n+1}{n}\mtx\Sigma_b+\frac pn\mtx\Sigma_e\right)\right\},\\
\widehat{\mtx\Sigma}_e-\mtx\Sigma_e
&=\frac{p+n+1}{n+1}
\left\{\frac1n\mtx Y^\top\mtx Y-(\mtx\Sigma_b+\mtx\Sigma_e)\right\}\\
&\quad-\frac n{n+1}
\left\{\frac1{n^2}\mtx Y^\top\mtx X\mtx X^\top\mtx Y
 -\left(\frac{p+n+1}{n}\mtx\Sigma_b+\frac pn\mtx\Sigma_e\right)\right\}.
\end{align*}
Their coefficients are bounded because $p/n=O(1)$. Hence
\eqref{eq:fe_sample_moment_mse} and Markov's inequality prove
\eqref{eq:fe_matrix_mse}.

It remains to pass from matrix errors to the three covariance functionals.
For $a\in\{b,e\}$,
write $\mtx\Delta_a=\widehat{\mtx\Sigma}_a-\mtx\Sigma_a$. The preceding
Frobenius rate and the population trace bounds give
\[
\|\mtx\Delta_a\|_F=O_P(q/\sqrt n),
\qquad
\|\mtx\Sigma_a\|_F=O(\sqrt{\ell_nq}).
\]
For any $a,c\in\{b,e\}$, expand the product and use Cauchy--Schwarz:
\begin{align}
\label{eq:fe_trace_product_bound}
\frac1{q^2}\left|
\tr(\widehat{\mtx\Sigma}_a\widehat{\mtx\Sigma}_c)
-\tr(\mtx\Sigma_a\mtx\Sigma_c)\right|
&\le\frac{
\|\mtx\Sigma_a\|_F\|\mtx\Delta_c\|_F
+\|\mtx\Sigma_c\|_F\|\mtx\Delta_a\|_F
+\|\mtx\Delta_a\|_F\|\mtx\Delta_c\|_F}{q^2}\notag\\
&=O_P\!\left(\sqrt{\frac{\ell_n}{nq}}+\frac1n\right).
\end{align}
Taking $(a,c)=(b,b),(e,b),(e,e)$ proves
\eqref{eq:fe_trace_relative_rates}. Since the design coefficients in
\eqref{eq:fe_V_matrix_decomp} are bounded,
\[
\|\widehat{\mtx V}_{\mathrm{fe}}-\mtx V_{\mathrm{fe}}\|
=O_P\!\left(\sqrt{\frac{\ell_n}{nq}}+\frac1n\right)
=o_P(q^{-1}).
\]
The last step compares the error with $q^{-1}$ and uses
$\sqrt{\ell_nq/n}+q/n\to0$. This proves
\eqref{eq:fe_V_relative_rate}.
\end{proof}

\begin{proof}[Proof of Theorem~\ref{thm:asymptotics_fixed_mt}]
We verify the scalar normal approximation and total-variance consistency
required by Lemma~\ref{lem:common_studentization}. Write $\vct\Delta=(\hat\sigma^2-\sigma^2,\hat\rho^2-\rho^2)^\top$.
Apply Lemma~\ref{lem:fe_chaos_clt} in the deterministic direction
$\vct a=\vct d/\|\vct d\|_2$. Lemma~\ref{lem:fixed_variance} gives
\[
\left|\frac{n\Var(\vct d^\top\vct\Delta)}
{\vct d^\top\mtx V_{\mathrm{fe}}\vct d}-1\right|
\le\frac{(C/n)\|\vct d\|_2^2}{(c/q)\|\vct d\|_2^2}
=O(q/n)=o(1),
\]
so the exact-variance normal approximation implies
\begin{equation}
\label{eq:fe_linear_snr_clt}
\frac{\sqrt n\,\vct d^\top\vct\Delta}
{(\vct d^\top\mtx V_{\mathrm{fe}}\vct d)^{1/2}}
\Longrightarrow\mathcal N(0,1).
\end{equation}
For the total variance, $\hat\tau=\tr(\mtx Y^\top\mtx Y)/(nq)$, and
\eqref{eq:fe_T_var1} gives the exact identities and bound
\[
\E\hat\tau=\tau,\qquad
\Var(\hat\tau)=\frac{2}{nq^2}
\tr\{(\mtx\Sigma_b+\mtx\Sigma_e)^2\}
\le\frac{C\ell_n}{nq}.
\]
Since $\tau\asymp1$ and $\ell_n/(nq)\to0$, Chebyshev's inequality gives
$\hat\tau/\tau\to_P1$. The scale bound
\eqref{eq:fe_exact_and_leading_scale} and the plug-in bound
\eqref{eq:fe_V_relative_rate} now verify
Lemma~\ref{lem:common_studentization} with
$\mtx V_n=\mtx V_{\mathrm{fe}}$ and
$\widehat{\mtx V}_n=\widehat{\mtx V}_{\mathrm{fe}}$.
Together with \eqref{eq:sigma_r_fixed_mt}, the lemma proves the
studentized limit, ratio consistency and positivity of the estimated
variance, and the stated coverage.
\end{proof}

\section{Proof of Theorem~\ref{thm:asymptotics_random_fixed_homo_mt}}
\label{sec:supp_fixed_design_theorem}
With the design fixed, the response vector is jointly Gaussian and the two
moment estimators are exactly conditionally unbiased. Their covariance and
normal approximation therefore follow from Gaussian quadratic forms,
without the Wishart averaging needed in S2. We retain the normalized Gram
matrix used in the main paper throughout the argument.

\subsection{Conditional moments, normal approximation, and studentization}
\label{sec:proof_random_fixed_homo_mt}

Throughout S3, probabilities and expectations are conditional on the
supplied deterministic design. Write
\[
\mtx H=\frac{n\mtx X\mtx X^\top}{\|\mtx X\|_F^2},
\qquad
\mtx\Sigma_b=\frac{\|\mtx X\|_F^2}{np}\mtx\Sigma_\beta,
\qquad
g_k=\frac1n\tr(\mtx H^k),\quad g_0=g_1=1.
\]
Thus $g_2,g_3,g_4$ are exactly the main-text spectral summaries. Conditional
on $\mtx X$, the response vector is centered Gaussian with covariance
\begin{equation}
\label{eq:fixed_homo_direct_covariance}
\mtx\Omega_X
=\frac1p\mtx\Sigma_\beta\otimes\mtx X\mtx X^\top
 +\mtx\Sigma_e\otimes\mtx I_n
=\mtx\Sigma_b\otimes\mtx H+\mtx\Sigma_e\otimes\mtx I_n.
\end{equation}
Put $\ell_n=1\vee\|\mtx\Sigma_b\|\vee\|\mtx\Sigma_e\|$.
The theorem gives $\|\mtx H\|\le C$, $g_2-1\ge c$ and $g_k\le C$ for
$k\le4$. The trace bounds give $\ell_n\le Cq$, and the response-spectrum
condition is $\ell_nq/n\to0$.

Define the two matrix moments by
\[
\mtx C_j=\mtx H^j,\qquad
\mtx M_j=\frac1n\mtx Y^\top\mtx C_j\mtx Y,\qquad j=0,1.
\]
Their conditional moment system and its inverse are
\begin{align}
\E\left[
\begin{pmatrix}\mtx M_0\\\mtx M_1\end{pmatrix}
\middle|\mtx X\right]
&=\begin{pmatrix}1&1\\1&g_2\end{pmatrix}
\begin{pmatrix}\mtx\Sigma_e\\\mtx\Sigma_b\end{pmatrix},
\label{eq:random_conditional_moment_system}\\
\widehat{\mtx A}_n
&=\frac1{g_2-1}\begin{pmatrix}g_2&-1\\-1&1\end{pmatrix}.
\label{eq:random_moment_inverse}
\end{align}
In particular, $\|\widehat{\mtx A}_n\|+
\|\widehat{\mtx A}_n^{-1}\|\le C$, and the component estimators in
\eqref{eq:homo_estimators_mt} satisfy
\begin{equation}
\label{eq:fixed_homo_matrix_moment_system}
\begin{pmatrix}\widehat{\mtx\Sigma}_e\\\widehat{\mtx\Sigma}_b\end{pmatrix}
=\widehat{\mtx A}_n\begin{pmatrix}\mtx M_0\\\mtx M_1\end{pmatrix}.
\end{equation}

The following matrix will be the exact covariance of the centered,
normalized moment traces; applying the same inverse moment transformation
then gives the covariance of the variance-component estimators:
\begin{align}
\mtx U_{\mathrm{re,h}}
&=\frac2{q^2}\Bigg\{
\|\mtx\Sigma_b\|_F^2\begin{pmatrix}g_2&g_3\\g_3&g_4\end{pmatrix}
+2\tr(\mtx\Sigma_e\mtx\Sigma_b)
 \begin{pmatrix}1&g_2\\g_2&g_3\end{pmatrix}
+\|\mtx\Sigma_e\|_F^2\begin{pmatrix}1&1\\1&g_2\end{pmatrix}
\Bigg\},
\label{eq:U_re_homo_fixed_mt}\\
\mtx V_{\mathrm{re,h}}
&=\widehat{\mtx A}_n\mtx U_{\mathrm{re,h}}
  \widehat{\mtx A}_n^\top.\notag
\end{align}
To estimate these covariances, replace the two response covariance
matrices in \eqref{eq:U_re_homo_fixed_mt} by their estimates, obtaining
$\widehat{\mtx U}_{\mathrm{re,h}}$, and put
$\widehat{\mtx V}_{\mathrm{re,h}}=
\widehat{\mtx A}_n\widehat{\mtx U}_{\mathrm{re,h}}
\widehat{\mtx A}_n^\top$. The inverse moment matrix depends only on the
fixed design and is unchanged in this substitution. The main-paper SNR
variances have the equivalent representations
\begin{equation}
\label{eq:sigma_r_homo_fixed_mt}
\sigma_{r,\mathrm{re,h}}^2
=\vct d^\top\mtx V_{\mathrm{re,h}}\vct d,
\qquad
\hat\sigma_{r,\mathrm{re,h}}^2
=\hat{\vct d}^\top\widehat{\mtx V}_{\mathrm{re,h}}\hat{\vct d}.
\end{equation}
To verify the equivalence with the scalar formula, use
$\tau=\rho^2+\sigma^2$ and the transformed gradient
\begin{equation}
\label{eq:fixed_homo_transformed_gradient}
\widehat{\mtx A}_n^\top\vct d
=\frac1{(g_2-1)\tau^2}
\begin{pmatrix}-(g_2\rho^2+\sigma^2)\\\tau\end{pmatrix}.
\end{equation}
To check the scalar coefficients explicitly, put
$\delta=g_2-1$ and $\gamma=g_2\rho^2+\sigma^2
=g_2\tau-\delta\sigma^2$. Contracting the signal, mixed and noise
moment matrices in \eqref{eq:U_re_homo_fixed_mt} gives numerators
$\gamma^2g_2-2\gamma\tau g_3+\tau^2g_4$,
$\gamma^2-2\gamma\tau g_2+\tau^2g_3$, and
$\gamma^2-2\gamma\tau+\tau^2g_2$, respectively. Their reductions are
\begin{align*}
\frac{\gamma^2g_2-2\gamma\tau g_3+\tau^2g_4}{\delta^2}
&=g_2\sigma^4+
 \frac{2(g_3-g_2^2)}\delta\sigma^2\tau
 +\frac{g_4-2g_2g_3+g_2^3}{\delta^2}\tau^2,\\
\frac{\gamma^2-2\gamma\tau g_2+\tau^2g_3}{\delta^2}
&=\sigma^4+\frac{g_3-g_2^2}{\delta^2}\tau^2,\\
\frac{\gamma^2-2\gamma\tau+\tau^2g_2}{\delta^2}
&=\rho^4+\frac{\tau^2}{\delta}.
\end{align*}
Multiplication by the common factor $2/(q^2\tau^4)$ and the respective
trace functionals gives \eqref{eq:fixed_homo_variance_function}.
These are algebraic identities, so the same substitution with estimates
gives \eqref{eq:fixed_homo_scalar_variance_main}.

The next lemma verifies the exact covariance representation, bounds its
scale, and controls the normal-approximation and plug-in errors. These are
the model-specific inputs to the ratio reduction in S1.

\begin{lemma}[Fixed-design random-effects approximation]
\label{lem:random_fixed_homo_approximation}
Under the conditions of
Theorem~\ref{thm:asymptotics_random_fixed_homo_mt}, the covariance estimators
are conditionally unbiased. Define
\begin{equation}
\label{eq:random_centered_moments}
\vct w=\frac{\sqrt n}{q}
\begin{pmatrix}
\tr(\mtx M_0)-\E\{\tr(\mtx M_0)\mid\mtx X\}\\
\tr(\mtx M_1)-\E\{\tr(\mtx M_1)\mid\mtx X\}
\end{pmatrix}.
\end{equation}
Then
\begin{equation}
\label{eq:fixed_homo_lemma_exact_relations}
\vct\Delta=\frac1{\sqrt n}\widehat{\mtx A}_n\vct w,
\qquad
\Cov(\vct w\mid\mtx X)=\mtx U_{\mathrm{re,h}},
\end{equation}
and
\begin{equation}
\label{eq:fixed_homo_lemma_matrix_scale}
\frac cq\mtx I_2\preceq\mtx U_{\mathrm{re,h}}
\preceq\frac{C\ell_n}{q}\mtx I_2,
\qquad
\frac cq\mtx I_2\preceq\mtx V_{\mathrm{re,h}}
\preceq\frac{C\ell_n}{q}\mtx I_2.
\end{equation}
Uniformly over nonzero $\vct b\in\mathbb R^2$,
\begin{equation}
\label{eq:fixed_homo_lemma_quadratic_clt}
\sup_t\left|
\PP\!\left(
\frac{\vct b^\top\vct w}
{(\vct b^\top\mtx U_{\mathrm{re,h}}\vct b)^{1/2}}
\le t\,\middle|\,\mtx X\right)-\Phi(t)\right|
\le\frac{C\ell_n}{\sqrt{nq}}=o(1).
\end{equation}
Finally,
\begin{align}
q^{-2}\|\widehat{\mtx\Sigma}_b-\mtx\Sigma_b\|_F^2
+q^{-2}\|\widehat{\mtx\Sigma}_e-\mtx\Sigma_e\|_F^2
&=O_P(n^{-1}),
\label{eq:fixed_homo_lemma_cov_est_rate}\\
\|\widehat{\mtx V}_{\mathrm{re,h}}-\mtx V_{\mathrm{re,h}}\|
&=o_P(q^{-1}).
\label{eq:fixed_homo_lemma_V_est_rate}
\end{align}
All constants and rates are uniform over the stated parameter sequences,
with either fixed or diverging $q$.
\end{lemma}

\begin{proof}
\proofblock{Exact conditional moments.}
For a deterministic symmetric $n\times n$ matrix $\mtx C$,
\eqref{eq:fixed_homo_direct_covariance} gives
\begin{equation}
\label{eq:fixed_homo_general_mean}
\E(\mtx Y^\top\mtx C\mtx Y\mid\mtx X)
=\tr(\mtx C\mtx H)\mtx\Sigma_b+\tr(\mtx C)\mtx\Sigma_e.
\end{equation}
Applying this with $\mtx C=\mtx I_n,\mtx H$, and using
$\tr(\mtx H^k)=ng_k$, proves
\eqref{eq:random_conditional_moment_system}. Inverting this system in
\eqref{eq:fixed_homo_matrix_moment_system} gives conditional unbiasedness;
taking normalized traces and centering gives the first identity in
\eqref{eq:fixed_homo_lemma_exact_relations}.

For the covariance calculation,
$\tr(\mtx Y^\top\mtx C\mtx Y)=
\vecop(\mtx Y)^\top(\mtx I_q\otimes\mtx C)\vecop(\mtx Y)$.
The centered Gaussian quadratic-form identity therefore gives
\begin{align}
&\Cov\{\tr(\mtx Y^\top\mtx C\mtx Y),
             \tr(\mtx Y^\top\mtx D\mtx Y)\mid\mtx X\}\notag\\
&\quad=2\tr\{(\mtx I_q\otimes\mtx C)\mtx\Omega_X
             (\mtx I_q\otimes\mtx D)\mtx\Omega_X\}.
\label{eq:fixed_homo_gaussian_covariance}
\end{align}
Expanding the Kronecker trace gives
\begin{align*}
&2\|\mtx\Sigma_b\|_F^2\tr(\mtx C\mtx H\mtx D\mtx H)
+2\tr(\mtx\Sigma_b\mtx\Sigma_e)\tr(\mtx C\mtx H\mtx D)\\
&\qquad+2\tr(\mtx\Sigma_e\mtx\Sigma_b)\tr(\mtx C\mtx D\mtx H)
+2\|\mtx\Sigma_e\|_F^2\tr(\mtx C\mtx D).
\end{align*}
For $\mtx C=\mtx H^j$ and $\mtx D=\mtx H^k$, the two mixed traces
coincide. Dividing by $nq^2$, as required by the definition of $\vct w$,
yields, for $j,k\in\{0,1\}$,
\begin{equation}
\label{eq:fixed_homo_indexed_covariance}
\Cov(w_j,w_k\mid\mtx X)
=\frac2{q^2}\left\{
\|\mtx\Sigma_b\|_F^2g_{j+k+2}
+2\tr(\mtx\Sigma_e\mtx\Sigma_b)g_{j+k+1}
+\|\mtx\Sigma_e\|_F^2g_{j+k}\right\},
\end{equation}
where $w_0,w_1$ denote the two entries of $\vct w$ in that order.
This is exactly \eqref{eq:U_re_homo_fixed_mt}.

\proofblock{Covariance scale and normal approximation.}
The noise contribution alone supplies a uniform lower covariance bound;
the signal and mixed contributions are positive semidefinite. To see this,
let $\lambda_1,\ldots,\lambda_n$ be the eigenvalues of $\mtx H$.
For $\ell=1,2$, the moment matrix
$\begin{pmatrix}g_\ell&g_{\ell+1}\\g_{\ell+1}&g_{\ell+2}\end{pmatrix}$
has quadratic form
\[
\frac1n\sum_{i=1}^n\lambda_i^\ell(u+v\lambda_i)^2\ge0.
\]
Thus the signal and mixed moment matrices are positive semidefinite. The
noise matrix $\begin{pmatrix}1&1\\1&g_2\end{pmatrix}$ has determinant
$g_2-1\ge c$ and bounded trace, so its eigenvalues are bounded above and
away from zero. For either response covariance $\mtx\Sigma_a$,
\[
 cq\le\frac{\tr^2(\mtx\Sigma_a)}q
 \le\|\mtx\Sigma_a\|_F^2
 \le\|\mtx\Sigma_a\|\tr(\mtx\Sigma_a)
 \le C\ell_nq,
\]
and $0\le\tr(\mtx\Sigma_e\mtx\Sigma_b)\le C\ell_nq$.
The noise term supplies the lower bound, and boundedness of $g_k$ supplies
the upper bound:
\begin{equation}
\label{eq:fixed_homo_U_scale}
\frac cq\mtx I_2\preceq\mtx U_{\mathrm{re,h}}
\preceq\frac{C\ell_n}{q}\mtx I_2.
\end{equation}
Since $\widehat{\mtx A}_n$ and its inverse have bounded norm,
conjugation gives
\begin{equation}
\label{eq:fixed_homo_V_scale}
\frac cq\mtx I_2\preceq\mtx V_{\mathrm{re,h}}
\preceq\frac{C\ell_n}{q}\mtx I_2.
\end{equation}

The lower bound just proved ensures that no direction has too small a
variance. We now compare it with an operator-norm bound for the associated
quadratic form. Write
$\vecop(\mtx Y)=\mtx\Omega_X^{1/2}\vct z$ with
$\vct z\sim\mathcal N(\vct0,\mtx I_{nq})$. Since
$\|\mtx\Omega_X\|\le\|\mtx\Sigma_b\|\|\mtx H\|
+\|\mtx\Sigma_e\|\le C\ell_n$, for
$\vct b=(b_0,b_1)^\top\ne\vct0$ the matrix
\[
\mtx Q_{\vct b}
=\frac1{\sqrt n\,q}\mtx\Omega_X^{1/2}
\{\mtx I_q\otimes(b_0\mtx I_n+b_1\mtx H)\}
\mtx\Omega_X^{1/2}
\]
satisfies
\[
\vct b^\top\vct w=\vct z^\top\mtx Q_{\vct b}\vct z
-\tr(\mtx Q_{\vct b}),\qquad
2\|\mtx Q_{\vct b}\|_F^2
=\vct b^\top\mtx U_{\mathrm{re,h}}\vct b.
\]
Moreover,
\[
\|\mtx Q_{\vct b}\|\le
\frac{C\ell_n\|\vct b\|_2}{\sqrt n\,q},\qquad
\|\mtx Q_{\vct b}\|_F\ge\frac{c\|\vct b\|_2}{\sqrt q}.
\]
Applying \eqref{eq:gaussian_quadratic_be} proves
\eqref{eq:fixed_homo_lemma_quadratic_clt}. Its error tends to zero because
$\ell_n\le Cq$ and $\ell_nq/n\to0$.

\proofblock{Estimating the covariance functionals.}
As in S2, a Frobenius error bound suffices because the leading variance uses
only three trace products. Here the entrywise variances follow directly
from the Gaussian covariance blocks of the response columns:
\[
\mtx\Omega_{kl}=(\mtx\Sigma_b)_{kl}\mtx H
               +(\mtx\Sigma_e)_{kl}\mtx I_n.
\]
The Gaussian bilinear-form covariance identity gives, for $j=0,1$,
\begin{equation}
\label{eq:fixed_homo_entry_variance}
\E\!\left[
\{(\mtx M_j)_{kl}-\E((\mtx M_j)_{kl}\mid\mtx X)\}^2
\,\middle|\,\mtx X\right]
=\frac1{n^2}\left\{
\tr(\mtx C_j\mtx\Omega_{kk}\mtx C_j\mtx\Omega_{ll})
+\tr(\mtx C_j\mtx\Omega_{kl}\mtx C_j\mtx\Omega_{lk})\right\}.
\end{equation}
All these row-space matrices commute because they are polynomials in
$\mtx H$. In an eigenbasis of $\mtx H$, write
$\mtx R_i=\lambda_i\mtx\Sigma_b+\mtx\Sigma_e\succeq0$.
Summing \eqref{eq:fixed_homo_entry_variance} over $k,l$ gives the exact
identity
\[
\E\{\|\mtx M_j-\E(\mtx M_j\mid\mtx X)\|_F^2\mid\mtx X\}
=\frac1{n^2}\sum_{i=1}^n\lambda_i^{2j}
 \{\tr^2(\mtx R_i)+\tr(\mtx R_i^2)\},
\]
where the weight is one for $j=0$. This also shows why off-diagonal
covariance blocks cause no sign difficulty. Since $\lambda_i^{2j}\le C$,
the two sums in the entrywise formula satisfy

\begin{align*}
\sum_{k,l}\tr(\mtx C_j\mtx\Omega_{kk}\mtx C_j\mtx\Omega_{ll})
&\le C\Big\{
\tr^2(\mtx\Sigma_b)\tr(\mtx H^2)
+2\tr(\mtx\Sigma_b)\tr(\mtx\Sigma_e)\tr(\mtx H)
+n\tr^2(\mtx\Sigma_e)\Big\}
=O(nq^2),\\
\sum_{k,l}\tr(\mtx C_j\mtx\Omega_{kl}\mtx C_j\mtx\Omega_{lk})
&\le C\Big\{
\|\mtx\Sigma_b\|_F^2\tr(\mtx H^2)
+2\tr(\mtx\Sigma_e\mtx\Sigma_b)\tr(\mtx H)
+n\|\mtx\Sigma_e\|_F^2\Big\}
=O(n\ell_nq).
\end{align*}
Here $\tr(\mtx H)=n$, $\tr(\mtx H^2)=ng_2=O(n)$, and the normalized
response traces are bounded. Since $\ell_n\le Cq$, summing
\eqref{eq:fixed_homo_entry_variance} yields
\[
\E\{\|\mtx M_j-\E(\mtx M_j\mid\mtx X)\|_F^2\mid\mtx X\}
\le Cq^2/n,\qquad j=0,1.
\]
The bounded inverse in \eqref{eq:fixed_homo_matrix_moment_system} and
Markov's inequality prove \eqref{eq:fixed_homo_lemma_cov_est_rate}.

We can now reuse \eqref{eq:fe_trace_product_bound}: its inputs are precisely
the preceding Frobenius rate and the population bound
$\|\mtx\Sigma_a\|_F=O(\sqrt{\ell_nq})$. With the present covariance
estimates, it gives
\[
\max_{a,c\in\{b,e\}}\frac1q
\left|\tr(\widehat{\mtx\Sigma}_a\widehat{\mtx\Sigma}_c)
-\tr(\mtx\Sigma_a\mtx\Sigma_c)\right|
=O_P\!\left(\sqrt{\frac{\ell_nq}{n}}+\frac qn\right)=o_P(1).
\]
The coefficients in \eqref{eq:U_re_homo_fixed_mt} are bounded, so
$\|\widehat{\mtx U}_{\mathrm{re,h}}-\mtx U_{\mathrm{re,h}}\|=o_P(q^{-1})$.
Finally, the exact identity
\[
\widehat{\mtx V}_{\mathrm{re,h}}-\mtx V_{\mathrm{re,h}}
=\widehat{\mtx A}_n
(\widehat{\mtx U}_{\mathrm{re,h}}-\mtx U_{\mathrm{re,h}})
\widehat{\mtx A}_n^\top
\]
proves \eqref{eq:fixed_homo_lemma_V_est_rate}.
\end{proof}

\begin{proof}[Proof of Theorem~\ref{thm:asymptotics_random_fixed_homo_mt}]
Only the SNR direction and consistency of the total variance remain to be
checked. Apply \eqref{eq:fixed_homo_lemma_quadratic_clt} with
$\vct b=\widehat{\mtx A}_n^\top\vct d$. The exact estimator relation
\eqref{eq:fixed_homo_lemma_exact_relations} gives
\begin{equation}
\label{eq:fixed_homo_linear_snr_clt}
\frac{\sqrt n\,\vct d^\top\vct\Delta}
{(\vct d^\top\mtx V_{\mathrm{re,h}}\vct d)^{1/2}}
\Longrightarrow\mathcal N(0,1)
\end{equation}
conditionally on $\mtx X$. Also
$\hat\tau=\tr(\mtx M_0)/q$, and the exact conditional moments give
\[
\E(\hat\tau\mid\mtx X)=\tau,\qquad
\Var(\hat\tau\mid\mtx X)
=\frac1n(\mtx U_{\mathrm{re,h}})_{11}
\le\frac{C\ell_n}{nq}.
\]
Thus $\hat\tau/\tau\to_P1$. The scale bound
\eqref{eq:fixed_homo_lemma_matrix_scale}, the plug-in bound
\eqref{eq:fixed_homo_lemma_V_est_rate}, and $\ell_n/(nq)\to0$ verify
Lemma~\ref{lem:common_studentization} conditionally on the fixed design.
Together with \eqref{eq:sigma_r_homo_fixed_mt}, it proves the studentized
limit and conditional coverage for both fixed and diverging $q$.
\end{proof}

\proofblock{Fixed projection for intercept removal.}
We finish by checking that the intercept removal in
Section~\ref{sec:yeast_application} preserves the homogeneous Gaussian
model in a lower-dimensional row space. Let $\mtx C$ be a fixed
full-column-rank matrix of $d$ nuisance covariates, and choose a fixed
$n\times m$ matrix $\mtx Q$, $m=n-d$, with
$\mtx Q^\top\mtx Q=\mtx I_m$ and $\mtx Q^\top\mtx C=0$.
Under the homogeneous fixed-design model
\[
\mtx Y=\mtx C\mtx A_0+\mtx X\mtx B+\mtx E,
\]
projection gives
$\mtx Q^\top\mtx Y=(\mtx Q^\top\mtx X)\mtx B+\mtx Q^\top\mtx E$.
Furthermore,
\[
\Cov\{\vecop(\mtx Q^\top\mtx E)\mid\mtx X\}
=\mtx\Sigma_e\otimes(\mtx Q^\top\mtx Q)
=\mtx\Sigma_e\otimes\mtx I_m.
\]
Joint Gaussianity gives independent rows of $\mtx Q^\top\mtx E$, still
independent of $\mtx B$. The signal covariance for this projected model is
\[
\mtx\Sigma_{b,Q}
=\frac{\|\mtx Q^\top\mtx X\|_F^2}{mp}\mtx\Sigma_\beta.
\]
Thus the fixed-design theorem applies with $m$ and $\mtx\Sigma_{b,Q}$
in place of $n$ and $\mtx\Sigma_b$, provided all of its conditions hold
for the projected model: $p\asymp m$, the design bounds, the normalized
trace bounds, and the response-spectrum growth condition. The projected
signal covariance need not equal the unprojected one. For intercept
removal, $\mtx C=\boldsymbol1_n$. Independence concerns the $m$ projected
noise coordinates, not the $n$ rows of $\mtx Q\mtx Q^\top\mtx Y$.

In the yeast application, scaling each projected marker to squared length
$m$ gives average predictor energy
$\|\mtx Q^\top\mtx X\|_F^2/(mp)=1$, although the theorem does not require
this normalization. The i.i.d.\ random-coefficient assumption is imposed on
that final marker scale: unequal predictor rescalings change the
coefficient model and do not preserve an arbitrary earlier one. The
responses remain on the chosen fixed scale. Neither a projection selected
from the responses nor heterogeneous noise generally preserves the
independent-row covariance used above.

\section{Proofs of Theorems~\ref{thm:asymptotics_mt} and~\ref{thm:studentized_random_mt}}
\label{sec:supp_random_design_theorem}
The random-design proof uses the same conditional Gaussian argument as S3,
but two additional terms must be controlled. The realized predictor energy
differs from its population value, creating a conditional bias, and
observation-specific noise produces weighted design traces in the
conditional covariance. We first bound these design terms and prove the
oracle limit. The final subsection estimates the extra covariance
dispersion under scalar heterogeneity and completes studentization.

\subsection{Design bounds and weighted trace reductions}
\label{sec:proof_of_asymptotics_mt}

Throughout S4, the model is $\mtx Y=\mtx X\mtx B+\mtx E$. To distinguish
the population predictor scale from its observed counterpart, put
\[
\mu_x=\frac{\tr(\mtx\Sigma)}p,\qquad
\hat\mu_x=\frac{\|\mtx X\|_F^2}{np},\qquad
a_n=\frac{\hat\mu_x}{\mu_x},\qquad
\ell_n=1\vee\|\mtx\Sigma_b\|\vee\max_i\|\mtx\Sigma_i\|.
\]
Thus $\mtx\Sigma_b=\mu_x\mtx\Sigma_\beta$ and
$\|\mtx\Sigma\|\le C\mu_x$. The ratio $a_n$ measures the difference
between observed and population energy; the scale $\mu_x$ itself need not
be bounded. The estimator uses the normalized design weights
\begin{equation}
\label{eq:re_direct_design_weights}
\mtx K=\mtx X\mtx X^\top,\qquad
\mtx H=\frac{n\mtx K}{\tr(\mtx K)}
=\frac{\mtx K}{p\hat\mu_x},\qquad
t_i=H_{ii},\qquad
g_k=\frac1n\tr(\mtx H^k).
\end{equation}
Here $g_0=g_1=1$, and $g_2,g_3,g_4$ are the main-paper summaries in
\eqref{eq:random_design_gk_main}. The diagonal weight $t_i$ records the
relative predictor energy of observation $i$. All design statements are
made on $\hat\mu_x>0$, an event whose probability tends to one.

Define the two matrix moments and their inverse transformation by
\begin{equation}
\label{eq:re_direct_moment_system}
\mtx M_0=\frac1n\mtx Y^\top\mtx Y,\qquad
\mtx M_1=\frac1n\mtx Y^\top\mtx H\mtx Y,\qquad
\widehat{\mtx A}_n
=\frac1{g_2-1}
\begin{pmatrix}g_2&-1\\-1&1\end{pmatrix}.
\end{equation}
The main-paper estimators satisfy the exact identity
\begin{equation}
\label{eq:re_direct_estimators}
\begin{pmatrix}\widehat{\bar{\mtx\Sigma}}_e\\
\widehat{\mtx\Sigma}_b\end{pmatrix}
=\widehat{\mtx A}_n
\begin{pmatrix}\mtx M_0\\\mtx M_1\end{pmatrix}.
\end{equation}
This is an algebraic identity for the observed estimator, not an assertion
of conditional unbiasedness. The conditional means will involve both
$a_n-1$ and the alignment between $t_i$ and the noise covariance.

\proofblock{Design bounds.}
The first lemma controls the global and rowwise energies, bounds the
spectral moments, and ensures that the inverse moment matrix remains
bounded.

\begin{lemma}[Design energy and spectral gap]
\label{lem:gkhat_bounds}
Under the design conditions of Theorem~\ref{thm:asymptotics_mt},
\begin{equation}
\label{eq:empirical_global_scale_rate}
a_n=1+O_P(n^{-1}),\qquad
\max_{i\le n}|t_i-1|
=O_P\!\left(\sqrt{\frac{\log n}{n}}\right).
\end{equation}
Moreover, there are fixed $c,C>0$ such that, on an event of probability
tending to one, $1/2\le a_n\le2$ and
\begin{equation}
\label{eq:common_design_event}
\|\mtx H\|\le C,\qquad
\max_{k\le4}g_k\le C,\qquad
g_2-1\ge c.
\end{equation}
\end{lemma}

\begin{proof}
Independence of the rows and the quadratic-form variance bound give
\[
\E\hat\mu_x=\mu_x,\qquad
\Var(\hat\mu_x)
=\frac{\Var(\|\vct x_1\|_2^2)}{np^2}
\le\frac{C\tr(\mtx\Sigma^2)}{np^2}
\le\frac{C\mu_x^2}{np}.
\]
Since $p\asymp n$, this gives the $O_P(n^{-1})$ relative error for the
global energy in \eqref{eq:empirical_global_scale_rate}. For the rowwise
energies, Hanson--Wright concentration and a union bound yield
\[
\max_i\big|\|\vct x_i\|_2^2-p\mu_x\big|
=O_P\{\mu_x\sqrt{p\log n}\}.
\]
Dividing by $p\hat\mu_x$ proves the row-weight bound. Also,
\eqref{eq:subg_matrix_norm_tail} gives
$\|\mtx X\|=O_P(\sqrt{n\mu_x})$ with a fixed-constant bound holding
with probability tending to one. Thus $\|\mtx H\|\le C$ on such an event,
and $g_k\le\|\mtx H\|^k$.

To obtain the lower bound on $g_2-1$, it is enough to retain the
nonnegative off-diagonal contribution. Indeed,
\[
g_2-1=\frac1n\|\mtx H-\mtx I_n\|_F^2
=\frac1n\sum_i(t_i-1)^2
+\frac1{n(p\hat\mu_x)^2}\sum_{i\ne j}(\vct x_i^\top\vct x_j)^2.
\]
It is therefore enough to control the off-diagonal energy
\[
J_n=\frac1{n(p\mu_x)^2}\sum_{i\ne j}(\vct x_i^\top\vct x_j)^2,
\qquad g_2-1\ge a_n^{-2}J_n.
\]
Independence of the rows gives
\begin{equation}
\label{eq:design_offdiagonal_expectation}
\E J_n=\frac{(n-1)\tr(\mtx\Sigma^2)}{p^2\mu_x^2}
\ge\frac{n-1}{p},
\end{equation}
since $\tr(\mtx\Sigma^2)\ge p\mu_x^2$.
For the variance, condition on one row and use the sub-Gaussian fourth-moment
and quadratic-form bounds:
\begin{align}
\E(\vct x_1^\top\vct x_2)^4
&\le C\E(\vct x_1^\top\mtx\Sigma\vct x_1)^2
\le C\|\mtx\Sigma\|^2\E\|\vct x_1\|_2^4
\le Cp^2\mu_x^4,\notag\\
\Cov\{(\vct x_1^\top\vct x_2)^2,(\vct x_1^\top\vct x_3)^2\}
&=\Var(\vct x_1^\top\mtx\Sigma\vct x_1)
\le C\tr(\mtx\Sigma^4)\le Cp\mu_x^4.
\label{eq:design_pair_covariances}
\end{align}
Disjoint pairs are independent. There are $O(n^2)$ equal or reversed pairs
and $O(n^3)$ pairs sharing one index, so
\[
\Var(J_n)
\le\frac{C(n^2p^2+n^3p)\mu_x^4}{n^2p^4\mu_x^4}
=C\left(\frac1{p^2}+\frac n{p^3}\right)=O(n^{-2}).
\]
Chebyshev's inequality, $p\asymp n$ and
\eqref{eq:design_offdiagonal_expectation} imply that $J_n$ is bounded away
from zero with probability tending to one. Since $a_n\to_P1$, the same is
true of $g_2-1$, completing the proof.
\end{proof}

\proofblock{Weighted traces.}
The conditional covariance contains design traces weighted by deterministic
noise quantities. The next lemma replaces a linear weight by its average.
Its sharper rate for $k=1$ will also control the conditional bias; the rates
for $k=2,3$ are needed only for covariance approximation. The weights in
this lemma are arbitrary bounded scalars and need not be noise scales.

\begin{lemma}[Weighted design traces]
\label{lem:leave one out analysis}
Let $\nu_i$ be deterministic scalar weights with
$\max_i|\nu_i|\le C$ and $\bar\nu=n^{-1}\sum_i\nu_i$.
Under the design conditions of Theorem~\ref{thm:asymptotics_mt},
\begin{align}
\frac1n\tr\{\mtx H\diag(\nu_i)\}
&=\bar\nu+O_P(n^{-1}),
\label{equation: leave one out eq 2}\\
\frac1n\tr\{\mtx H^k\diag(\nu_i)\}
&=\bar\nu g_k+O_P(n^{-1/2}),\qquad k=2,3.
\label{equation: leave one out eq 3}
\end{align}
\end{lemma}

\begin{proof}
Center the weights by writing $c_i=\nu_i-\bar\nu$, so their sum is zero.
For $k=1$, the centered numerator is
$\sum_i c_i\|\vct x_i\|_2^2$, with variance at most $Cnp\mu_x^2$.
After division by $np\hat\mu_x$, its order is $O_P(n^{-1})$.

For $k=2,3$, independence of the diagonal entries is unavailable, but
exchangeability and the zero-sum weights suffice. Set
$Z_i=(\mtx K^k)_{ii}
=\vct x_i^\top(\mtx X^\top\mtx X)^{k-1}\vct x_i$.
Exchangeability and $\sum_i c_i=0$ give
\[
\E\left(\sum_i c_iZ_i\right)^2
=\sum_i c_i^2\{\E Z_1^2-\E(Z_1Z_2)\}
\le\sum_i c_i^2\,\E Z_1^2,
\]
where $Z_i\ge0$ was used in the last step.
The matrix-norm tail bound and the row-moment bound in S1 imply
\[
\E\|\mtx X\|^{8k-8}\le C(n\mu_x)^{4k-4},\qquad
\E\|\vct x_1\|_2^8\le C(p\mu_x)^4.
\]
Hence, by Cauchy--Schwarz,
\[
\E Z_1^2
\le\{\E\|\mtx X\|^{8k-8}\}^{1/2}
   \{\E\|\vct x_1\|_2^8\}^{1/2}
\le C(n\mu_x)^{2k-2}(p\mu_x)^2
\le C(p\mu_x)^{2k}.
\]
It follows that
$\{n(p\mu_x)^k\}^{-1}\sum_i c_iZ_i=O_P(n^{-1/2})$.
Multiplication by $a_n^{-k}=O_P(1)$ gives the asserted bound for
$n^{-1}\tr\{\mtx H^k\diag(c_i)\}$.
\end{proof}

A noise--noise term contains two centered weights rather than one. The
following bound shows that its leading contribution comes from the
diagonal terms. Applied to vectorized covariance deviations, this will
produce $\kappa_{\mathrm{tot}}$.

\begin{lemma}[Centered quadratic trace]
\label{lem:centered_trace_reduction}
Under the design conditions of Theorem~\ref{thm:asymptotics_mt}, let
$b_n\ge1$ be deterministic and let deterministic vectors
$\vct u_i$ satisfy
$\sum_i\vct u_i=0$ and $\max_i\|\vct u_i\|_2^2\le Cb_n$.
Then
\begin{equation}
\label{eq:centered_trace_reduction}
\frac1n\sum_{i,j}H_{ij}^2\vct u_i^\top\vct u_j
=\frac1n\sum_i\|\vct u_i\|_2^2
+O_P\!\left(b_n\sqrt{\frac{\log n}{n}}\right).
\end{equation}
\end{lemma}

\begin{proof}
We separate the diagonal contribution from the off-diagonal remainder.
For the diagonal terms, \eqref{eq:empirical_global_scale_rate} gives
\[
\left|\frac1n\sum_i(t_i^2-1)\|\vct u_i\|_2^2\right|
=O_P\!\left(b_n\sqrt{\frac{\log n}{n}}\right).
\]
For the off-diagonal terms, put $c_{ij}=\vct u_i^\top\vct u_j$.
Then $|c_{ij}|\le Cb_n$ and
$\sum_{i\ne j}c_{ij}=-\sum_i\|\vct u_i\|_2^2$.
In the raw numerator
$\sum_{i\ne j}c_{ij}(\vct x_i^\top\vct x_j)^2$, insert and subtract
$\tr(\mtx\Sigma^2)$. The deterministic part, divided by
$n(p\hat\mu_x)^2$, is $O_P(b_n/p)$.

For the centered part, we use the same pair classification as in
\eqref{eq:design_pair_covariances}. Equal or reversed pairs have covariance
bounded by $Cp^2\mu_x^4$. Pairs sharing exactly one index have covariance
\[
\Cov\{(\vct x_i^\top\vct x_j)^2,
       (\vct x_i^\top\vct x_k)^2\}
=\Var(\vct x_i^\top\mtx\Sigma\vct x_i)
\le Cp\mu_x^4,
\]
and disjoint pairs are independent.
There are $O(n^2)$ pairs of the first type and $O(n^3)$ of the second.
Thus the centered numerator, divided by $n(p\mu_x)^2$, has variance
at most
\[
\frac{Cb_n^2(n^2p^2+n^3p)\mu_x^4}
     {n^2p^4\mu_x^4}
=O(b_n^2/n^2).
\]
Replacing $\mu_x$ by $\hat\mu_x$ multiplies the centered expression by
$a_n^{-2}=O_P(1)$. The off-diagonal terms are therefore $O_P(b_n/n)$,
which completes the proof.
\end{proof}

\proofblock{Aggregating the noise contributions.}
We now apply the two weighted-trace bounds to the response covariance
matrices. The goal is to express all terms using three homogeneous
covariance functionals and a single dispersion term. Write $\mtx E_0=\bar{\mtx\Sigma}_e$ and
$\mtx D_i=\mtx\Sigma_i-\mtx E_0$, and set
\[
F_b=\|\mtx\Sigma_b\|_F^2,\qquad
F_e=\|\mtx E_0\|_F^2,\qquad
F_{be}=\tr(\mtx E_0\mtx\Sigma_b),\qquad
\kappa=\kappa_{\mathrm{tot}}.
\]
These four functionals are $O(\ell_nq)$; for example,
$n^{-1}\sum_i\tr(\mtx\Sigma_i^2)
\le\ell_n\tr(\mtx E_0)=O(\ell_nq)$. The errors below are controlled at
the smaller scale $o_P(q)$: the covariance calculation multiplies each
trace by $q^{-2}$, so this yields the required $o_P(q^{-1})$ error.

For the signal--noise terms, let $\mtx\Lambda_{kl}=\diag\{(\mtx\Sigma_i)_{kl}:1\le i\le n\}$ and
$\mtx T=\diag\{\tr(\mtx\Sigma_i\mtx\Sigma_b):1\le i\le n\}$.
The normalized weights needed below are uniformly bounded, even though
individual noise traces need not be $O(q)$:
\begin{align*}
0\le\tr(\mtx\Sigma_i\mtx\Sigma_b)
&\le\|\mtx\Sigma_i\|\tr(\mtx\Sigma_b)\le Cq\ell_n,\\
\tr(\mtx\Sigma_i^2)&\le q\ell_n^2,\qquad
\|\mtx D_i\|_F^2\le4q\ell_n^2.
\end{align*}
Here $\|\mtx E_0\|\le\ell_n$, so $\|\mtx D_i\|\le2\ell_n$.
The weighted-trace lemma, applied to
$q^{-1}\tr(\mtx\Sigma_i\mtx\Sigma_b)/\ell_n$, therefore gives
\begin{equation}
\label{eq:weighted_gk_reduction}
\frac1n\tr(\mtx H^k\mtx T)=g_kF_{be}+o_P(q),
\qquad k=1,2,3.
\end{equation}
The errors are $O_P(q\ell_n/n)$ for $k=1$ and
$O_P(q\ell_n/\sqrt n)$ for $k=2,3$; the oracle growth conditions make
both $o_P(q)$.

The three noise--noise terms are
\begin{align}
\frac1n\sum_{k,l}\tr(\mtx\Lambda_{kl}^2)
&=F_e+\kappa,\notag\\
\frac1n\sum_{k,l}\tr(\mtx H\mtx\Lambda_{kl}^2)
&=F_e+\kappa+o_P(q),
\label{eq:aggregate_trace_reduction_2}\\
\frac1n\sum_{k,l}\tr(\mtx H\mtx\Lambda_{kl}
                         \mtx H\mtx\Lambda_{kl})
&=g_2F_e+\kappa+o_P(q).
\label{eq:aggregate_trace_reduction_1}
\end{align}
The first identity follows from
$n^{-1}\sum_i\tr(\mtx\Sigma_i^2)=F_e+\kappa$.
For the second, apply \eqref{equation: leave one out eq 2} to
$\tr(\mtx\Sigma_i^2)/(q\ell_n^2)$; the error is
$O_P(q\ell_n^2/n)=o_P(q)$.

For \eqref{eq:aggregate_trace_reduction_1}, the off-diagonal interactions
between observations must also be retained. Its left-hand side is
\[
\frac1n\sum_{i,j}H_{ij}^2\tr(\mtx\Sigma_i\mtx\Sigma_j).
\]
Substituting $\mtx\Sigma_i=\mtx E_0+\mtx D_i$ splits this expression into
\[
g_2F_e
+\frac2n\tr\!\left[
\mtx H^2\diag\{\tr(\mtx E_0\mtx D_i)\}\right]
+\frac1n\sum_{i,j}H_{ij}^2\tr(\mtx D_i\mtx D_j).
\]
The middle term has centered weights because $\sum_i\mtx D_i=0$, and is
$O_P(q\ell_n^2/\sqrt n)$ by
\eqref{equation: leave one out eq 3}. For the final term, apply
\eqref{eq:centered_trace_reduction} to
$\vct u_i=\vecop(\mtx D_i)$ and $b_n=q\ell_n^2$ to obtain
$\kappa+O_P\{q\ell_n^2\sqrt{\log n/n}\}$ for the last term.
Both errors are $o_P(q)$ because $p\asymp n$ and
$\ell_n^2\sqrt{\log p/n}\to0$. Thus covariance heterogeneity remains in
the leading term through $\kappa$, while the other effects of the
observation-specific weights are negligible at the required scale.

\subsection{Oracle normal approximation}
\label{sec:proof_random_design_oracle}

\begin{proof}[Proof of Theorem~\ref{thm:asymptotics_mt}]
We first center the estimator at its conditional mean and use its exact
conditional covariance. The bounds from the preceding subsection then
allow us to replace that covariance by the main-paper oracle variance and
the conditional mean by the population target. Work on the design event
\eqref{eq:common_design_event}, where
$\|\widehat{\mtx A}_n\|+\|\widehat{\mtx A}_n^{-1}\|\le C$; its
probability tends to one.

\proofblock{Conditional bias.}
The conditional covariance of $\vecop(\mtx Y)$ is
\begin{equation}
\label{eq:re_heterogeneous_omega}
\mtx\Omega_X
=\frac1p\mtx\Sigma_\beta\otimes\mtx K
+[\mtx\Lambda_{kl}]_{k,l=1}^q
=a_n\mtx\Sigma_b\otimes\mtx H
+[\mtx\Lambda_{kl}]_{k,l=1}^q.
\end{equation}
The factor $a_n$ appears because $\mtx H$ uses the observed predictor
energy, whereas $\mtx\Sigma_b$ uses its population value. No change to the
coefficient distribution is made. The two conditional moment equations are
\[
\E(\mtx M_0\mid\mtx X)=\bar{\mtx\Sigma}_e+a_n\mtx\Sigma_b,\qquad
\E(\mtx M_1\mid\mtx X)
=\bar{\mtx\Sigma}_e+a_ng_2\mtx\Sigma_b+\mtx R_X,
\]
where
\[
\mtx R_X=\frac1n\sum_i(t_i-1)\mtx\Sigma_i,\qquad
r_X=\frac1q\tr(\mtx R_X).
\]
Because $q^{-1}\tr(\mtx\Sigma_i)\le\ell_n$, the weighted-trace bound
\eqref{equation: leave one out eq 2} gives $r_X=O_P(\ell_n/n)$.
Using \eqref{eq:re_direct_estimators}, we obtain the exact conditional bias
\begin{equation}
\label{eq:re_cond_bias_vector}
\E(\widehat{\vct\theta}\mid\mtx X)-\vct\theta
=\begin{pmatrix}0\\(a_n-1)\rho^2\end{pmatrix}
+\widehat{\mtx A}_n\begin{pmatrix}0\\r_X\end{pmatrix},
\qquad
\widehat{\vct\theta}
=\begin{pmatrix}\hat\sigma^2\\\hat\rho^2\end{pmatrix}.
\end{equation}
The first term is the global-energy error; the second is the effect of
unequal row weights on heterogeneous noise. Since $a_n-1=O_P(n^{-1})$
and $\rho^2$ is bounded,
\begin{equation}
\label{eq:conditional_mean_consistency}
\|\E(\widehat{\vct\theta}\mid\mtx X)-\vct\theta\|_2
=O_P(\ell_n/n).
\end{equation}

\proofblock{Leading conditional covariance.}
We distinguish the exact conditional covariance $\mtx U_n$ from its
leading approximation $\widetilde{\mtx U}_n$. Define the centered moment
vector
\begin{equation}
\label{eq:re_direct_centered_moments}
\vct w=\frac{\sqrt n}{q}
\begin{pmatrix}
\tr(\mtx M_0)-\E\{\tr(\mtx M_0)\mid\mtx X\}\\
\tr(\mtx M_1)-\E\{\tr(\mtx M_1)\mid\mtx X\}
\end{pmatrix},
\qquad
\mtx U_n=\Cov(\vct w\mid\mtx X).
\end{equation}
For $j,k\in\{0,1\}$, the Gaussian quadratic-form covariance identity gives
\begin{align}
(\mtx U_n)_{j+1,k+1}
&=\frac{2}{nq^2}
\tr\{(\mtx I_q\otimes\mtx H^j)\mtx\Omega_X
      (\mtx I_q\otimes\mtx H^k)\mtx\Omega_X\}\notag\\
&=\frac{2a_n^2F_b}{q^2}g_{j+k+2}
+\frac{4a_n}{nq^2}\tr(\mtx H^{j+k+1}\mtx T)
+\frac{2}{nq^2}\sum_{r,s}
\tr(\mtx H^j\mtx\Lambda_{rs}\mtx H^k\mtx\Lambda_{rs}).
\label{eq:re_general_conditional_covariance}
\end{align}
The two mixed products coincide by cyclicity of the trace and because
the powers of $\mtx H$ commute. Thus the three terms correspond to
signal--signal, signal--noise, and noise--noise contributions.

Replacing the weighted traces by the reductions proved above gives the
following candidate for the leading covariance and its transformation to
variance-component coordinates:
\begin{align}
\widetilde{\mtx U}_n
&=\frac{2}{q^2}\Bigg\{
F_b\begin{pmatrix}g_2&g_3\\g_3&g_4\end{pmatrix}
+2F_{be}\begin{pmatrix}1&g_2\\g_2&g_3\end{pmatrix}
+F_e\begin{pmatrix}1&1\\1&g_2\end{pmatrix}
+\kappa\begin{pmatrix}1&1\\1&1\end{pmatrix}
\Bigg\},
\label{eq:re_U_decomposition}\\
\mtx V_{\mathrm{re}}
&=\widehat{\mtx A}_n\widetilde{\mtx U}_n
  \widehat{\mtx A}_n^\top.
\label{eq:re_V_matrix_transform}
\end{align}
Substituting \eqref{eq:weighted_gk_reduction}--%
\eqref{eq:aggregate_trace_reduction_1} into
\eqref{eq:re_general_conditional_covariance} gives this leading covariance,
apart from the factors $a_n^2$ and $a_n$ in the first two terms.
Replacing these factors by one contributes
$O_P\{\ell_n/(nq)\}=o_P(q^{-1})$, since
$F_b+F_{be}=O(\ell_nq)$ and $a_n-1=O_P(n^{-1})$.
Keeping the three sources of error separate gives
\begin{equation}
\label{eq:re_cov_relative_approx}
q\|\mtx U_n-\widetilde{\mtx U}_n\|
=O_P\!\left(\frac{\ell_n}{n}+\frac{\ell_n}{\sqrt n}
       +\ell_n^2\sqrt{\frac{\log n}{n}}\right)=o_P(1).
\end{equation}
The first term comes from $a_n-1$, the second from the mixed weighted
traces, and the last bounds the noise-trace errors. The final equality
uses $p\asymp n$, $\ell_n\ge1$, and the second oracle growth condition.
The approximation includes both the weighted-trace errors and the
population-to-sample energy difference. We next show that its
$o_P(q^{-1})$ error is small relative to every directional variance.

The first two moment matrices in \eqref{eq:re_U_decomposition} are
positive semidefinite: if $\lambda_i$ are the eigenvalues of $\mtx H$,
their quadratic forms are averages of
$\lambda_i^2(u+v\lambda_i)^2$ and
$\lambda_i(u+v\lambda_i)^2$, respectively.
The noise matrix has determinant $g_2-1\ge c$ and bounded trace.
Since $F_e\ge\tr^2(\bar{\mtx\Sigma}_e)/q\ge cq$, this term supplies the
lower covariance bound. The four functionals are $O(\ell_nq)$, so
\begin{equation}
\label{eq:re_U_scale}
\frac cq\mtx I_2\preceq\widetilde{\mtx U}_n
\preceq\frac{C\ell_n}{q}\mtx I_2.
\end{equation}
In particular, the covariance comparison is uniform over directions:
\[
\sup_{\vct b\ne\vct0}
\left|\frac{\vct b^\top\mtx U_n\vct b}
{\vct b^\top\widetilde{\mtx U}_n\vct b}-1\right|
\le Cq\|\mtx U_n-\widetilde{\mtx U}_n\|=o_P(1).
\]
Thus $\mtx U_n$ obeys the same scale bounds, with adjusted constants,
on an event of probability tending to one. The bounded norms of
$\widehat{\mtx A}_n$ and its inverse then imply
\begin{equation}
\label{eq:re_variance_scale}
\frac cq\mtx I_2\preceq\mtx V_{\mathrm{re}}
\preceq\frac{C\ell_n}{q}\mtx I_2
\end{equation}
with probability tending to one.

The transformed covariance gives the design-dependent oracle variance
used in the theorem:
\begin{equation}
\label{eq:sigma_r_mt}
\vct d^\top\mtx V_{\mathrm{re}}\vct d
=\mathcal V_{\mathrm{re}}(\mtx\Sigma_b,\bar{\mtx\Sigma}_e;
\rho^2,\sigma^2;\kappa_{\mathrm{tot}},g_2,g_3,g_4)
=\sigma_{r,\mathrm{re}}^2.
\end{equation}
To see this, put $\tau=\rho^2+\sigma^2$ and use
\begin{equation}
\label{eq:random_transformed_gradient}
\widehat{\mtx A}_n^\top\vct d
=\frac1{(g_2-1)\tau^2}
\begin{pmatrix}-(g_2\rho^2+\sigma^2)\\\tau\end{pmatrix}
\end{equation}
in \eqref{eq:re_U_decomposition}.
The identity $\widehat{\mtx A}_n(1,1)^\top=(1,0)^\top$ shows that the
extra dispersion enters the noise-component coordinate after the moment
transformation. Its contribution in the SNR direction is therefore
$2\kappa_{\mathrm{tot}}\rho^4/(q^2\tau^4)$. The three homogeneous terms use precisely the scalar contractions
calculated after \eqref{eq:fixed_homo_transformed_gradient}, with
$\mtx\Sigma_e$ replaced by $\bar{\mtx\Sigma}_e$. Adding the dispersion
contribution gives \eqref{eq:random_variance_function}.
Since $c\le\|\vct d\|_2\le C$, \eqref{eq:re_variance_scale} implies
$c/q\le\sigma_{r,\mathrm{re}}^2\le C\ell_n/q$.

\proofblock{Conditional normal approximation.}
Conditional Gaussianity now gives the same quadratic-form representation
as in S3, with the heterogeneous response covariance. Conditional on
$\mtx X$, write
$\vecop(\mtx Y)=\mtx\Omega_X^{1/2}\vct z$ with
$\vct z\sim\mathcal N(\vct0,\mtx I_{nq})$.
The noise part of \eqref{eq:re_heterogeneous_omega} is
permutation-similar to $\diag(\mtx\Sigma_1,\ldots,\mtx\Sigma_n)$.
Thus $\|\mtx\Omega_X\|\le C\ell_n$ on the design event.
For $\vct b=(b_0,b_1)^\top\ne\vct0$,
\[
\vct b^\top\vct w
=\vct z^\top\mtx Q_{\vct b}\vct z-\tr(\mtx Q_{\vct b}),\qquad
\mtx Q_{\vct b}
=\frac1{\sqrt n\,q}\mtx\Omega_X^{1/2}
 \{\mtx I_q\otimes(b_0\mtx I_n+b_1\mtx H)\}
 \mtx\Omega_X^{1/2}.
\]
Its operator norm is at most
$C\ell_n\|\vct b\|_2/(\sqrt n\,q)$, whereas
\[
2\|\mtx Q_{\vct b}\|_F^2
=\vct b^\top\mtx U_n\vct b
\ge c\|\vct b\|_2^2/q.
\]
The conditional version of \eqref{eq:gaussian_quadratic_be} therefore yields
\begin{equation}
\label{eq:re_uniform_quadratic_clt}
\sup_{\vct b\ne\vct0}\sup_t
\left|\PP\!\left(
\frac{\vct b^\top\vct w}{(\vct b^\top\mtx U_n\vct b)^{1/2}}
\le t\,\middle|\,\mtx X\right)-\Phi(t)\right|
=O_P\!\left(\frac{\ell_n}{\sqrt{nq}}\right)=o_P(1).
\end{equation}
The first oracle growth condition implies the last equality. Uniformity
is needed because the SNR direction in moment coordinates is
$\widehat{\mtx A}_n^\top\vct d$, which depends on the observed design.

\proofblock{Removing the conditioning and forming the ratio.}
Equation~\eqref{eq:re_direct_estimators} gives
\begin{equation}
\label{eq:re_theta_w_exact}
\sqrt n\{\widehat{\vct\theta}
-\E(\widehat{\vct\theta}\mid\mtx X)\}
=\widehat{\mtx A}_n\vct w.
\end{equation}
Take the $\mtx X$-measurable direction
$\vct b=\widehat{\mtx A}_n^\top\vct d$ in
\eqref{eq:re_uniform_quadratic_clt}. Its norm is bounded above and away from
zero on the design event, and
\eqref{eq:re_cov_relative_approx}--\eqref{eq:sigma_r_mt} give
\[
\vct b^\top\mtx U_n\vct b
=\vct d^\top\mtx V_{\mathrm{re}}\vct d+o_P(q^{-1})
=\sigma_{r,\mathrm{re}}^2\{1+o_P(1)\}.
\]
To make the removal of conditioning explicit, put
\[
Z_n=\frac{\vct b^\top\vct w}
{(\vct b^\top\mtx U_n\vct b)^{1/2}},\qquad
\epsilon_n(\mtx X)
=\sup_t|\PP(Z_n\le t\mid\mtx X)-\Phi(t)|.
\]
Assign $Z_n=0$ when its denominator vanishes; the covariance lower
bound makes this event negligible. The conditional bound gives
$\epsilon_n(\mtx X)\to_P0$, and
$0\le\epsilon_n(\mtx X)\le1$. For every $\epsilon>0$,
\[
\sup_t|\PP(Z_n\le t)-\Phi(t)|
\le\E\epsilon_n(\mtx X)
\le\epsilon+\PP\{\epsilon_n(\mtx X)>\epsilon\}\longrightarrow\epsilon.
\]
Letting $\epsilon\downarrow0$ proves $Z_n\Longrightarrow\mathcal N(0,1)$.
By \eqref{eq:re_theta_w_exact},
\[
\frac{\sqrt n\,\vct d^\top\{\widehat{\vct\theta}
-\E(\widehat{\vct\theta}\mid\mtx X)\}}{\sigma_{r,\mathrm{re}}}
=Z_n\left\{\frac{\vct b^\top\mtx U_n\vct b}
{\sigma_{r,\mathrm{re}}^2}\right\}^{1/2}.
\]
The second factor tends to one in probability, so Slutsky's theorem
applies without independence between it and $Z_n$. It remains to replace
conditional centering by the population target. On the same standard-deviation scale, the bias is
\[
\frac{\sqrt n\,|\vct d^\top
\{\E(\widehat{\vct\theta}\mid\mtx X)-\vct\theta\}|}
{\sigma_{r,\mathrm{re}}}
=O_P\!\left(\ell_n\sqrt{\frac qn}\right)=o_P(1),
\]
by \eqref{eq:conditional_mean_consistency},
$\sigma_{r,\mathrm{re}}^2\ge c/q$, and the first oracle growth condition.
Thus, with $\vct\Delta=\widehat{\vct\theta}-\vct\theta$,
\begin{equation}
\label{eq:re_linear_snr_clt}
\frac{\sqrt n\,\vct d^\top\vct\Delta}{\sigma_{r,\mathrm{re}}}
\Longrightarrow\mathcal N(0,1).
\end{equation}

To transfer this linear limit to the ratio, we need only consistency of
its denominator. Here $\hat\tau=\tr(\mtx M_0)/q$, and its conditional
moments satisfy
\[
\E(\hat\tau\mid\mtx X)=\sigma^2+a_n\rho^2,\qquad
\Var(\hat\tau\mid\mtx X)=\frac1n(\mtx U_n)_{11}
=O_P\!\left(\frac{\ell_n}{nq}\right).
\]
A conditional variance bound in probability is sufficient here; we do
not need an unconditional expectation bound for that random variance.
Indeed, for $v_n^2=\ell_n/(nq)$ and fixed $R,M>0$, conditional
Chebyshev gives
\begin{equation}
\label{eq:conditional_chebyshev_localization}
\PP\{|\hat\tau-\E(\hat\tau\mid\mtx X)|>Mv_n\}
\le\PP\{\Var(\hat\tau\mid\mtx X)>Rv_n^2\}+\frac{R}{M^2}.
\end{equation}
First choose $R$ using the $O_P(v_n^2)$ variance bound, and then choose
$M$. This proves a fluctuation of order $O_P(v_n)$. Together with
$a_n-1=O_P(n^{-1})$, it gives
\begin{equation}
\label{eq:re_total_variance_consistency}
\hat\tau-\tau
=O_P\!\left(\sqrt{\frac{\ell_n}{nq}}+\frac1n\right)=o_P(1),
\qquad \frac{\hat\tau}{\tau}\to_P1.
\end{equation}
The first oracle growth condition implies $\ell_n/(nq)\to0$.
The oracle part of Lemma~\ref{lem:common_studentization}, using
\eqref{eq:re_variance_scale}, \eqref{eq:re_linear_snr_clt} and the exact
identity \eqref{eq:exact_ratio_identity}, therefore gives
\begin{equation}
\label{eq:re_oracle_snr_clt}
\frac{\sqrt n(\hat r^2-r^2)}{\sigma_{r,\mathrm{re}}}
\Longrightarrow\mathcal N(0,1).
\end{equation}
This proves Theorem~\ref{thm:asymptotics_mt}. The normalization remains
design-dependent throughout; no deterministic limit for the spectral
moments is required.
\end{proof}

\subsection{Feasible variance estimation and studentization}
\label{sec:proof_studentized_random_main}

The oracle proof leaves four response-covariance functionals to estimate:
$F_b,F_e,F_{be}$ and $\kappa_{\mathrm{tot}}$. The first three are controlled
by the covariance estimates under general heterogeneity. Scalar
heterogeneity is used only to estimate the fourth, through an additional
response moment.

Put $\varepsilon_n=\sqrt{(q+4\log p)/n}$.
Let $\widehat{\widetilde{\mtx U}}_n$ be obtained from
\eqref{eq:re_U_decomposition} by replacing
$(\mtx\Sigma_b,\bar{\mtx\Sigma}_e,\kappa_{\mathrm{tot}})$ with
$(\widehat{\mtx\Sigma}_b,\widehat{\bar{\mtx\Sigma}}_e,
\hat\kappa_{\mathrm{tot}})$, and define
\begin{equation}
\label{eq:re_plugin_matrix}
\widehat{\mtx V}_{\mathrm{re}}
=\widehat{\mtx A}_n\widehat{\widetilde{\mtx U}}_n
 \widehat{\mtx A}_n^\top.
\end{equation}
The same algebra as in \eqref{eq:sigma_r_mt} gives the exact identity
\begin{equation}
\label{eq:hat_sigma_r_mt}
\hat{\vct d}^{\top}\widehat{\mtx V}_{\mathrm{re}}\hat{\vct d}
=\hat\sigma_{r,\mathrm{re}}^2
\end{equation}
for the plug-in variance in \eqref{eq:random_plugin_variance_main}.

The following lemma separates these tasks. Part~(i) controls the covariance
matrices under general heterogeneity; part~(ii) estimates the dispersion
under scalar heterogeneity; part~(iii) transfers the functional errors to
the variance estimate.

\begin{lemma}[Random-effects variance estimation]
\label{lem:studentization_random}
Under the conditions of Theorem~\ref{thm:asymptotics_mt}, the following hold.
\begin{enumerate}[(i)]
\item The covariance estimators satisfy
\begin{equation}
\label{eq:re_covariance_estimation_rate}
\|\widehat{\mtx\Sigma}_b-\mtx\Sigma_b\|
+\|\widehat{\bar{\mtx\Sigma}}_e-\bar{\mtx\Sigma}_e\|
=O_P(\ell_n\varepsilon_n)=o_P(1).
\end{equation}
\item Under the scalar heterogeneity model and the additional conditions of
Theorem~\ref{thm:studentized_random_mt},
\[
\hat\eta-\eta=O_P(\ell_n^2\varepsilon_n)=o_P(1),\qquad
q^{-1}|\hat\kappa_{\mathrm{tot}}-\kappa_{\mathrm{tot}}|=o_P(1).
\]
\item More generally, if $\ell_n^2\varepsilon_n\to0$ and
$q^{-1}|\hat\kappa_{\mathrm{tot}}-\kappa_{\mathrm{tot}}|=o_P(1)$, then
\begin{equation}
\label{eq:re_plugin_consistency}
\|\widehat{\mtx V}_{\mathrm{re}}-\mtx V_{\mathrm{re}}\|
=o_P(q^{-1}).
\end{equation}
\end{enumerate}
\end{lemma}

\begin{proof}
\proofblock{Covariance estimation: part (i).}
We separate the fluctuation about the conditional mean from the difference
between that mean and the population moment. The fluctuation of either
matrix moment is controlled by the same conditional Gaussian
quadratic-form argument. Set
$\mtx C_j=\mtx H^j$, $j=0,1$.
For a unit vector $\vct u\in\mathbb R^q$,
\[
\vct u^\top\{\mtx M_j-\E(\mtx M_j\mid\mtx X)\}\vct u
=\vct z^\top\mtx Q_{j,\vct u}\vct z-\tr(\mtx Q_{j,\vct u}),
\qquad
\mtx Q_{j,\vct u}
=\frac1n\mtx\Omega_X^{1/2}
(\vct u\vct u^\top\otimes\mtx C_j)\mtx\Omega_X^{1/2}.
\]
On the design event,
$\|\mtx\Omega_X\|\le C\ell_n$ and $\|\mtx C_j\|\le C$.
Moreover, $\mtx Q_{j,\vct u}$ has rank at most $n$, so
\[
\|\mtx Q_{j,\vct u}\|\le C\ell_n/n,\qquad
\|\mtx Q_{j,\vct u}\|_F\le C\ell_n/\sqrt n.
\]
Hanson--Wright and a $1/4$-net of the $q$-dimensional unit sphere
therefore give, with conditional probability at least $1-Cp^{-3}$,
\begin{equation}
\label{eq:re_direct_matrix_concentration}
\max_{j=0,1}\|\mtx M_j-\E(\mtx M_j\mid\mtx X)\|
\le C\ell_n\left\{\sqrt{\frac{q+4\log p}{n}}
+\frac{q+4\log p}{n}\right\}.
\end{equation}
For a fixed direction, Hanson--Wright gives failure probability
$2e^{-ct}$ at threshold
$C\ell_n\{\sqrt{t/n}+t/n\}$. The two moments and a net of at most
$9^q$ points therefore have total failure probability at most
$4\cdot9^q e^{-ct}$. Taking $t=A(q+4\log p)$ for a sufficiently
large fixed $A$, and applying
Lemma~\ref{lem:operator norm bound by net cover}, gives the stated
$Cp^{-3}$ bound. The oracle assumptions imply
\[
\ell_n\varepsilon_n
\le\ell_n\sqrt{q/n}+2\ell_n\sqrt{\log p/n}\to0,
\]
since $\ell_n\ge1$ and $\ell_n^2\sqrt{\log p/n}\to0$.
In particular $\varepsilon_n\to0$, so the linear term
$\ell_n\varepsilon_n^2$ in \eqref{eq:re_direct_matrix_concentration}
is absorbed by $\ell_n\varepsilon_n$.

For the difference in means, the oracle proof has already identified the
two terms $\mtx R_X$ and $(a_n-1)\mtx\Sigma_b$. Set
$d_X=\max_i|t_i-1|$. Since each noise covariance is positive semidefinite,
\[
-d_X\bar{\mtx\Sigma}_e\preceq\mtx R_X
\preceq d_X\bar{\mtx\Sigma}_e,
\qquad
\|\mtx R_X\|\le d_X\|\bar{\mtx\Sigma}_e\|
=O_P(\ell_n\varepsilon_n).
\]
The matrix inequality, rather than a sum of individual operator norms,
is what uses the average noise covariance. Also
$|a_n-1|\|\mtx\Sigma_b\|=O_P(\ell_n/n)$.
Writing $\mtx J_j=\mtx M_j-\E(\mtx M_j\mid\mtx X)$, the exact
estimator identities give
\begin{align*}
\widehat{\mtx\Sigma}_b-\mtx\Sigma_b
&=\frac{\mtx J_1-\mtx J_0+\mtx R_X}{g_2-1}
 +(a_n-1)\mtx\Sigma_b,\\
\widehat{\bar{\mtx\Sigma}}_e-\bar{\mtx\Sigma}_e
&=\frac{g_2\mtx J_0-\mtx J_1-\mtx R_X}{g_2-1}.
\end{align*}
Combining these formulas with the design bounds and
\eqref{eq:re_direct_matrix_concentration} proves part~(i).

Both remaining parts use the induced errors in the three covariance
functionals. Write
$\mtx\Delta_b=\widehat{\mtx\Sigma}_b-\mtx\Sigma_b$ and
$\mtx\Delta_e=\widehat{\bar{\mtx\Sigma}}_e-\bar{\mtx\Sigma}_e$.
For this calculation, $\mtx\Sigma_a$ denotes either $\mtx\Sigma_b$ or
$\bar{\mtx\Sigma}_e$, with its corresponding error $\mtx\Delta_a$. The
bound
\[
\frac1q\big|\|\mtx\Sigma_a+\mtx\Delta_a\|_F^2
               -\|\mtx\Sigma_a\|_F^2\big|
\le2\|\mtx\Sigma_a\|\|\mtx\Delta_a\|+\|\mtx\Delta_a\|^2
\]
and its mixed-product analogue yield
\begin{equation}
\label{eq:re_normalized_covariance_functionals}
\frac1q\left(
|\widehat F_b-F_b|+|\widehat F_e-F_e|
+|\widehat F_{be}-F_{be}|\right)
=O_P(\ell_n^2\varepsilon_n).
\end{equation}
Hats denote substitution of the covariance estimates into the corresponding
functionals. Part~(i) also gives
$|\hat\rho^2-\rho^2|+|\hat\sigma^2-\sigma^2|
=O_P(\ell_n\varepsilon_n)$.

\proofblock{Estimating scalar heterogeneity: part (ii).}
We first identify the fourth-moment equation, then control its sampling
fluctuation and the effect of plugging in the covariance estimates.
Under \eqref{eq:scalar_heterogeneity_mt},
$\bar{\mtx\Sigma}_e=\mtx\Sigma_e$,
$\kappa_{\mathrm{tot}}=\eta F_e$, and $\eta=O(1)$ because
$\max_i\nu_i\le C$.
Define
\[
\hat m_4=\frac1n\sum_i(\vct y_i^\top\vct y_i)^2,\qquad
b_{2n}=\frac1n\sum_i t_i^2,\qquad
c_{\nu n}=\frac1n\sum_i t_i\nu_i,
\]
and
\[
L_b=2F_b+q^2\rho^4,\qquad
L_e=2F_e+q^2\sigma^4,\qquad
L_{be}=2F_{be}+q^2\sigma^2\rho^2.
\]
The quantities $L_b,L_e,L_{be}$ collect the signal, noise and mixed terms
in the Gaussian fourth-moment expansion. In this notation, the main-paper
estimator \eqref{eq:eta_hat_scalar_main} is
\begin{equation}
\label{eq:eta_hat_scalar_mt}
\hat\eta
=\frac{\hat m_4-b_{2n}\widehat L_b-2\widehat L_{be}}
       {\widehat L_e}-1.
\end{equation}
The design bounds give $b_{2n}=O_P(1)$, and
\eqref{equation: leave one out eq 2} gives
$c_{\nu n}=1+O_P(n^{-1})$.

For each row,
$\Cov(\vct y_i\mid\mtx X)=a_nt_i\mtx\Sigma_b+\nu_i\mtx\Sigma_e$.
For a centered Gaussian vector with covariance $\mtx C$,
$\E\|\vct y\|_2^4=2\tr(\mtx C^2)+\tr^2(\mtx C)$.
Inserting the row covariance gives the individual-row identity
\[
\E(\|\vct y_i\|_2^4\mid\mtx X)
=a_n^2t_i^2L_b+\nu_i^2L_e+2a_nt_i\nu_iL_{be}.
\]
Averaging this identity does not require independence between rows, and gives
\begin{equation}
\label{eq:eta_conditional_fourth_moment}
\E(\hat m_4\mid\mtx X)
=(1+\eta)L_e+a_n^2b_{2n}L_b+2a_nc_{\nu n}L_{be}.
\end{equation}
The noise coefficient is $1+\eta=n^{-1}\sum_i\nu_i^2$; the signal and
mixed coefficients retain the realized design weights. All three
$L$-functionals are $O(q^2)$, so replacing $a_n$ and $c_{\nu n}$ by one
creates an $O_P(q^2/n)$ error. This yields the population equation matched
by the estimator, up to a negligible design error.

To control the sampling fluctuation, we cannot treat the fourth moments
as independent across rows: the rows still share the random coefficients
conditional on $\mtx X$. Instead, apply
\eqref{eq:gaussian_poincare_covariance} to the full response vector,
conditionally on $\mtx X$. Its covariance has operator norm $O_P(\ell_n)$, and
\[
\nabla_{\vct y_i}\hat m_4=\frac4n\|\vct y_i\|_2^2\vct y_i,\qquad
\|\nabla\hat m_4\|_2^2
=\frac{16}{n^2}\sum_i\|\vct y_i\|_2^6.
\]
The row covariance norms are uniformly $O_P(\ell_n)$.
For a standard Gaussian $\vct g\in\mathbb R^q$,
$\E\|\vct g\|_2^6=q(q+2)(q+4)$, so
$\max_i\E(\|\vct y_i\|_2^6\mid\mtx X)=O_P(\ell_n^3q^3)$.
Consequently,
\begin{equation}
\label{eq:eta_poincare_variance}
\Var(\hat m_4\mid\mtx X)
\le\|\mtx\Omega_X\|\,
\E(\|\nabla\hat m_4\|_2^2\mid\mtx X)
=O_P(\ell_n^4q^3/n).
\end{equation}
Set $v_n=\ell_n^2/\sqrt{nq}$. The conditional variance of
$q^{-2}\hat m_4$ is $O_P(v_n^2)$ by
\eqref{eq:eta_poincare_variance}. The same localization as in
\eqref{eq:conditional_chebyshev_localization} gives, for $R,M>0$,
\[
\PP\{q^{-2}|\hat m_4-\E(\hat m_4\mid\mtx X)|>Mv_n\}
\le\PP\{q^{-4}\Var(\hat m_4\mid\mtx X)>Rv_n^2\}
 +\frac{R}{M^2}.
\]
Choosing $R$ and then $M$ proves
\begin{equation}
\label{eq:eta_empirical_fourth_concentration}
q^{-2}\{\hat m_4-\E(\hat m_4\mid\mtx X)\}
=O_P\!\left(\frac{\ell_n^2}{\sqrt{nq}}\right).
\end{equation}
Combining this bound with \eqref{eq:eta_conditional_fourth_moment} gives
\begin{equation}
\label{eq:eta_fourth_moment_reduced}
\hat m_4
=(1+\eta)L_e+b_{2n}L_b+2L_{be}
+O_P\!\left[q^2\left\{\frac{\ell_n^2}{\sqrt{nq}}+\frac1n\right\}\right].
\end{equation}

We finally replace the population fourth-moment coefficients by their
estimates. Part~(i) gives $\hat\rho^2,\hat\sigma^2=O_P(1)$ and, for
example,
\[
|\hat\rho^4-\rho^4|
=|\hat\rho^2-\rho^2|\,|\hat\rho^2+\rho^2|
=O_P(\ell_n\varepsilon_n).
\]
The corresponding identity for $\hat\sigma^4$ and the two-term expansion
of $\hat\rho^2\hat\sigma^2-\rho^2\sigma^2$, together with
\eqref{eq:re_normalized_covariance_functionals}, imply
\[
q^{-2}\left(
|\widehat L_b-L_b|+|\widehat L_e-L_e|
+|\widehat L_{be}-L_{be}|\right)
=O_P(\ell_n^2\varepsilon_n).
\]
Since $L_e\ge q^2\sigma^4\ge cq^2$ and
$\ell_n^2\varepsilon_n\to0$, the estimated denominator satisfies
$\widehat L_e\ge cq^2/2$ with probability tending to one. On the event
that it vanishes, $\hat\eta$ may be assigned any fixed value without
affecting the limit. Subtracting $\eta$ from
\eqref{eq:eta_hat_scalar_mt} gives the exact error decomposition
\begin{align*}
\widehat L_e(\hat\eta-\eta)
&=\hat m_4-\{(1+\eta)L_e+b_{2n}L_b+2L_{be}\}\\
&\quad-b_{2n}(\widehat L_b-L_b)
       -2(\widehat L_{be}-L_{be})
       -(1+\eta)(\widehat L_e-L_e).
\end{align*}
Divide by $\widehat L_e$, use $b_{2n}=O_P(1)$ and $\eta=O(1)$, and
apply \eqref{eq:eta_fourth_moment_reduced}. Since
$\ell_n^2/\sqrt{nq}+n^{-1}=O(\ell_n^2\varepsilon_n)$, this proves
\[
\hat\eta-\eta=O_P(\ell_n^2\varepsilon_n)=o_P(1).
\]
Truncation preserves this rate: for $\hat\eta_+=\max(\hat\eta,0)$,
$|\hat\eta_+-\eta|\le|\hat\eta-\eta|$ and $\hat\eta_+=O_P(1)$.
To estimate the variance, however, we need the product
$\kappa_{\mathrm{tot}}=\eta F_e$ at scale $o_P(q)$, not only consistency
of $\hat\eta$. Since $F_e/q=O(\ell_n)$,
\begin{align}
\frac1q|\hat\kappa_{\mathrm{tot}}-\kappa_{\mathrm{tot}}|
&\le|\hat\eta_+-\eta|\frac{F_e}{q}
+\hat\eta_+\frac{|\widehat F_e-F_e|}{q}\notag\\
&=O_P(\ell_n^3\varepsilon_n)+O_P(\ell_n^2\varepsilon_n)
=o_P(1).
\label{eq:kappa_relative_consistency}
\end{align}
This is where the additional factor $\ell_n$ enters, explaining the
feasible growth condition $\ell_n^3\varepsilon_n\to0$ and completing
part~(ii).

\proofblock{Plug-in variance: part (iii).}
Under the assumptions of part~(iii),
\eqref{eq:re_normalized_covariance_functionals} controls the first three
functional errors in \eqref{eq:re_U_decomposition}, and consistency of
$\hat\kappa_{\mathrm{tot}}$ controls the fourth.
All four design matrices have bounded operator norm on
\eqref{eq:common_design_event}, so
\[
\|\widehat{\widetilde{\mtx U}}_n-\widetilde{\mtx U}_n\|
=o_P(q^{-1}).
\]
The exact identity
\[
\widehat{\mtx V}_{\mathrm{re}}-\mtx V_{\mathrm{re}}
=\widehat{\mtx A}_n
(\widehat{\widetilde{\mtx U}}_n-\widetilde{\mtx U}_n)
\widehat{\mtx A}_n^\top
\]
and $\|\widehat{\mtx A}_n\|=O_P(1)$ prove part~(iii).
\end{proof}

\begin{samepage}
\begin{proof}[Proof of Theorem~\ref{thm:studentized_random_mt}]
The oracle result already supplies the scalar normal approximation and
total-variance consistency. Under scalar heterogeneity, parts~(ii) and~(iii)
of Lemma~\ref{lem:studentization_random} add
$\|\widehat{\mtx V}_{\mathrm{re}}-\mtx V_{\mathrm{re}}\|=o_P(q^{-1})$.
The covariance scale \eqref{eq:re_variance_scale}, the linear limit
\eqref{eq:re_linear_snr_clt}, and total-variance consistency
\eqref{eq:re_total_variance_consistency} therefore verify all remaining
conditions of Lemma~\ref{lem:common_studentization}; the oracle growth
conditions already imply $\ell_n/(nq)\to0$. Using the scalar identities
\eqref{eq:sigma_r_mt} and \eqref{eq:hat_sigma_r_mt}, the lemma gives
\[
\frac{\hat\sigma_{r,\mathrm{re}}^2}{\sigma_{r,\mathrm{re}}^2}\to_P1,
\qquad
\PP(\hat\sigma_{r,\mathrm{re}}^2>0)\to1,\qquad
\frac{\sqrt n(\hat r^2-r^2)}{\hat\sigma_{r,\mathrm{re}}}
\Longrightarrow\mathcal N(0,1).
\]
This proves Theorem~\ref{thm:studentized_random_mt} with the standard error
defined in the main paper.
\end{proof}
\end{samepage}

\section{Implementation details for the yeast analysis}
\label{sec:yeast_supp}

We average all finite replicate values in the released \texttt{pheno\_raw}
object and align the resulting strain means with the rows of
\texttt{gdata} \citep{bloom2015data}. Requiring a mean in all 20 environments
leaves 3474 strains and excludes 916; no phenotype is imputed. Lactate and
Lactose have the most missing strain means, 628 and 574, respectively;
these counts can overlap. Among retained strain--environment cells,
14,535 have one valid measurement, 51,693 have two, 605 have three, and
2647 have four. Only YPD has more than two. The missingness mechanism is
not identified by these counts. The released measurements already follow
the source protocol's colony filtering, within-plate spatial adjustment,
and adjustment for control-medium growth
\citep{bloom2013finding,bloom2015interactions}; we do not repeat those steps.

Let $n=3474$, $m=n-1=3473$, and
$\mtx H=\mtx I_n-\vct1_n\vct1_n^\top/n$.
We center the response means, $\mtx Y_c=\mtx H\mtx Y$, without rescaling
their columns. Each marker is centered and divided by
$\{m^{-1}\sum_i(x_{ij}-\bar x_j)^2\}^{1/2}$, giving a final matrix
$\mtx X_c$ with $\|\mtx X_c\|_F^2=mp$. All 28,220 markers have
values $-1$ or $1$, no missing genotypes, and nonzero variance in the
retained sample. The random-coefficient model uses this final marker scale.
Writing $\mtx H=\mtx Q\mtx Q^\top$ with
$\mtx Q^\top\mtx Q=\mtx I_m$, the fixed-projection argument in
Section~\ref{sec:supp_fixed_design_theorem} applies Section~3 to
$(\mtx Q^\top\mtx X_c,\mtx Q^\top\mtx Y)$ under its stated model
conditions. Here $m$ is the dimension after removing the intercept.
We compute the equivalent centered cross-products without constructing
$\mtx Q$:
\[
\mtx K=\frac{\mtx X_c\mtx X_c^\top}{p},\qquad
\mtx M_0=\frac{\mtx Y_c^\top\mtx Y_c}{m},\qquad
\mtx M_1=\frac{\mtx Y_c^\top\mtx K\mtx Y_c}{m},\qquad
g_k=\frac{\tr(\mtx K^k)}m.
\]
Thus $\tr(\mtx K)=m$, and the covariance estimates are
\[
\widehat{\mtx\Sigma}_b=\frac{\mtx M_1-\mtx M_0}{g_2-1},\qquad
\widehat{\mtx\Sigma}_e=\mtx M_0-\widehat{\mtx\Sigma}_b.
\]
The earlier stored kernel is $(p/m)\mtx K$; the analysis script accounts
for this factor. For Table~\ref{tab:yeast_environment_contributions}, let
$\hat s_j=(\widehat{\mtx\Sigma}_b)_{jj}$ and
$\hat v_j=(\mtx M_0)_{jj}$. The table gives $\hat v_j$,
$100\hat v_j/\sum_k\hat v_k$, and $\hat s_j/\hat v_j$.
Using the same sample and design for each environment gives
\[
\hat r^2=\frac{\sum_j\hat s_j}{\sum_j\hat v_j}
=\sum_j\frac{\hat v_j}{\sum_k\hat v_k}\frac{\hat s_j}{\hat v_j},
\]
and each $\hat s_j/\hat v_j$ equals its separate single-environment
moment estimate.

For the Section~4 illustration, we apply
\eqref{eq:eta_hat_scalar_main} and \eqref{eq:random_plugin_variance_main}
using the $n$ original centered strain rows and
$\mtx K_n=n\mtx X_c\mtx X_c^\top/\|\mtx X_c\|_F^2$,
with all second, fourth, and design moments normalized by $n$.
The fourth moment uses strain rows, not rotated contrast coordinates.
This $n$-based convention accounts for the small point-estimate difference
from the primary $m$-based fit. It is a working-model calculation after
empirical centering and scaling; the homogeneous-noise fixed-projection
argument above does not justify heterogeneous-noise inference.

\clearpage

\fi

\end{document}